\documentclass[12pt]{article}
\usepackage{amssymb,amsmath,amsthm,tikz,multirow,nccrules,float,pst-solides3d,enumerate,graphicx,subfig}
\usetikzlibrary{arrows,calc}

\title{Tilings of the sphere by congruent quadrilaterals I: edge combination $a^2bc$}
\author{Yixi Liao, Pinren Qian, Erxiao Wang\thanks{Corresponding author (wang.eric@zjnu.edu.cn). Research was supported by Key projects of Zhejiang Natural Science Foundation No. LZ22A010003 and ZJNU Shuang-Long Distinguished Professorship Fund No. YS304319159.}, Yingyun Xu \\
	Zhejiang Normal University}

\allowdisplaybreaks[4]  

\newcommand{\sub}{\subset}
\newcommand{\pa}{\partial}
\newcommand{\mb}{\mathbf}
\newcommand{\mc}{\mathcal}
\newcommand{\bb}{\mathbb}
\newcommand\aaa{\alpha}
\newcommand\bbb{\beta}
\newcommand\ccc{\gamma}
\newcommand\ddd{\delta}

\newcommand{\dash}{\hspace{0.1em}\dashrule{0.7}{2.4 1 2.4 1 2.4}\hspace{0.1em}} 
\newcommand{\thin}{\hspace{0.1em}\rule{0.7pt}{0.8em}\hspace{0.1em}}
\newcommand{\thick}{\hspace{0.1em}\rule{1.5pt}{0.8em}\hspace{0.1em}}

\newtheorem{theorem}{Theorem}
\newtheorem{lemma}[theorem]{Lemma}

\newtheorem{proposition}[theorem]{Proposition}
\newtheorem*{theorem*}{Theorem}

\theoremstyle{definition}
\newtheorem*{definition*}{Definition}
\newtheorem*{case*}{Case}
\newtheorem*{subcase*}{Subcase}

\theoremstyle{remark}

\numberwithin{equation}{section}

\begin{document}
	
	\date{}
	\maketitle
	
	\begin{abstract}
		Edge-to-edge tilings of the sphere by congruent $a^2bc$-quadrilaterals are classified as $3$ classes: a sequence of $2$-parameter families of $2$-layer earth map tilings with $2n$ $(n\ge3)$ tiles, a $1$-parameter family of quadrilateral subdivisions of the octahedron with $24$ tiles together with a flip modification for a special parameter, and a sequence of $3$-layer earth map tilings with $8n$ $(n\ge2)$ tiles together with two flip modifications for odd $n$. We also describe the moduli and calculate the geometric data. 
		
	\end{abstract}

\section{Introduction} \label{introduction}
Tiling problems have existed for thousands of years, and their modern studies by scientists have lasted for hundreds of years. However, a full classification of monohedral convex tilings of the plane has been completed only recently, see \cite{rao} for the hardest pentagon case and see \cite{zong} for a  recent  survey. There are not as many studies on spherical tilings as the planar ones. Recall that for edge-to-edge tilings of the sphere by congruent simple polygons with all vertices' degrees being $\ge 3$, the tile must be triangle, quadrilateral, or pentagon (see \cite{ua2}, for example). The study of triangle case was started by Sommerville \cite{so} in 1924, initially classified by Davies \cite{da} in 1967 and completed with full details by Ueno and Agaoka \cite{ua} in 2002. Recent works of Wang, Yan and Akama \cite{wy1,wy2,wy3,awy} studied pentagon case. However the quadrilateral case has remained largely open after some early explorations \cite{ua2,sa,ac}, and we will give its full classification in this series of $3$ papers (see also  \cite{lw,lpwx}). There is also a simultaneous independent study by a Hong Kong group Cheung, Luk and Yan via quite different strategies, and it is great that the readers can compare and combine two outputs for a more complete view of such a complicated classification program. 

\begin{figure}[htp]
	\centering
	\begin{tikzpicture}[>=latex,scale=0.8]       
		\draw (0,0) -- (0,2) 
		(0,0) -- (2,0);
		\draw[dashed]  (0,2)--(2,2);
		\draw[line width=1.5] (2,0)--(2,2);
		\node at (0.25,0.25) {\small $\aaa$};
		\node at (1.75,0.25) {\small $\bbb$};
		\node at (0.25,1.75) {\small $\ccc$};
		\node at (1.75,1.75) {\small $\ddd$};
	\end{tikzpicture}\hspace{50pt}
	\begin{tikzpicture}[>=latex,scale=1.2]
		\draw
		(0,0.4) -- node[above=-2] {\small $a$} ++(1,0);
		
		\draw[line width=1.5]
		(0,-0.1) -- node[above=-2] {\small $b$} ++(1,0);
		
		\draw[dashed]
		(0,-0.6) -- node[above=-2] {\small $c$} ++(1,0);
	\end{tikzpicture}
	\caption{Quadrilaterals with the edge combination $a^2bc$.}
	\label{quadrilateral}
\end{figure}
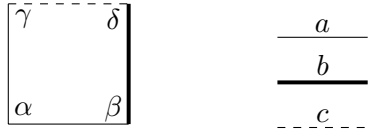

The lengths of $4$ edges of the quadrilateral in our tiling may have $4$ possible arrangements:  $a^2bc,a^2b^2,a^3b,a^4$ (see \cite{ua2} or our Lemma \ref{edge_combo}). Here $a^2bc$ means four edge lengths are $a,a,b,c$ in order with $a,b,c$ distinct. In particular a quadrilateral with $4$ mutually distinct edge lengths does not admit any edge-to-edge tiling. Sakano and Akama \cite{sa} classified  tilings for $a^2b^2$ and $a^4$ via the list of triangular tilings in \cite{ua}. Akama and Cleemput \cite{ac} had some partial study for $a^3b$ assuming convexity. We will classify tilings for $a^3b$ in the subsequent papers \cite{lw,lpwx} of this series. In this paper, we classify $a^2bc$. The quadrilateral is given by Fig. \ref{quadrilateral}, where $a,b,c$ are the normal, thick, and dashed lines. Throughout this paper, an $a^2bc$-tiling is always an {\em edge-to-edge} tiling of the sphere by congruent simple quadrilaterals in Fig. \ref{quadrilateral}, such that all vertices have degree $\ge 3$.

\begin{theorem*}
	There are exactly three classes of $a^2bc$-tilings:
	\begin{enumerate}
		\item A sequence of $2$-parameter families of $2$-layer earth map tilings 
		\newline  $T(2n \, \bbb\ccc\ddd, 2\aaa^n)$ with $2n$ tiles for any integer $n\ge3$; 
		
		\item A $1$-parameter family of quadrilateral subdivisions of the octahedron $T(8\aaa^3, 6\ddd^4, 12\bbb^2\ccc^2)$ with $24$ tiles, among which the special $\bbb=\frac{\pi}{3}$ case admits a  flip modification:  $T(2\aaa^3, 6\aaa\ccc^2, 6\ddd^4,6\bbb^2\ccc^2,6\aaa^2\bbb^2)$; 
		
		\item A sequence of $3$-layer earth map tilings (each has a unique quadrilateral) $T(4n\, \aaa\ccc^2, 2\bbb^{2n},2n\,\ddd^4,2n\,\aaa^2\bbb^2)$ with $8n$ tiles for any $n\ge2$, among which each odd $n=2m+1$ case admits exactly two flip modifications: 
		\begin{itemize}
			\item $T((8m+4)\aaa\ccc^2,4\aaa\bbb^{2m+2},(4m+2)\ddd^4,4m\,\aaa^2\bbb^2)$;
			\item $T((8m+2)\aaa\ccc^2,2\aaa\bbb^{2m+2},2\bbb^{2m}\ccc^2,(4m+2)\ddd^4,(4m+2)\aaa^2\bbb^2)$.
		\end{itemize}		
		
	\end{enumerate}
	The $2$nd and $3$rd classes have a unique quadrilateral in common with $\ddd=\frac{\pi}{2}$, $\aaa=\ccc=2\bbb=\frac{2\pi}{3}$, which is simply half of any pentagonal face of the regular dodecahedron and admits five different tilings in total. 				
\end{theorem*} 

\begin{figure}[htp]
	\centering
	\includegraphics[scale=0.197]{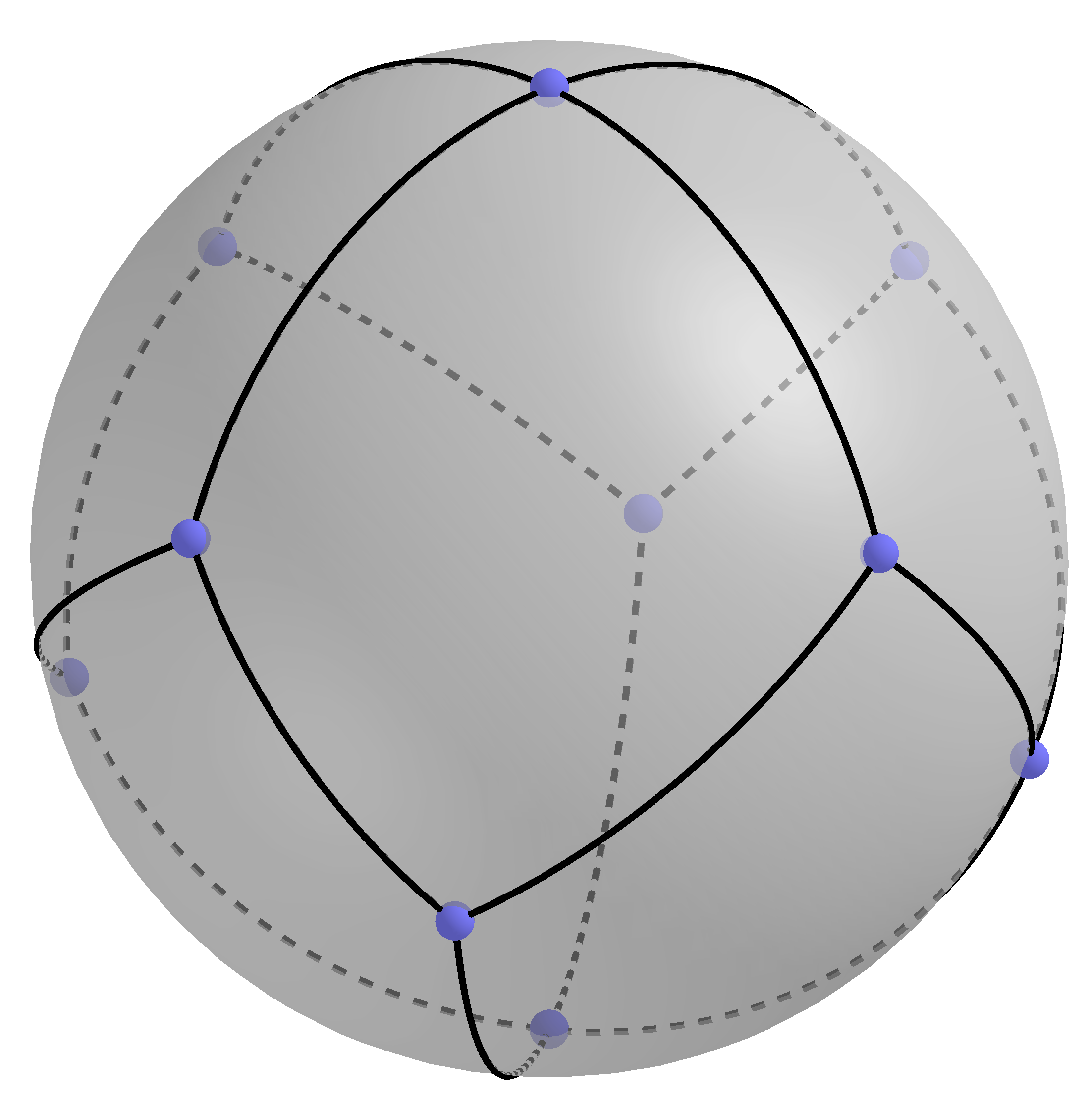} \hspace{9pt}   
	\includegraphics[scale=0.188]{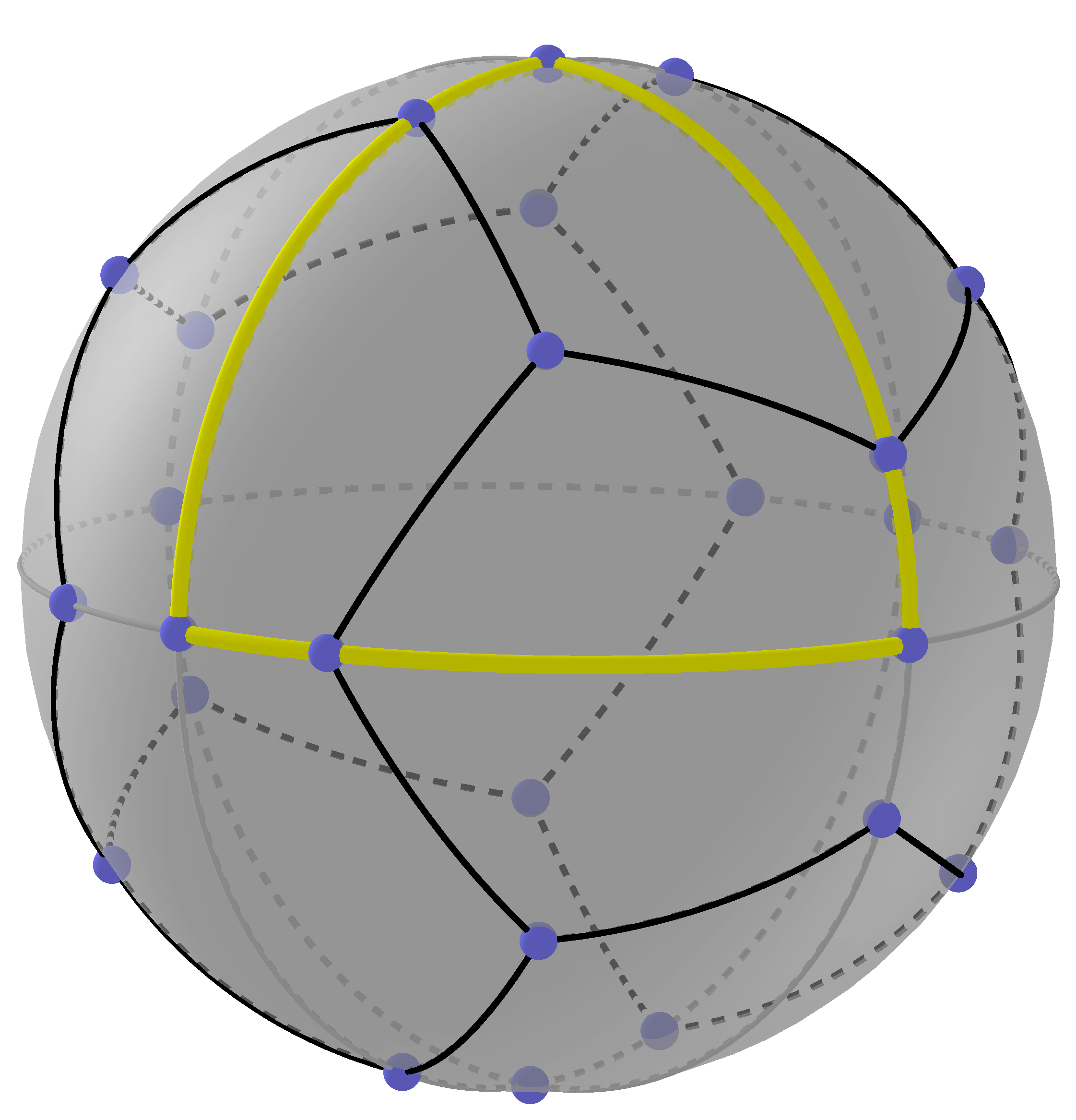}  \hspace{9pt}
	\includegraphics[scale=0.188]{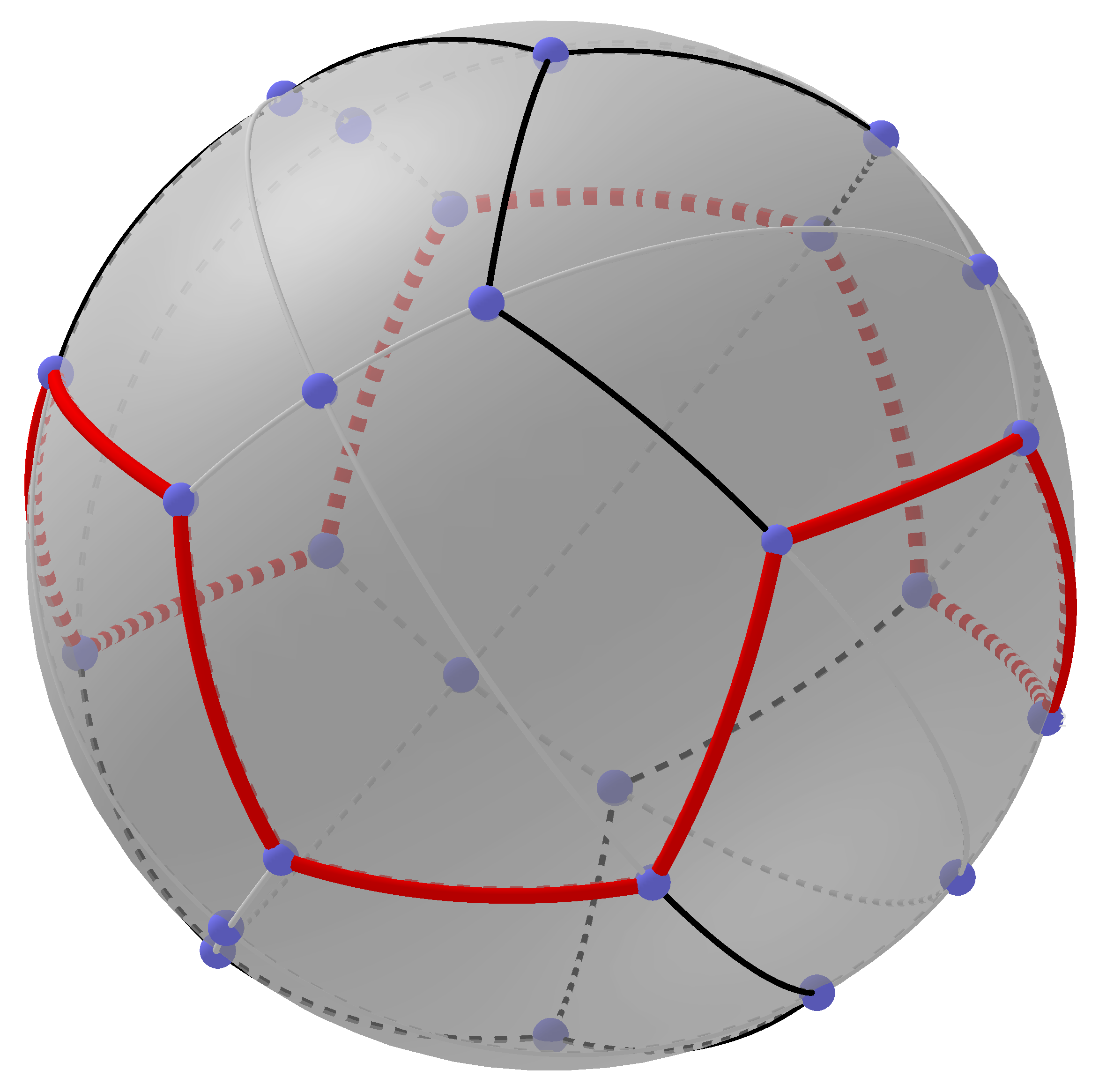}
	
	\includegraphics[scale=0.20]{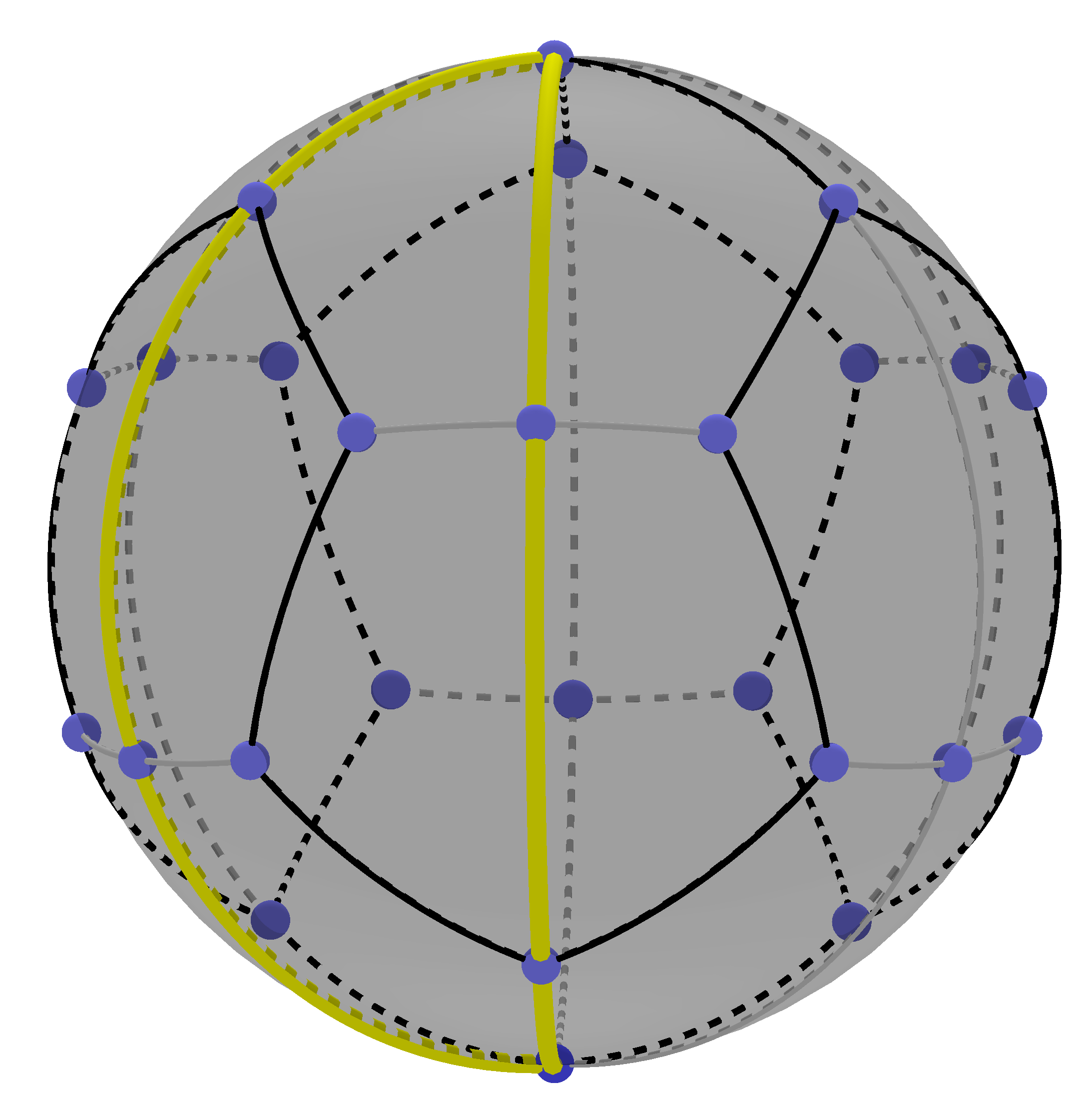}  \hspace{9pt}
	\includegraphics[scale=0.187]{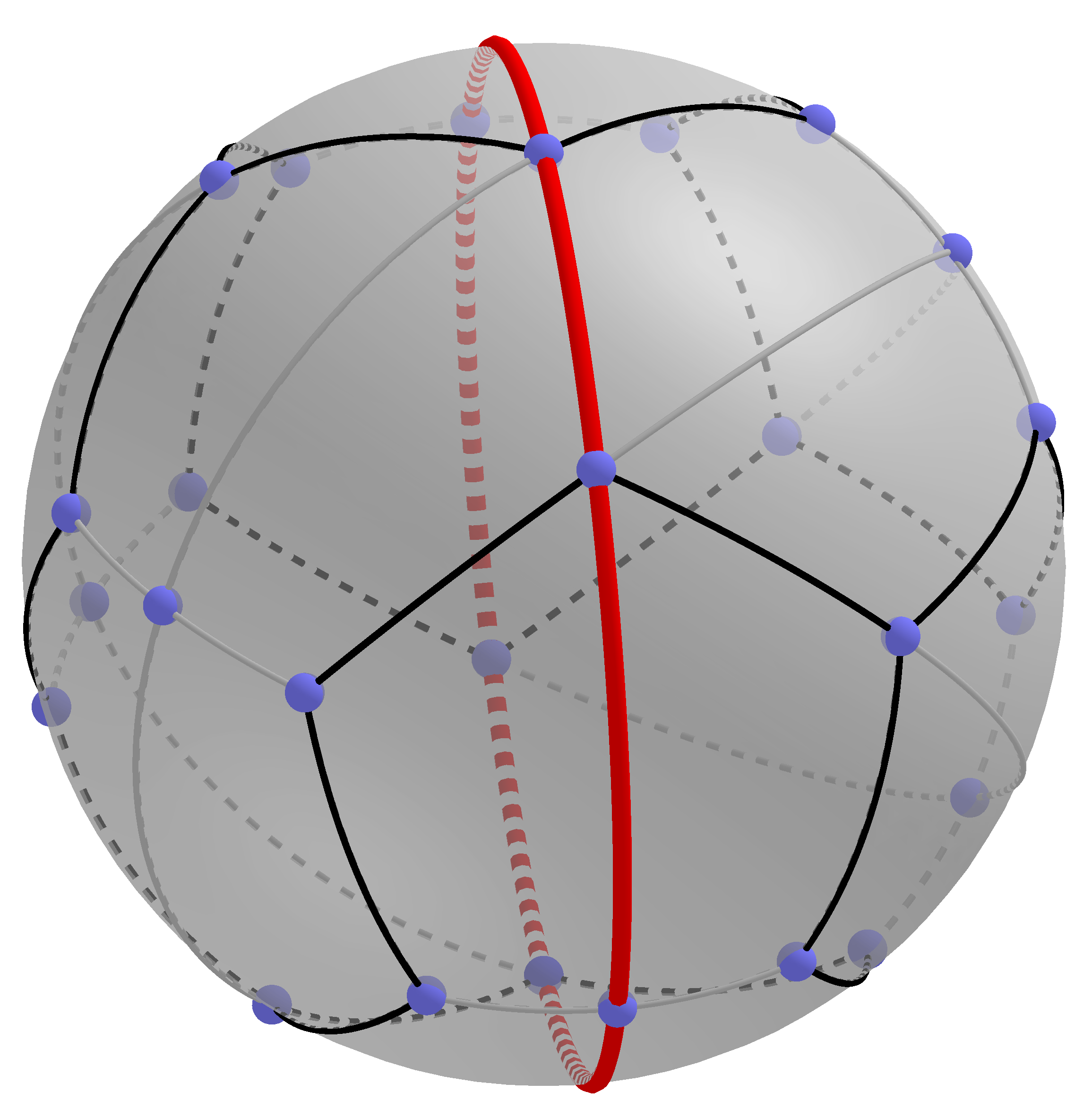}  \hspace{9pt}
	\includegraphics[scale=0.198]{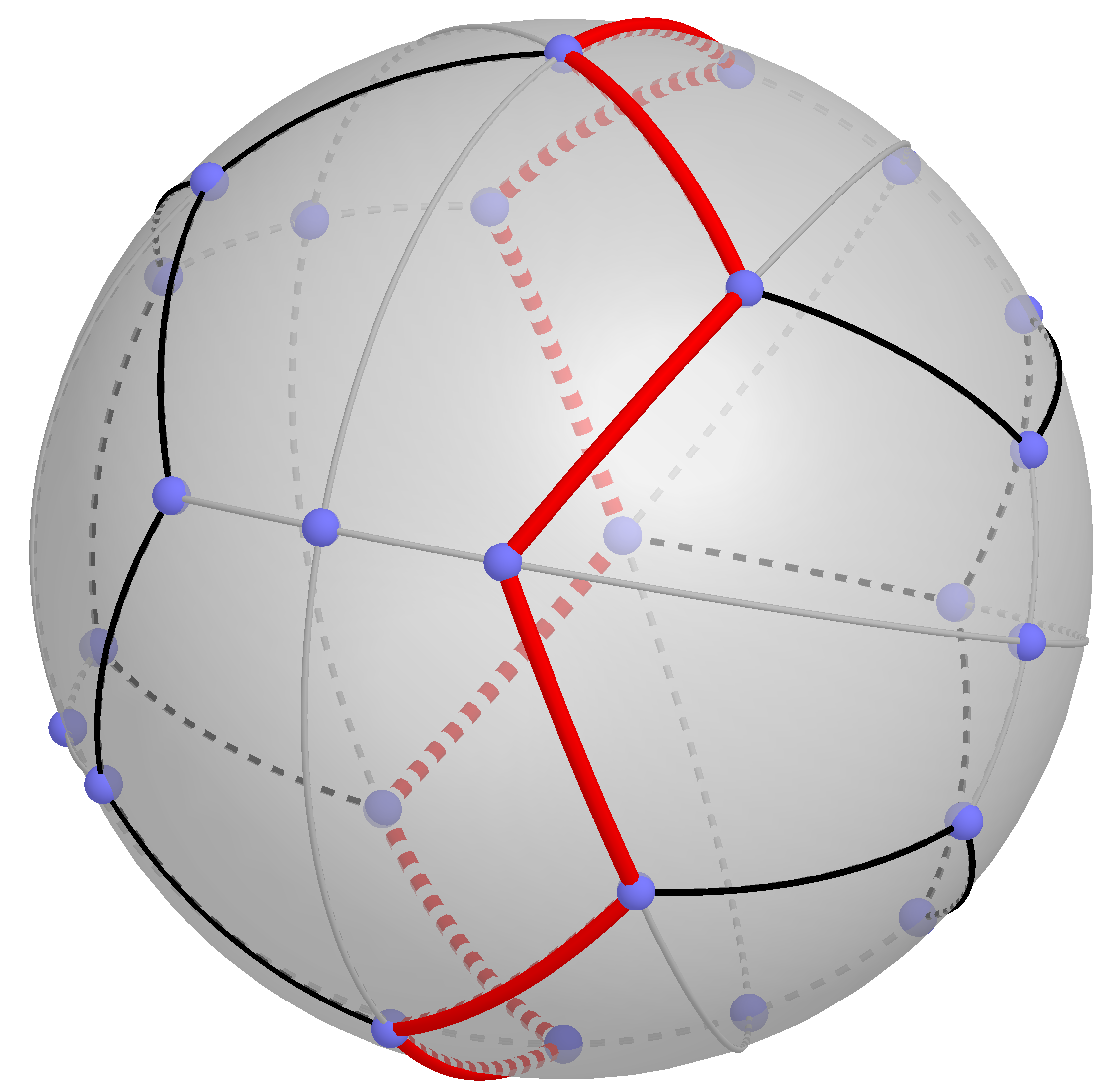}
	\caption{Six types of $a^2bc$-tilings drawn by GeoGebra.}
	\label{real figure}  
\end{figure}

The notation $T(2n \, \bbb\ccc\ddd, 2\aaa^n)$ means the tiling has exactly $2$ vertices $\aaa^n$ and $2n$ vertices $\bbb\ccc\ddd$, and is uniquely determined by them. Fig. \ref{real figure} shows some authentic pictures for all types of $a^2bc$-tilings: the first picture is a $2$-layer earth map tiling; the second is a quadrilateral subdivision of the octahedron; the third is the flip of half of the second with thick red boundary; the fourth is a $3$-layer earth map tiling; the fifth and sixth are two flips of half of the fourth with thick red boundaries. The second to the sixth pictures are five different tilings of the same quadrilateral in the end of the Theorem.

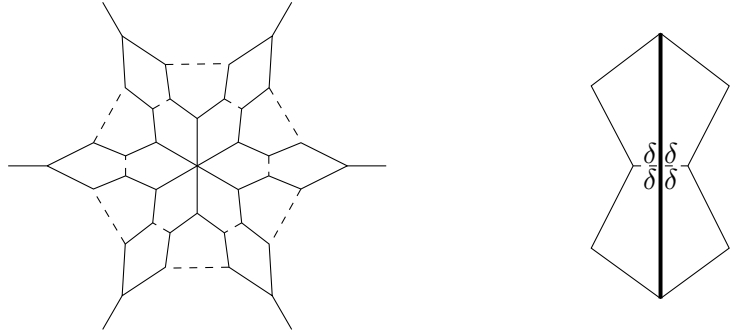
\begin{figure}[htp]
	\centering
	\begin{tikzpicture}[>=latex,scale=0.7]
		
		\foreach \b in {0,1,2,3,4,5}
		{
			\begin{scope}[rotate=60*\b ] 
				\draw (0,0)--(0.78,0.45)--(1.36,0.19)--(1.97,0.44)--(2.85,0)--(3.59,0)
				(0.78,-0.45)--(1.36,-0.19)--(1.97,-0.44)--(2.85,0);
				
				\draw[dashed] 
				(1.36,0.19)--(1.36,-0.19)
				(1.97,-0.44)--(1.39,-1.48);

			\end{scope}
		}

		\begin{scope}[>=latex,scale=0.4,xshift=22 cm]
			
			\draw (1.3,0)--(3.26,3.8)--(0,6.3)--(-3.29,3.8)--(-1.3,0)--(-3.29,-3.91)--(0,-6.3)--(3.26,-3.91)--(1.3,0);

			\draw[line width=1.5] (0,6.3)--(0,-6.3);

			\draw[dashed] (-1.3,0)--(1.3,0);

			\node at (0.5,0.6){\small $\ddd$}; \node at (-0.5,0.6){\small $\ddd$};
			\node at (-0.5,-0.6){\small $\ddd$}; \node at (0.5,-0.6){\small $\ddd$};
			
		\end{scope}	
	\end{tikzpicture}

	\caption{The induced earth map tilings by congruent symmetric pentagons.} \label{a4b}
\end{figure}

There is an interesting observation: All $\ddd$-vertices are $\ddd^4$ in the $2$nd and $3$rd classes, and four quadrilaterals of $\ddd^4$ in the right of Fig. \ref{a4b} can be viewed as a pair of symmetric $a^4b$-pentagons. All symmetric $a^4b$-tilings can be obtained in this way. For example, $3$-layer quadrilateral earth map tilings (and two flips) induce pentagonal earth map tilings (and two flips) in the left. 

\subsubsection*{Outline of the Paper}

The classification for $a^2bc$ is mainly the analysis of the neighborhood of a special tile (see Lemma \ref{base_tile}) with four vertices of degree $333d$, $334d$, $335d$, or $344d$. However, our subsequent papers \cite{lw,lpwx} on the classification of Type $a^3b$ are more about allowable combinations of angles at degree $3$ and $4$ vertices, and will apply interesting new techniques of cyclotomic field and trigonometric Diophantine equation when all angles of the quadrilateral are rational multiples of $\pi$. 

This paper is organized as follows. Section \ref{basic_facts} develops basic techniques needed for the classification work. This includes general results for all quadrilateral tilings of the sphere and some technical results specific to $a^2bc$. All other sections analyze the neighborhood of a special tile and complete the classification. Along the way we describe the moduli of $2$-layer earth map tilings and the quadrilateral subdivisions, and also provide exact calculations for the unique quadrilaterals in the $3$-layer earth map tilings. 


\subsubsection*{Acknowledgment}

We would like to thank Professor YAN Min for very helpful discussions on our early preprint of this work during his visit of our university in May 2021. This work has been announced in our submission to TJCDCGGG2021 (The 23rd Thailand-Japan Conference on Discrete and Computational Geometry, Graphs, and Games) on July 10, 2021. We would like to thank the organizers of this conference especially during this hard time of COVID-19 pandemic. Thank two junior students Fangbin Chen and Nan Zhang for showing us how to draw Fig. 2 using GeoGebra. Lastly we thank one referee for the long and detailed suggestions very much, which essentially improved the writings.

	\section{Basic Facts}\label{basic_facts}

	\subsection*{Vertex}
	Let $v,e,f$ be the numbers of vertices, edges, and tiles. Let $v_k$ be the number of vertices of degree $k$. We have Euler's formula and some basic counting: 
	\begin{align*}
	2&=v-e+f, \\
	2e=4f
	&=\sum_{k=3}^{\infty}kv_k=3v_3+4v_4+5v_5+\cdots, \\
	v
	&=\sum_{k=3}^{\infty}v_k=v_3+v_4+v_5+\cdots.
	\end{align*}
	Then it is easy to derive $v=f+2$ and  
	\begin{align}
	f &=6+\sum_{k=4}^{\infty}(k-3)v_k=6+v_4+2v_5+3v_6+\cdots, \label{vcountf} \\
	v_3 &=8+\sum_{k=5}^{\infty}(k-4)v_k=8+v_5+2v_6+3v_7+\cdots. \label{vcountv}
	\end{align}
	So $f\ge 6$ and $v_3 \ge 8$. These equalities show that there are many degree $3$ vertices, much more than the total of all degree $\ge 5$ vertices, but the number of degree $4$ vertices is uncertain. 

For $a^2bc$ or $a^3b$ quadrilateral tilings, each $b$-edge is shared by exactly two tiles. Then $f$ is twice of the number of $b$-edges, and is therefore even. 
	
	\begin{lemma}\label{base_tile}
		Any edge-to-edge quadrilateral tiling of the sphere with all vertices' degrees being $\ge 3$ has a special tile, whose four vertices have degree $333d$($d \ge 3$), $334d$($4\le d \le 11$), $335d$($d=5,6,7$), or $344d$($d=4,5$).
	\end{lemma}
	
	\begin{proof}
		Denote the degrees of four vertices of any tile $T$ by $d_1,d_2,d_3,d_4$. Counting the total number of vertices via each tile's contribution, we get
		\[
		\sum_{\textrm{all $f$ tiles } T} (\frac{1}{d_1}+\frac{1}{d_2}+\frac{1}{d_3}+\frac{1}{d_4}) = v = f+2.
		\]
		So there must exists a special tile $T$ such that $\frac{1}{d_1}+\frac{1}{d_2}+\frac{1}{d_3}+\frac{1}{d_4}>1$. The integer $\ge 3$ solutions $d_1,d_2,d_3,d_4$ are exactly $333d$($d \ge 3$), $334d$($4\le d \le 11$), $335d$($d=5,6,7$), or $344d$($d=4,5$). 
	\end{proof}
	
	When there are two tiles with different vertex-degrees in Lemma \ref{base_tile}, we will only call the earlier one in the above order its \textit{special tile}. For example, the tiling in the fifth picture of Fig. \ref{real figure} has tiles with vertex-degrees $3344$, $3345$, and $3445$. But we only say that it has \textit{special}  $3344$-Tile.

	\subsection*{Angle}
	
	The sum of all angles ({\em angle sum}) at a vertex is $2\pi$. The following is the {\em angle sum for quadrilateral}.
	
	\begin{lemma}\label{anglesum}
		If all tiles in a tiling of the sphere by $f$ quadrilaterals have the same four angles $\aaa,\bbb,\ccc,\ddd$, then 
		\[
		\aaa+\bbb+\ccc+\ddd = \left(2+\tfrac{4}{f}\right) \pi , 
		\]
		ranging in $(2\pi,\frac83 \pi]$. In particular no vertex contains all four angles.
	\end{lemma}
	
	\begin{proof}
		The angle sum at each vertex is $2\pi$,  and the total sum of all angles in a tiling is $2\pi v$. When the sum of four angles in each tile is the same $\Sigma=\alpha+\beta+\gamma+\delta$, the total sum of all angles is also $\Sigma f$. Then $2\pi v = \Sigma f$. By $v=f+2$, we get $\Sigma=(2+\frac{4}{f}) \pi$. So $2\pi<\Sigma\le\frac83 \pi$ for $f\ge 6$. 
	\end{proof}     
	
	Henceforth, to be concise and fluent, we will often use this angle sum lemma without mentioning it. 
	
	\subsection*{Edge}
	
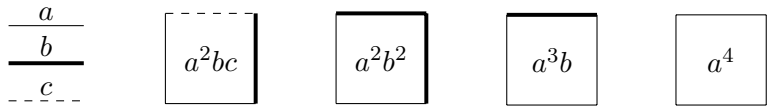
\begin{figure}[htp]
		\centering
		\begin{tikzpicture}[>=latex,scale=1]
		
		\begin{scope}[xshift=-2.2cm]
		
		\draw
		(0,0.4) -- node[above=-2] {\small $a$} ++(1,0);
		
		\draw[line width=1.5]
		(0,-0.1) -- node[above=-2] {\small $b$} ++(1,0);
		
		\draw[dashed]
		(0,-0.6) -- node[above=-2] {\small $c$} ++(1,0);
		
		\end{scope}
		\end{tikzpicture}\hspace{26pt}
		\begin{tikzpicture}[>=latex,scale=0.6]       
		\draw (0,0) -- (0,2) 
		(0,0) -- (2,0);
		\draw[dashed]  (0,2)--(2,2);
		\draw[line width=1.5] (2,0)--(2,2);
		\node at (1,1) {\small $a^2bc$};
		\end{tikzpicture}\hspace{25pt}
		\begin{tikzpicture}[>=latex,scale=0.6]       
		\draw (0,0) -- (2,0) 
		(0,0) -- (0,2);
		\draw[line width=1.5] (2,0)--(2,2)
		(0,2)--(2,2);
		\node at (1,1) {\small $a^2b^2$};
		\end{tikzpicture}\hspace{25pt}
		\begin{tikzpicture}[>=latex,scale=0.6]       
		\draw (0,0) -- (2,0) 
		(2,0)--(2,2)
		(0,0) -- (0,2);
		\draw[line width=1.5] (0,2)--(2,2);
		\node at (1,1) {\small $a^3b$};
		\end{tikzpicture}\hspace{25pt}
		\begin{tikzpicture}[>=latex,scale=0.6]       
		\draw (0,0) -- (2,0) 
		(2,0)--(2,2)
		(0,0) -- (0,2)
		(0,2)--(2,2);
		\node at (1,1) {\small $a^4$};
		\end{tikzpicture}
		\caption{Edge arrangements suitable for tiling.}
		\label{edges1}
	\end{figure}

	We restate four types of edge length arrangement for the quadrilaterals in our spherical tilings (see \cite{ua2}) here for the reader's convenience:

	\begin{lemma}\label{edge_combo}
		In a tiling of the sphere by congruent quadrilaterals (all vertices have degree $\ge 3$), the edge lengths of any tile are arranged in one of the four ways in Fig. \ref{edges1}, with distinct edge lengths $a,b,c$. 
	\end{lemma}

	\begin{proof}
		There are only five possible edge combinations ($a,b,c,d$ are distinct)  
		\[
		abcd,\;
		a^2bc,\;
		a^2b^2,\;
		a^3b,\;
		a^4.
		\]
		For $abcd$, without loss of generality, we may assume that the edges are arranged as in the first of Fig. \ref{edges2} and the vertex shared by $b,c$ has degree $3$ (recall $v_3 \ge 8$). Let $x$ be the third edge at the vertex. Then $x,b$ are adjacent in a tile, and $x,c$ are adjacent in another tile. Since there is no edge in the quadrilateral that is adjacent to both $b$ and $c$, we get a contradiction.  
		
		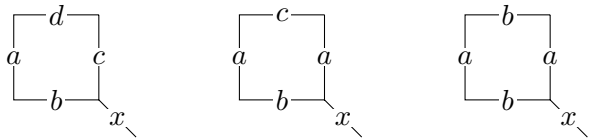
\begin{figure}[htp]
			\centering
			\begin{tikzpicture}[>=latex,scale=1]


			\draw
			(45:0.8) -- node[fill=white,inner sep=1] {\small $d$}
			(135:0.8) -- node[fill=white,inner sep=1] {\small $a$} 
			(225:0.8) -- node[fill=white,inner sep=1] {\small $b$} 
			(-45:0.8) -- node[fill=white,inner sep=1] {\small $c$} 
			(45:0.8)
			(-45:0.8) -- node[fill=white,inner sep=1] {\small $x$} 
			(-45:1.5);

			\begin{scope}[xshift=3cm]
			
			\draw
			(45:0.8) -- node[fill=white,inner sep=1] {\small $c$}
			(135:0.8) -- node[fill=white,inner sep=1] {\small $a$} 
			(225:0.8) -- node[fill=white,inner sep=1] {\small $b$} 
			(-45:0.8) -- node[fill=white,inner sep=1] {\small $a$} 
			(45:0.8)
			(-45:0.8) -- node[fill=white,inner sep=1] {\small $x$} 
			(-45:1.5);

			\end{scope}

			
			\begin{scope}[xshift=6cm]
			
			\draw
			(45:0.8) -- node[fill=white,inner sep=1] {\small $b$}
			(135:0.8) -- node[fill=white,inner sep=1] {\small $a$} 
			(225:0.8) -- node[fill=white,inner sep=1] {\small $b$} 
			(-45:0.8) -- node[fill=white,inner sep=1] {\small $a$} 
			(45:0.8)
			(-45:0.8) -- node[fill=white,inner sep=1] {\small $x$} 
			(-45:1.5);
			
			\end{scope}
			

			\end{tikzpicture}
			\caption{Not suitable for tiling.}
			\label{edges2}
		\end{figure}
		
		Similar contradictions occur in the second and the third of Fig. \ref{edges2}, since there is no edge in the quadrilateral that is adjacent to both $a$ and $b$. Here the first (adjacent $a$) of Fig. \ref{edges1} and second (separated $a$) of Fig. \ref{edges2} are two possible arrangements for the combination $a^2bc$;  the second of Fig. \ref{edges1} and third of Fig. \ref{edges2} are two  arrangements for the combination $a^2b^2$.  		
	\end{proof}

	\subsection*{Basic Techniques}
	
	We use the notations and techniques in \cite[Section 2]{wy1}, and add some discussion specific to $a^2bc$.
	
	\begin{lemma}\label{geometry3}
		An $a^2bc$-quadrilateral in Fig. \ref{quadrilateral} with $\aaa<\pi$ has $\bbb=\ccc$ if and only if $\ddd=\pi$. In other words, if $\ddd\neq\pi$, then $\bbb\neq\ccc$.  
	\end{lemma}
	
\begin{figure}[htp]
		\centering
		\begin{tikzpicture}[>=latex,scale=0.8]		
			\begin{scope}[xshift=-6 cm,scale=1] 		
				\draw (0,0)--(2,2)--(4,0);
				\draw[dashed] (2.5,0)--(4,0);
				\draw[line width=1.5] (0,0)--(2.5,0);
				\node at (0.6,0.25){\footnotesize $\bbb$};
				\node at (2,1.6){\footnotesize $\aaa$};
				\node at (3.4,0.25){\footnotesize $\ccc$};
				\fill (2.5,0) circle (0.03);
				\node at (0.8,1.3){$a$};  	\node at (3.2,1.3){$a$};		\node at (2,-0.3){$b+c$};  				
			\end{scope}    	 	
			\begin{scope}[xshift=0 cm]
				
				\draw[dotted]
				(0,0) -- (2.65,1.53);
				
				\node at (0.4,2){$A$};  \node at (-0.4,0){$B$};  \node at (4.3,0){$D$}; 
				\node at (2.8,1.8){$C$}; 
				
			\end{scope}     	
			\draw
			(0,0) -- (70:2) -- ++(-10:2);
			\draw[dashed] (2.65,1.53) -- (4,0);

			\draw[line width=1.5]
			(0,0) -- (3.95,0);	
			
			\node at (0.8,1.5) {\small $\aaa$};
			\node at (2.05,1.4) {\small $\theta$};
			\node at (0.35,0.45) {\small $\theta$};
			\node at (2.55,1.2) {\small $\ccc'$};
			\node at (0.8,0.25) {\small $\bbb'$};
			\node at (3.4,0.25) {\small $\ddd$};
			
			\begin{scope}[scale=0.5,xshift=16 cm]
				\draw (0,4)--(4,0)
				(0,4)--(-4,0);
				\draw[line width=1.5]
				(-4,0)--(-0.5,1);
				\draw[dashed]
				(4,0)--(-0.5,1);
				\draw[dotted]
				(4,0)--(-4,0);
				
				\node at (0,3.4) {\small $\aaa$};
				\node at (-0.5,1.5) {\small $\ddd$};
				\node at (-3.2,0.3) {\small $\theta$};
				\node at (3.2,0.3) {\small $\theta$};
				
				\node at (-1.5,0.2) {\small $\bbb'$};
				\node at (1.5,0.2) {\small $\ccc'$};
			\end{scope}
			
		\end{tikzpicture}
		\caption{Proof of Lemma \ref{geometry3} and \ref{geometry4}.}
		\label{g3proof}
	\end{figure}
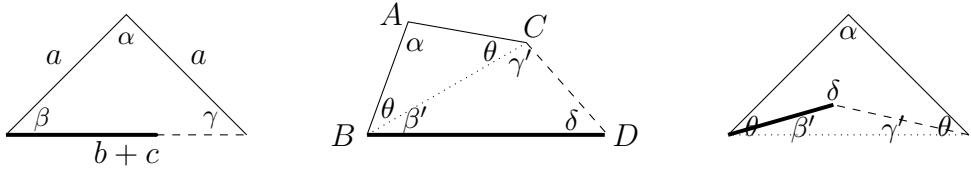

	\begin{proof}     	
		If $\ddd=\pi$, we get an isosceles triangle in the first picture of Fig. \ref{g3proof}, thus $\bbb=\ccc$. If $\bbb=\ccc$ and $\ddd\neq\pi$, then $\bbb'=|\bbb-\theta|=|\ccc-\theta|=\ccc'$ as shown in the second and third pictures of Fig. \ref{g3proof}. So we get $b=c$, a contradiction. 
	\end{proof}
	
	\begin{lemma}\label{geometry4}
		A convex $a^2bc$-quadrilateral in Fig. \ref{quadrilateral} has $\aaa+2\bbb>1$, $\aaa+2\ccc>1$.
	\end{lemma}
	
	\begin{proof}
		The line $BC$ is inside the quadrilateral in the second picture of Fig. \ref{g3proof}.  Thus $\theta<\bbb,\ccc$. This implies  $\aaa+2\bbb>\aaa+2\theta>1$ and $\aaa+2\ccc>\aaa+2\theta>1$.     	
	\end{proof}

	\begin{lemma}[Parity Lemma]\label{edge_combo1}
		
		In an $a^2bc$-tiling, the respective numbers of $\bbb,\ccc,\ddd$ at any vertex  have the same parity, i.e. they are all odd or all even.   
		
	\end{lemma}
	
	\begin{proof}
		The total number of $\bbb,\ddd$ together at a vertex is twice the number of $b$-edges at the vertex. So the respective numbers of $\bbb$ and $\ddd$ must have the same parity. Similar argument applies to $\ccc,\ddd$. 
	\end{proof}

	We call a vertex \textit{even} or \textit{odd} whenever the degrees of $\bbb,\ccc,\ddd$ are even or odd. Then Lemma \ref{anglesum} implies any $\aaa$-vertex $\aaa \cdots$ is always even.

	\begin{lemma}\label{ad}
		In an $a^2bc$-tiling, a vertex without $\bbb,\ccc$ must be $\alpha^k$ or $\delta^k$.
	\end{lemma}
	
	\begin{proof}
		If a vertex has only $a$-edge, then it has only $a^2$-angles $\aaa$. Therefore the vertex is $\aaa^k$. If a vertex has no $a$-edge, then it has only $bc$-angle $\delta$. Therefore the vertex is $\ddd^k$. In all other cases, it has $ab$-angle $\bbb$ or $ac$-angle $\ccc$.		
	\end{proof}

		In a tiling of the sphere by $f$ congruent tiles, each angle of the tile appears $f$ times in total. If one vertex has more $\aaa$ than $\bbb$, there must exist another vertex with more $\bbb$ than $\aaa$.
	    Such global counting induces many interesting and useful results. 
	
	\begin{lemma}[Balance Lemma]\label{balance}
		In an $a^2bc$-tiling with $f$ tiles, one of $\bbb^2\cdots$, $\ccc^2\cdots$, $\ddd^2\cdots$ is a vertex if and only if all three are vertices. Moreover,  if all three are not vertices, then $\aaa^{f/2}$ and $\bbb\ccc\ddd$ are the only vertices. 
	\end{lemma}

	\begin{proof}						
		If  $\bbb^2\cdots$ is not a vertex, then any vertex $\aaa^k\bbb^l\ccc^m\ddd^n$ has $l=0,1$. If $l=0$, then $m\ge l$. If $l=1$, then $m$ is odd by  Parity Lemma, which also implies $m\ge 1=l$. Then $m=l\le 1$ at every vertex. So $\ccc^2\cdots$ is never a vertex. Similarly $n=l\le 1$ at every vertex and $\ddd^2\cdots$ is never a vertex. When $l=m=n=0$, the vertex is $\aaa^k$ ($k\ge 3$). When $l=m=n=1$, Lemma \ref{anglesum} forces $k=0$ and the vertex is $\bbb\ccc\ddd$. Finally both $\aaa^k$ and $\bbb\ccc\ddd$ must appear. Similar proof works when $\ccc^2\cdots$ or $\ddd^2\cdots$ is not a vertex. By $\bbb\ccc\ddd$ and Lemma \ref{anglesum}, we get $k=f/2$. 		
	\end{proof}
	
	\begin{lemma}\label{lem3}
		In an $a^2bc$-tiling, there are only four possible types of degree $3$ vertices $\aaa^{3}$, $\aaa\bbb^2$, $\aaa\ccc^2$ and $\bbb\ccc\ddd$ shown in Fig. \ref{deg3}. 
	\end{lemma}
	\begin{proof}
		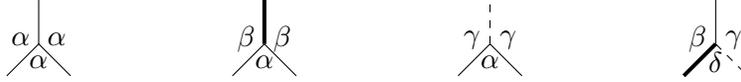
\begin{figure}[htp]
			\centering
			\begin{tikzpicture}[>=latex,scale=0.6]	
				\begin{scope}

					\draw (0,0)--(0,1)
					(0,0)--({-sqrt(2)/2},{-sqrt(2)/2})
					(0,0)--({sqrt(2)/2},{-sqrt(2)/2});
					\node at (0,-0.4){\small $\aaa$};
					\node at (0.4,0.1){\small $\aaa$};
					\node at (-0.4,0.1){\small $\aaa$};
				\end{scope}
			    \begin{scope}[xshift=5cm]
			    	\draw 
			    	(0,0)--({-sqrt(2)/2},{-sqrt(2)/2})
			    	(0,0)--({sqrt(2)/2},{-sqrt(2)/2});
			    	\node at (0,-0.4){\small $\aaa$};
			    	\node at (0.4,0.1){\small $\bbb$};
			    	\node at (-0.4,0.1){\small $\bbb$};
					\draw[line width=1.5] (0,0)--(0,1);
			    \end{scope}	
					
				\begin{scope}[xshift=10cm]
					\draw 
					(0,0)--({-sqrt(2)/2},{-sqrt(2)/2})
					(0,0)--({sqrt(2)/2},{-sqrt(2)/2});
					\node at (0,-0.4){\small $\aaa$};
					\node at (0.4,0.1){\small $\ccc$};
					\node at (-0.4,0.1){\small $\ccc$};
					\draw[dashed] (0,0)--(0,1);
				\end{scope}		
			\begin{scope}[xshift=15cm]
				\draw (0,0)--(0,1);
				\node at (0,-0.4){\small $\ddd$};
				\node at (0.4,0.1){\small $\ccc$};
				\node at (-0.4,0.1){\small $\bbb$};
				\draw[dashed] (0,0)--({sqrt(2)/2},{-sqrt(2)/2});
				\draw[line width=1.5] (0,0)--({-sqrt(2)/2},{-sqrt(2)/2});
			\end{scope}						
				
			\end{tikzpicture}
			\caption{Four possible types of degree $3$ vertices.}
			\label{deg3}
		\end{figure}
		Since there is neither $bb$-angle nor $cc$-angle, the $3$ edges at any degree $3$ vertex must be $aaa$, $aab$, $aac$ or $abc$ in Fig. \ref{deg3}, which determine four degree $3$ vertices uniquely. 
	\end{proof}

    \begin{lemma}\label{lem4}
    	In an $a^2bc$-tiling, besides $\aaa^4$, $\bbb^4$, $\ccc^4$, $\ddd^4$, there are only five possible types of degree $4$ vertices $\aaa^2\bbb^2$, $\aaa^2\ccc^2$, $\bbb^2\ccc^2$, $\bbb^2\ddd^2$, $\ccc^2\ddd^2$ shown in Fig. \ref{deg4}. Each of them is uniquely determined by two different angles in it.
    \end{lemma}
    \begin{proof}
    	One proof is to list all possibilities of $4$ edges at a vertex. Note that $\aaa\bbb\ccc\ddd$ is never a vertex by Lemma \ref{anglesum}. Another proof is to apply Parity Lemma. A degree $4$ vertex must be even, so it has to be $\theta^4$ or $\theta^2 \phi^2$ for any angles $\theta \neq \phi$ of the quadrilateral. Note that there is no $\aaa^2\ddd^2$ by Lemma \ref{ad}. 
    	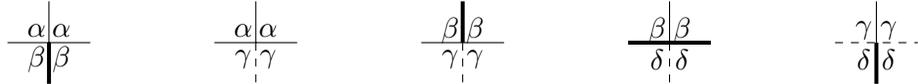
\begin{figure}[htp]
    		\centering
    		\begin{tikzpicture}[>=latex,scale=0.55]	
    			\begin{scope}   				  				
    				\draw (0,0)--(0,1)
    				(-1,0)--(1,0);
    				\draw[line width=1.5] (0,0)--(0,-1);
    				\node at (-0.3,-0.4){\small $\bbb$};
    				\node at (0.3,-0.4){\small $\bbb$};
    				\node at (0.3,0.3){\small $\aaa$};
    				\node at (-0.3,0.3){\small $\aaa$};
    			\end{scope}
    			\begin{scope}[xshift=5cm]
    				\draw (0,0)--(0,1)
    				(-1,0)--(1,0);
    				\draw[dashed] (0,0)--(0,-1);
    				\node at (-0.3,-0.4){\small $\ccc$};
    				\node at (0.3,-0.4){\small $\ccc$};
    				\node at (0.3,0.3){\small $\aaa$};
    				\node at (-0.3,0.3){\small $\aaa$};
    			\end{scope}	
    			
    			\begin{scope}[xshift=10cm]
    				\draw 
    				(-1,0)--(1,0);
    				\draw[dashed] (0,0)--(0,-1);
    				\draw[line width=1.5] (0,0)--(0,1);
    				\node at (-0.3,-0.4){\small $\ccc$};
    				\node at (0.3,-0.4){\small $\ccc$};
    				\node at (0.3,0.3){\small $\bbb$};
    				\node at (-0.3,0.3){\small $\bbb$};
    			\end{scope}		
    			\begin{scope}[xshift=15cm]
    				\draw 
    				(0,0)--(0,1);
    				\draw[dashed] (0,0)--(0,-1);
    				\draw[line width=1.5] (-1,0)--(1,0);
    				\node at (-0.3,-0.4){\small $\ddd$};
    				\node at (0.3,-0.4){\small $\ddd$};
    				\node at (0.3,0.3){\small $\bbb$};
    				\node at (-0.3,0.3){\small $\bbb$};
    			\end{scope}	
    		\begin{scope}[xshift=20cm]
    			\draw 
    			(0,0)--(0,1);
    			\draw[dashed] (-1,0)--(1,0);
    			\draw[line width=1.5] (0,0)--(0,-1);
    			\node at (-0.3,-0.4){\small $\ddd$};
    			\node at (0.3,-0.4){\small $\ddd$};
    			\node at (0.3,0.3){\small $\ccc$};
    			\node at (-0.3,0.3){\small $\ccc$};
    		\end{scope}						
    			
    		\end{tikzpicture}
    		\caption{Five possible types of degree $4$ vertices with two different angles.}
    		\label{deg4}
    	\end{figure}
    \end{proof}

    \begin{proposition}\label{lembc3}
    	There is no $a^2bc$-tiling having both vertices $\aaa\bbb^2$ and $\aaa\ccc^2$. 
    \end{proposition}
    \begin{proof}
    	If $\aaa\bbb^2$ and $\aaa\ccc^2$ are both vertices, then $\bbb=\ccc$ and $\ddd=\pi$ by Lemma \ref{geometry3}. So $\ddd^2\cdots$ is not a vertex. Then Balance Lemma \ref{balance} implies that $\aaa^{f/2}$ and $\bbb\ccc\ddd$ are the only vertices, a contradiction. 
    \end{proof}
	
	The very useful tool \textit{adjacent angle deduction} (abbreviated as AAD) has been introduced in \cite[Section 2.5]{wy1}. The following is \cite[Lemma 10]{wy1}.

	\begin{lemma}\label{aadlemma}
		The AAD of the vertex $\alpha^n$ has the following properties: 
		\begin{itemize}
			\item If $\beta\thin\beta\cdots$ or $\gamma\thin\gamma\cdots$ is not a vertex, then $\alpha^n$ has the unique AAD $\thin^{\beta}\alpha^{\gamma}\thin^{\beta}\alpha^{\gamma}\thin^{\beta}\alpha^{\gamma}\thin\cdots$.
			\item If $n$ is odd, then we have the AAD $\thin^{\beta}\alpha^{\gamma}\thin^{\beta}\alpha^{\gamma}\thin$ at $\alpha^n$.
		\end{itemize}
	\end{lemma}

	We remark that $\theta^n$ for $\theta=\thin^{\aaa}\bbb^{\delta}\thick$, $\thin^{\aaa}\ccc^{\delta}\dash$ or $\thick^{\beta}\delta^{\ccc}\dash$ has a unique AAD.

	\begin{lemma} \label{lembd}
		In an $a^2bc$-tiling, if $\bbb>\frac{\pi}{2},\bbb+\ddd>\pi$, then $\bbb\thin\bbb\cdots$ is not a vertex. If $\ddd>\frac{\pi}{2},\bbb+\ddd>\pi$, then $\ddd\dash\ddd\cdots$ is not a vertex. 
	\end{lemma}
	
	\begin{proof}
		We have $\thick\bbb\thin\bbb\thick\cdots=\theta\thick\bbb\thin\bbb\thick\rho\cdots$ with $\theta,\rho=\bbb$ or $\ddd$, and $\theta,\rho$ are not the same angle (i.e., the vertex is not degree $3$). Then the sum of angles is $>2\pi$ by $\bbb>\frac{\pi}{2}$ and $\bbb+\ddd>\pi$, a contradiction.	The case $\ddd\dash\ddd\cdots$ is similar.	
	\end{proof}
 			
	\begin{lemma}\label{lemac2}
		In an $a^2bc$-tiling, if $\aaa\ccc^2$ is a vertex, then $\aaa \cdots=\aaa\ccc^2$ or $\aaa^k \bbb^{2t}$. Furthermore $\aaa^k \bbb^{2t}$ for some $k \ge 1,t\ge 0$ must appear. 
	\end{lemma}
	
	\begin{proof}
		Recall that any $\aaa$-vertex $\aaa \cdots$ is even. Let  $\aaa\cdots=\aaa^k\bbb^l\ccc^m\ddd^n$, and $l,m,n$ are all even.  
		
		When $m \ge 2$, we get $\aaa\cdots=\aaa\ccc^2 \cdots= \aaa\ccc^2$. 
		
		When $m =0$, $\aaa\cdots=\aaa^k\bbb^l\ddd^n=\aaa^k\bbb^{2t}\ddd^{2s}$. If $s > 0$, then Lemma \ref{ad} implies $t>0$, and we have $\aaa+2\bbb+2\ddd \le 2\pi$. By $\aaa+2\ccc= 2\pi$, we deduce that $\aaa+\bbb+\ccc+\ddd \le 2\pi$, contradicting Lemma \ref{anglesum}. So $s=0$ and the vertex is $\aaa^k \bbb^{2t}$ for some $k \ge 1,t\ge 0$. 
		
		Finally $\aaa\ccc^2$ imply the existence of a different vertex with more $\alpha$ than $\ccc$, which has to be $\aaa^k \bbb^{2t}$. (There may exist several $\aaa^k \bbb^{2t}$ with different $k,t$.)				
	\end{proof}

\noindent	We will use Lemma/Proposition $n'$ to denote the use of Lemma/Proposition $n$ after exchanging $\bbb\leftrightarrow\ccc$.

	\section{$333d$-Tile}
	\label{2abc 1}

This section classifies all tilings with a special $333d$-Tile (Lemma \ref{base_tile}) as $2$-layer earth map tilings. To facilitate discussion, we denote by $T_i$ the tile labeled $i$, by $E_{ij}$ the edge shared by $T_i,T_j$. We denote by $\theta_i$ the angle $\theta$ in $T_i$. We say a tile is {\em determined} when we know all the edges and angles of the tile.	
		
	\begin{proposition}\label{333d}  
		For an $a^2bc$-tiling, the following statements are equivalent: 
		\begin{enumerate}[$(1)$]
			\item Every tile is a $333d$-Tile. 
			\item There exists a $333d$-Tile. 
			\item The $bc$-angle $\delta$ appears at some degree $3$ vertex (or $\bbb\ccc\ddd$ is a vertex).
			\item It is the $2$-layer earth map tiling $T(f \, \bbb\ccc\ddd, 2\aaa^{f/2})$ ($f \ge 6$) in Fig. \ref{bcdpic}. 
		\end{enumerate} 
	\end{proposition}
	\begin{proof}
		$(1)\Rightarrow(2)$ is trivial.
		
		$(2)\Rightarrow(3)$: If $\ddd$ never appears in degree $3$ vertices, then the $333d$-Tile has both $\bbb$ and $\ccc$ as degree $3$ vertices without $\ddd$, which must be $\aaa\bbb^2$ and $\aaa\ccc^2$ respectively by Lemma \ref{lem3}. But this contradicts Proposition \ref{lembc3}. 
		
		$(3)\Rightarrow(4)$:   
		For any degree $3$ vertex containing $\thick\ddd\dash$, the third edge can only be $a$-edge,  and  the vertex must be $\bbb\ccc\ddd$. 
		
		Next we show $\bbb\ddd\cdots=\bbb\ccc\ddd$. Let $\bbb\ddd\cdots=\aaa^k\bbb^l\ccc^m\ddd^n$. If $m\ge1$, we have $\bbb\ddd\cdots=\bbb\ccc\ddd$. If $m=0$, Parity Lemma implies $l,n\ge2$. So we have $\bbb+\ddd\le\pi$, and then $\ccc\ge\pi$ by $\bbb\ccc\ddd$. However, the unique AAD $\thick^{\bbb}\ddd^{\ccc}\dash^{\ccc}\ddd^{\bbb}\thick\cdots$ of $\aaa^k\bbb^l\ddd^n$ gives $\ccc^2\cdots$, a contradiction. Therefore, $\bbb\ddd\cdots=\bbb\ccc\ddd$. 
		
		Similarly, $\ccc\ddd\cdots=\bbb\ccc\ddd$. In Fig. \ref{bcdpic}, $\bbb_1\ccc_3\ddd_2$ determines $T_1,T_2,T_3$. Then $\ccc_2\ddd_3\cdots=\bbb_4\ccc_2\ddd_3$ determines $T_4$;  $\bbb_3\ddd_4\cdots=\bbb_3\ccc_5\ddd_4$ determines $T_5$. The argument started at $\bbb_1\ccc_3\ddd_2$ can be repeated at $\bbb_3\ccc_5\ddd_4$. More repetitions give the unique tiling of $f=2d$ tiles with $2\aaa^d$ ($d\ge3$) and $2d \, \bbb\ccc\ddd$.

		\begin{figure}[htp]
			\centering
			\begin{tikzpicture}[>=latex,scale=0.6] 
			\foreach \a in {0,1,2}
			{
				\begin{scope}[xshift=2*\a cm] 
				\draw (0,0)--(0,-2)
				(2,0)--(2,-2)
				(0.5,-3)--(0.5,-5);
				\draw[dashed]  (0,-2)--(0.5,-3);
				\draw[line width=1.5] (0.5,-3)--(2,-2);
				\node at (0.6,-2.55){\small $\ddd$};
				\node at (0.2,-3.2){\small $\ccc$};
				\node at (0.8,-3.2){\small $\bbb$};
				\node at (1,0){\small $\aaa$};
				\node at (1.5,-5){\small $\aaa$};
				
				\node at (1.9,-2.45){\small $\ddd$};
				\node at (1.7,-1.8){\small $\bbb$};
				\node at (2.3,-1.8){\small $\ccc$};
				\end{scope}
			}
			
			\draw[dashed] (6,-2)--(6.5,-3);
			\draw  (6.5,-3)--(6.5,-5);
			
			\draw[dotted] (-1,-2.5)--(8,-2.5);
			
			\fill (8,-2) circle (0.05); 
			\fill (8.2,-2) circle (0.05);
			\fill (8.4,-2) circle (0.05);
			
			\node[draw,shape=circle, inner sep=0.5] at (1,-1) {\small $1$};
			\node[draw,shape=circle, inner sep=0.5] at (3,-1) {\small $3$};
			\node[draw,shape=circle, inner sep=0.5] at (5,-1) {\small $5$};
			\node[draw,shape=circle, inner sep=0.5] at (1.5,-4) {\small $2$};
			\node[draw,shape=circle, inner sep=0.5] at (3.5,-4) {\small $4$};
			\node[draw,shape=circle, inner sep=0.5] at (5.5,-4) {\small $6$};
			
			\node at (6.2,-3.2){\small $\ccc$};
			\node at (0.3,-1.8){\small $\ccc$};
			
			\end{tikzpicture}
			\caption{ The $2$-layer earth map tiling $T(f \, \bbb\ccc\ddd, 2\aaa^{f/2})$.} \label{bcdpic}
		\end{figure}
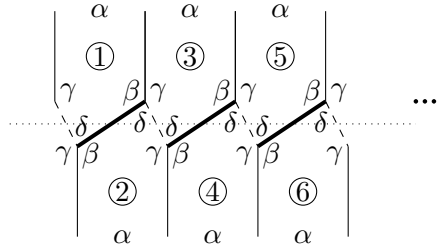
	$(4)\Rightarrow(1)$: Any tile in the $2$-layer earth map tiling is a $333d$-Tile. 	
\end{proof}
	

    \begin{proposition}\label{tuilun0}
    	For an $a^2bc$-tiling, if $\ddd=\pi$ (or equivalently $\bbb=\ccc$), then it is a $2$-layer earth map tiling. 
    \end{proposition}
    \begin{proof}
    	By Lemma \ref{geometry3}, $\bbb=\ccc$ if and only if $\ddd=\pi$. But $\ddd=\pi$ implies that   $\ddd^2\cdots$ is not a vertex. Then Balance Lemma \ref{balance} implies that $\ddd\cdots=\bbb\ccc\ddd$, which determines a $2$-layer earth map tiling by Proposition \ref{333d}.
    \end{proof}
    

\subsection*{Geometric realization and the moduli of $T(2n \, \bbb\ccc\ddd, 2\aaa^n)$}

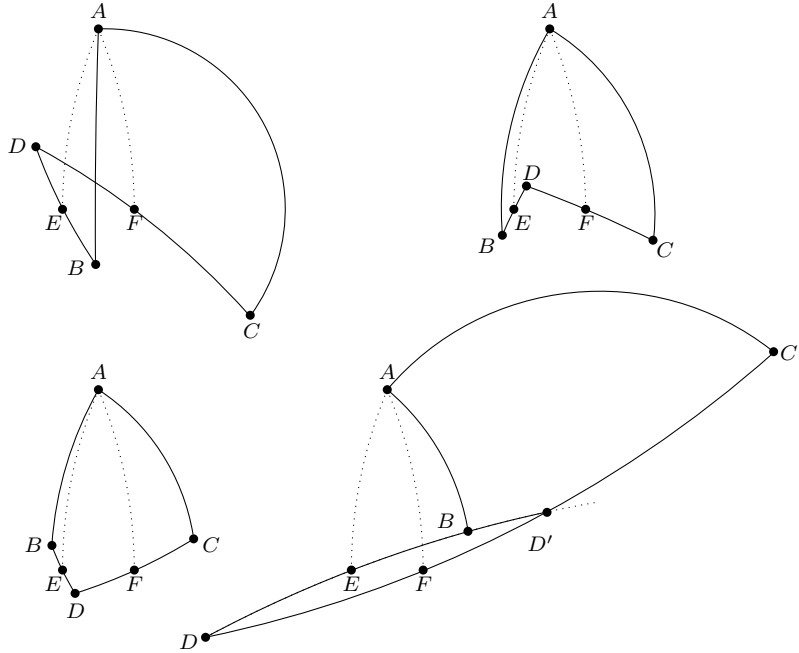
\begin{figure}[htp]
	\centering
	\begin{tikzpicture}[>=latex,scale=1.2] 	
		\begin{scope}[xshift=0 cm]
			\draw[dotted] ({-2*sin(180/8)/(1+cos(180/8))},0) arc(180:157.3:{2*sin(180/8)/(1-cos(180/8)*cos(180/8))});
			\draw[dotted]
			({2*sin(180/8)/(1+cos(180/8))},0) arc(0:22.7:{2*sin(180/8)/(1-cos(180/8)*cos(180/8))});
			
			\fill (0,2) circle (0.05);  \fill ({2*sin(180/8)/(1+cos(180/8))},0) circle (0.05); 
			\fill ({-2*sin(180/8)/(1+cos(180/8))},0) circle (0.05); 
			
			\fill (1.144,-0.343) circle (0.05); 
			\fill (-0.525,-0.288) circle (0.05);
			
			\fill ({-2/(4+sqrt(14))},{2/(4+sqrt(14))}) circle (0.05); 
			
			\draw({-2/(4+sqrt(14))},{2/(4+sqrt(14))}) arc(151:157:5.86);
			
			\draw({-2/(4+sqrt(14))},{2/(4+sqrt(14))}) arc(70:63.3:13.375);
			
			\draw(0,2) arc(149.7:185:3.925); 
			
			\draw(0,2) arc(59.1:-6:2.396); 
			
			\node at (0,2.2){\scriptsize $A$}; \node at (-0.3,-0.15){\scriptsize $E$}; \node at (0.4,-0.15){\scriptsize $F$}; \node at (-0.7,-0.4){\scriptsize $B$}; \node at (1.28,-0.45){\scriptsize $C$};  \node at (-0.2,0.4){\scriptsize $D$};

		\end{scope}   
		\begin{scope}[xshift=-5cm]
			\draw[dotted] ({-2*sin(180/8)/(1+cos(180/8))},0) arc(180:157.3:{2*sin(180/8)/(1-cos(180/8)*cos(180/8))});
			\draw[dotted]
			({2*sin(180/8)/(1+cos(180/8))},0) arc(0:22.7:{2*sin(180/8)/(1-cos(180/8)*cos(180/8))});
			
			\fill (0,2) circle (0.05);  \fill ({2*sin(180/8)/(1+cos(180/8))},0) circle (0.05); 
			\fill ({-2*sin(180/8)/(1+cos(180/8))},0) circle (0.05);
			
			\fill ({-2*sqrt(5)/(4+sqrt(6))},{2*sqrt(5)/(4+sqrt(6))}) circle (0.05);
			
			\fill (-0.03,-0.612) circle (0.05); 
			
			\fill (1.682,-1.175) circle (0.05); 
			
			\draw({-2*sqrt(5)/(4+sqrt(6))},{2*sqrt(5)/(4+sqrt(6))}) arc(126.8+73:142+72.6:5.85);
			
			\draw({-2*sqrt(5)/(4+sqrt(6))},{2*sqrt(5)/(4+sqrt(6))}) arc(61.6:41.8:8.75);			
			
			\draw(0,2) arc(178:180.5:60.44); 
			
			\draw(0,2) arc(92:-36:2);
			
			\node at (0,2.2){\scriptsize $A$}; \node at (-0.5,-0.15){\scriptsize $E$}; \node at (0.4,-0.15){\scriptsize $F$}; \node at (-0.25,-0.65){\scriptsize $B$}; \node at (1.7,-1.35){\scriptsize $C$};  \node at (-0.9,0.7){\scriptsize $D$};
			
		\end{scope}
		\begin{scope}[xshift=-5 cm, yshift=-4 cm]
			
			\draw[dotted] ({-2*sin(180/8)/(1+cos(180/8))},0) arc(180:157.3:{2*sin(180/8)/(1-cos(180/8)*cos(180/8))});
			\draw[dotted]
			({2*sin(180/8)/(1+cos(180/8))},0) arc(0:22.7:{2*sin(180/8)/(1-cos(180/8)*cos(180/8))});
			
			\fill (0,2) circle (0.05);  \fill ({2*sin(180/8)/(1+cos(180/8))},0) circle (0.05); 
			\fill ({-2*sin(180/8)/(1+cos(180/8))},0) circle (0.05); 
			
			\fill (-0.515,0.273) circle (0.05); 
			
			\fill (1.055,0.346) circle (0.05); 
			
			\fill ({-2/(4+sqrt(14))},{-2/(4+sqrt(14))}) circle (0.05); 
			
			\draw({-2/(4+sqrt(14))},{-2/(4+sqrt(14))}) arc(155.5+56:144+56:3.053);
			
			\draw({-2/(4+sqrt(14))},{-2/(4+sqrt(14))}) arc(64.6-137:79-137:5.862);
			
			\draw(0,2) arc(152.6-2:176:4.077); 
			
			\draw(0,2) arc(56.8:9.5:2.39); 
			
			\node at (0,2.2){\scriptsize $A$}; \node at (-0.5,-0.15){\scriptsize $E$}; \node at (0.4,-0.15){\scriptsize $F$}; \node at (-0.72,0.28){\scriptsize $B$}; \node at (1.25,0.28){\scriptsize $C$};  \node at (-0.25,-0.45){\scriptsize $D$};

		\end{scope}	
		\begin{scope}[xshift=-1.8cm,yshift=-4 cm] 	
			\draw[dotted] ({-2*sin(180/8)/(1+cos(180/8))},0) arc(180:157.3:{2*sin(180/8)/(1-cos(180/8)*cos(180/8))});
			\draw[dotted]
			({2*sin(180/8)/(1+cos(180/8))},0) arc(0:22.7:{2*sin(180/8)/(1-cos(180/8)*cos(180/8))});
			
			\fill (0,2) circle (0.05);  \fill ({2*sin(180/8)/(1+cos(180/8))},0) circle (0.05); 
			\fill ({-2*sin(180/8)/(1+cos(180/8))},0) circle (0.05);
			
			\fill ({(-2)*sqrt(350)/20/(1-sqrt(2)/20)},{(-2)*sqrt(48)/20/(1-sqrt(2)/20)}) circle (0.05);
			
			
			\fill (0.894,0.430) circle (0.05);
			
			\fill (4.280,2.421) circle (0.05); 
			
			\fill (1.77,0.64) circle (0.05); 
			
			
			\draw({(-2)*sqrt(350)/20/(1-sqrt(2)/20)},{(-2)*sqrt(48)/20/(1-sqrt(2)/20)}) arc(56+62.2:102.1:14.42);
			
			\draw[dotted] ({(-2)*sqrt(350)/20/(1-sqrt(2)/20)},{(-2)*sqrt(48)/20/(1-sqrt(2)/20)}) arc(56+62.2:100:14.42);
			
			\draw({(-2)*sqrt(350)/20/(1-sqrt(2)/20)},{(-2)*sqrt(48)/20/(1-sqrt(2)/20)}) arc(63.7+218:93+218.6:13.66);
			
			\draw(0,2) arc(50:8.4:2.61); 
			
			\draw(0,2) arc(74+65.7:52:3.09);
			
			\node at (0,2.2){\scriptsize $A$}; \node at (-0.4,-0.15){\scriptsize $E$}; \node at (0.4,-0.15){\scriptsize $F$}; \node at (0.65,0.55){\scriptsize $B$}; \node at (4.45,2.41){\scriptsize $C$};  \node at (-2.2,-0.8){\scriptsize $D$}; \node at (1.7,0.3){\scriptsize $D'$};

		\end{scope}
	\end{tikzpicture}
	\caption{Quadrilaterals constructed by $4$ typical positions of $D$. }\label{333dQuad}
\end{figure}

The symmetry of the $2$-layer earth map tiling $T(2n \, \bbb\ccc\ddd, 2\aaa^n)$ in Fig. \ref{bcdpic} implies that all $\aaa$-angles are assembled around north/south poles and the $2n$ middle points of all $b$-edges and $c$-edges distribute evenly on the equator with spacing $\frac{\pi}{n}$. This suggests the following geometric construction illustrated in Fig. \ref{333dQuad}. Fix a point $A$ on the sphere as the north pole, and take two points $E,F$ on the equator (i.e. $AE = AF = \frac{\pi}{2}$)  with $EF = \frac{\pi}{n}$.   the quadrilateral is then determined by the location of $D$ as follows: Extend $DE$ to $B$, such that $E$ is the middle point of $DB$. Extend $DF$ to $C$, such that $F$ is the middle point of $DC$. Then connect $A$ to $B$, $C$ to form the quadrilateral $\square ABDC$. Fig. \ref{333dQuad} shows four typical positions of $D$. We use the stereographic projection from the antipode of the middle point of $EF$ in both Fig. \ref{333dQuad} and \ref{333dModuli}. 

Thus the moduli is the possible locations of $D$, such that the boundary of $\square ABDC$ has no self intersection. In Fig. \ref{333dModuli}, we denote the south pole by $A'$. Extend $FE$ to $P$, such that $PE=\frac{\pi}{2}$.  Extend $EF$ to $Q$, such that $FQ=\frac{\pi}{2}$. Then we get the triangle $\triangle A'PQ$ with $PEFQ$ as one edge.  

\begin{theorem}
	The boundary of the quadrilateral $\square ABDC$ has no self intersection, if and only if $D$ lies in the interior of $\triangle AEF \cup \triangle A'PQ$  in Fig. \ref{333dModuli}, which describes the moduli of $2$-layer earth map tilings $T(2n \, \bbb\ccc\ddd, 2\aaa^n)$ for any $n\ge3$. Furthermore, $\square ABDC$ degenerates to a triangle if and only if $D$ lies in the interior of $EF$ $(\ddd=\pi)$, or $A'E$ $ (\bbb=\pi)$, or $A'F$ $ (\ccc=\pi)$. 
\end{theorem}

\begin{proof}
	When $D$ is in the northern hemisphere, two pictures in the first row of Fig. \ref{333dQuad} shows that the boundary of $\square ABDC$ has no self intersection if and only if $D$ lies in the interior of $\triangle AEF$. It is concave with $\ddd>\pi$, as shown in the second picture of Fig. \ref{real' figure}. 
	
	When $D$ is on the equator, the boundary of $\square ABDC$ has no self intersection if and only if $D$ lies in the interior of $EF$, and it degenerates to a triangle with $\ddd=\pi$. 
	
	When $D$ is in the southern hemisphere: $\square ABDC$ is simple and convex if and only if $D$ lies in the interior of $\triangle A'EF$ (shown in Fig. \ref{333dModuli} and the third picture of Fig. \ref{333dQuad}); $\square ABDC$ degenerates to a simple triangle if and only if $D$ lies in the interior of $A'E$ $ (\bbb=\pi)$ or $A'F$ $ (\ccc=\pi)$; $\square ABDC$ is simple and concave with $\bbb>\pi$ if and only if $D$ lies in the interior of $\triangle A'EP$; and symmetrically $\square ABDC$ is simple and concave with $\ccc>\pi$ if and only if $D$ lies in the interior of $\triangle A'FQ$. We will prove the case of $\triangle A'EP$, as shown in the fourth picture of Fig. \ref{333dQuad}, then the other cases follow easily. The key fact is that any two great arcs ($<2\pi$) starting from $D$ either intersect at its antipode $D'$ or never intersect. 
	
	When $D$ is on the left of the longitude $AEA'$, we have $DF>DE$ and $\bbb>\pi$. 	
	If $DE<\frac{\pi}{2}$, then $DB<\pi$ and it is too short to reach $D'$. So $DB$ does not intersect $DC$. 	
	If $DE\ge\frac{\pi}{2}$, then $DB\ge\pi$ and $DC=2DF>2DE\ge\pi$. So $DB$  meets $DC$ at $D'$. 	
	All such $D$ satisfying $DE=\frac{\pi}{2}$ form the great arc $A'P$. So $\square ABDC$ is simple and concave with $\bbb>\pi$ if and only if $D$ lies in the interior of $\triangle A'EP$.
\end{proof}

\begin{figure}[htp]
\centering
\begin{tikzpicture} [>=latex,scale=0.8] 
	
	\draw[color=gray!40] (0,0) circle (2); 
	\fill (0,-2) circle (1pt);
	\draw[dotted,color=gray!40] (0,0) ellipse (2 and 0.8);
	\draw[color=gray!40] (2,0) arc (0:-180:2 and 0.8);
	\fill (0,2) circle (1pt);
	\fill (-0.74,-0.74) circle (1pt);
	\fill (0.74,-0.74) circle (1pt);
	\fill (-1.16,0.64) circle (1pt);
	\fill (1.04,0.68) circle (1pt);
	
	
	\fill (-0.3,-1.2) circle (1pt);  \fill (-1.13,-0.08) circle (1pt);
	\fill (1.45,-0.08) circle (1pt);
	
	\draw (-0.3,-1.2) arc(230:203:3);
	\draw (-0.3,-1.2) arc(282.3:323:3);
	\draw (-1.13,-0.08) arc(174:129:3.1);
	\draw (1.45,-0.08) arc(10.5:59:3.1);
	
	\node at (-0.3,-1.5){\small $D$}; \node at (-1.35,-0.08){\small $B$}; \node at (1.7,-0.08){\small $C$};
	
	\node at (0,2.2){\small $A$}; \node at (0,-2.2){\small $A'$};
	\node at (-0.9,-1){\small $E$}; \node at (0.75,-1){\small $F$};
	\node at (-1.3,0.8){\small $P$}; \node at (0.7,0.7){\small $Q$};

	\begin{scope}[xshift=7 cm] 
		\draw[gray!50] (0,0) circle (2);

		\draw[dotted]  plot[smooth,samples=893,domain=-60:60] ({2*(sqrt(4*cos(\x)*cos(\x)-1)/(2*cos(\x)*cos(\x)))*sin(\x-180/6)/(1+sqrt(4*cos(\x)*cos(\x)-1)/(2*cos(\x)*cos(\x))*cos(\x-180/6))},{2*	(1-2*cos(\x)*cos(\x))/(2*cos(\x)*cos(\x))/(1+(sqrt(4*cos(\x)*cos(\x)-1)/(2*cos(\x)*cos(\x)))*cos(\x-180/6))});
		
		\draw[dotted]  plot[smooth,samples=893,domain=-60:60] ({-2*(sqrt(4*cos(\x)*cos(\x)-1)/(2*cos(\x)*cos(\x)))*sin(\x-180/6)/(1+sqrt(4*cos(\x)*cos(\x)-1)/(2*cos(\x)*cos(\x))*cos(\x-180/6))},{2*	(1-2*cos(\x)*cos(\x))/(2*cos(\x)*cos(\x))/(1+(sqrt(4*cos(\x)*cos(\x)-1)/(2*cos(\x)*cos(\x)))*cos(\x-180/6))});
		
		\draw[dotted]  (0,-2)--(0,2);

		\draw ({-2*sin(180/6)/(1+cos(180/6))},0) arc(180:150:{2*sin(180/6)/(1-cos(180/6)*cos(180/6))});
		\draw ({2*sin(180/6)/(1+cos(180/6))},0) arc(0:30:{2*sin(180/6)/(1-cos(180/6)*cos(180/6))});
		
		\draw({-2*cos(180/6)/(1-sin(180/6))},0) arc(180:300:{2*cos(180/6)/(1-sin(180/6)*sin(180/6))});
		\draw ({2*cos(180/6)/(1-sin(180/6))},0) arc(0:-120:{2*cos(180/6)/(1-sin(180/6)*sin(180/6))});
		
		\draw({-2*cos(180/6)/(1-sin(180/6))},0)--({-2*sin(180/6)/(1+cos(180/6))},0)
		({2*sin(180/6)/(1+cos(180/6))},0)--({2*cos(180/6)/(1-sin(180/6))},0);

		\node[rotate=90] at (-0.2,-1.4) {\scriptsize $b=c$};
		
		\node[rotate=10] at (-1,2.2) {\scriptsize $a=b$};
		\node[rotate=-10] at (1,2.2) {\scriptsize $a=c$};
		
		\node at (-0.05,-2.25) {\scriptsize $A'$}; \node at (0,2.2) {\scriptsize $A$};
		\node at (-0.8,0.2) {\scriptsize $E$}; \node at (0.8,0.2) {\scriptsize $F$};
		\node at (-3.3,0.2) {\scriptsize $P$}; \node at (3.3,0.2) {\scriptsize $Q$};
		
		\node at (0,-2.8) {\small $f=6$};
		
	\end{scope}
	\begin{scope}[yshift=-6 cm] 
		\draw[gray!50] (0,0) circle (2);

		\draw[dotted]  plot[smooth,samples=781,domain=-60:60] ({2*(sqrt(4*cos(\x)*cos(\x)-1)/(2*cos(\x)*cos(\x)))*sin(\x-180/8)/(1+sqrt(4*cos(\x)*cos(\x)-1)/(2*cos(\x)*cos(\x))*cos(\x-180/8))},{2*	(1-2*cos(\x)*cos(\x))/(2*cos(\x)*cos(\x))/(1+(sqrt(4*cos(\x)*cos(\x)-1)/(2*cos(\x)*cos(\x)))*cos(\x-180/8))});
		
		\draw[dotted]  plot[smooth,samples=781,domain=-60:60] ({-2*(sqrt(4*cos(\x)*cos(\x)-1)/(2*cos(\x)*cos(\x)))*sin(\x-180/8)/(1+sqrt(4*cos(\x)*cos(\x)-1)/(2*cos(\x)*cos(\x))*cos(\x-180/8))},{2*	(1-2*cos(\x)*cos(\x))/(2*cos(\x)*cos(\x))/(1+(sqrt(4*cos(\x)*cos(\x)-1)/(2*cos(\x)*cos(\x)))*cos(\x-180/8))});
		
		\draw[dotted]  (0,-2)--(0,2);

		\draw ({-2*sin(180/8)/(1+cos(180/8))},0) arc(180:157.3:{2*sin(180/8)/(1-cos(180/8)*cos(180/8))});
		\draw ({2*sin(180/8)/(1+cos(180/8))},0) arc(0:22.7:{2*sin(180/8)/(1-cos(180/8)*cos(180/8))});          	
		
		\draw({-2*cos(180/8)/(1-sin(180/8))},0) arc(180:292.4:{2*cos(180/8)/(1-sin(180/8)*sin(180/8))});
		\draw ({2*cos(180/8)/(1-sin(180/8))},0) arc(0:-112.4:{2*cos(180/8)/(1-sin(180/8)*sin(180/8))});
		
		\draw({-2*cos(180/8)/(1-sin(180/8))},0)--({-2*sin(180/8)/(1+cos(180/8))},0)
		({2*sin(180/8)/(1+cos(180/8))},0)--({2*cos(180/8)/(1-sin(180/8))},0);

		\node[rotate=90] at (-0.2,-1.4) {\scriptsize $b=c$};
		
		\node[rotate=10] at (-1,2.2) {\scriptsize $a=b$};
		\node[rotate=-10] at (1,2.2) {\scriptsize $a=c$};
		\node at (0,-2.8) {\small $f=8$};

	\end{scope}
	\begin{scope}[xshift=7 cm,yshift=-6 cm]

		\draw[gray!50] (0,0) circle (2);
		
		\draw[dotted]  plot[smooth,samples=777+4,domain=-60:60] ({2*(sqrt(4*cos(\x)*cos(\x)-1)/(2*cos(\x)*cos(\x)))*sin(\x-180/10)/(1+sqrt(4*cos(\x)*cos(\x)-1)/(2*cos(\x)*cos(\x))*cos(\x-180/10))},{2*	(1-2*cos(\x)*cos(\x))/(2*cos(\x)*cos(\x))/(1+(sqrt(4*cos(\x)*cos(\x)-1)/(2*cos(\x)*cos(\x)))*cos(\x-180/10))});	    
		
		\draw[dotted]  plot[smooth,samples=777+4,domain=-60:60] ({-2*(sqrt(4*cos(\x)*cos(\x)-1)/(2*cos(\x)*cos(\x)))*sin(\x-180/10)/(1+sqrt(4*cos(\x)*cos(\x)-1)/(2*cos(\x)*cos(\x))*cos(\x-180/10))},{2*	(1-2*cos(\x)*cos(\x))/(2*cos(\x)*cos(\x))/(1+(sqrt(4*cos(\x)*cos(\x)-1)/(2*cos(\x)*cos(\x)))*cos(\x-180/10))});	
	
		\draw[dotted]  (0,-2)--(0,2);	    		    		    	
		\draw ({-2*sin(180/10)/(1+cos(180/10))},0) arc(180:162:{2*sin(180/10)/(1-cos(180/10)*cos(180/10))});
		\draw ({2*sin(180/10)/(1+cos(180/10))},0) arc(0:18:{2*sin(180/10)/(1-cos(180/10)*cos(180/10))});

		\draw({-2*cos(180/10)/(1-sin(180/10))},0) arc(180:288:{2*cos(180/10)/(1-sin(180/10)*sin(180/10))});
		\draw ({2*cos(180/10)/(1-sin(180/10))},0) arc(0:-108:{2*cos(180/10)/(1-sin(180/10)*sin(180/10))});
		
		\draw({-2*cos(180/10)/(1-sin(180/10))},0)--({-2*sin(180/10)/(1+cos(180/10))},0)
		({2*sin(180/10)/(1+cos(180/10))},0)--({2*cos(180/10)/(1-sin(180/10))},0);
		
		\node[rotate=90] at (-0.2,-1.4) {\scriptsize $b=c$};
		
		\node[rotate=10] at (-1,2.2) {\scriptsize $a=b$};	    	
		\node[rotate=-10] at (1,2.2) {\scriptsize $a=c$};		
		\node at (0,-2.8) {\small $f\ge10$};

	\end{scope}			

\end{tikzpicture}	
\caption{The moduli  $(\triangle AEF \cup \triangle A'P_{EF}Q)^\circ$ and its projection. }\label{333dModuli}
\end{figure}
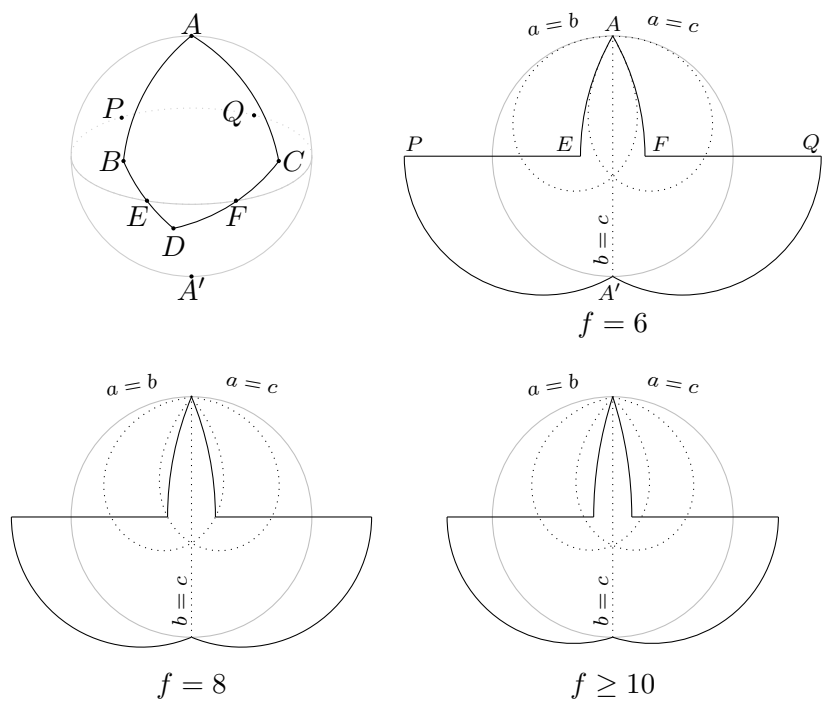	
\begin{figure}[htp]
	\centering
	
	\includegraphics[scale=0.22]{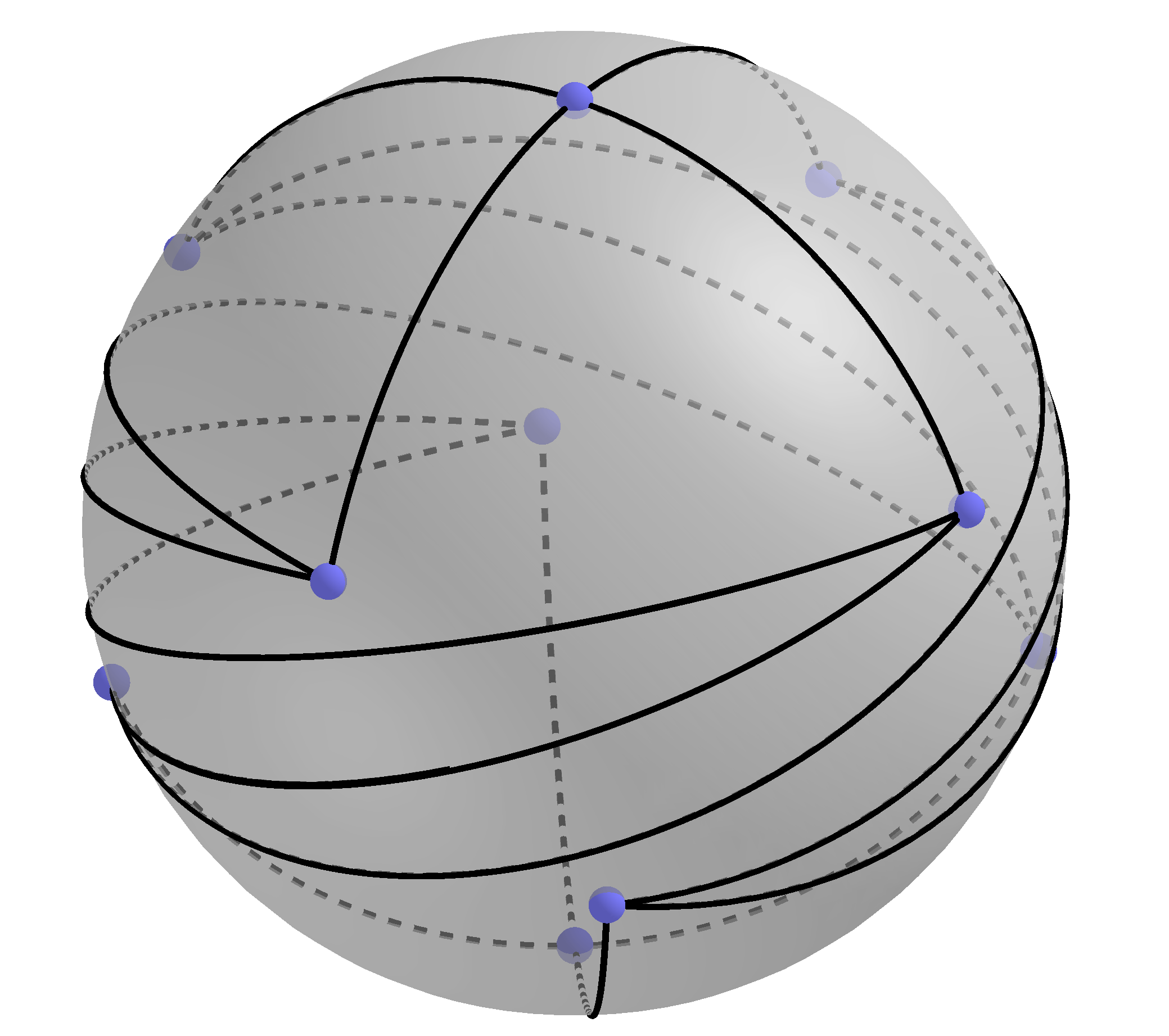}\hspace{60pt}
	\includegraphics[scale=0.19]{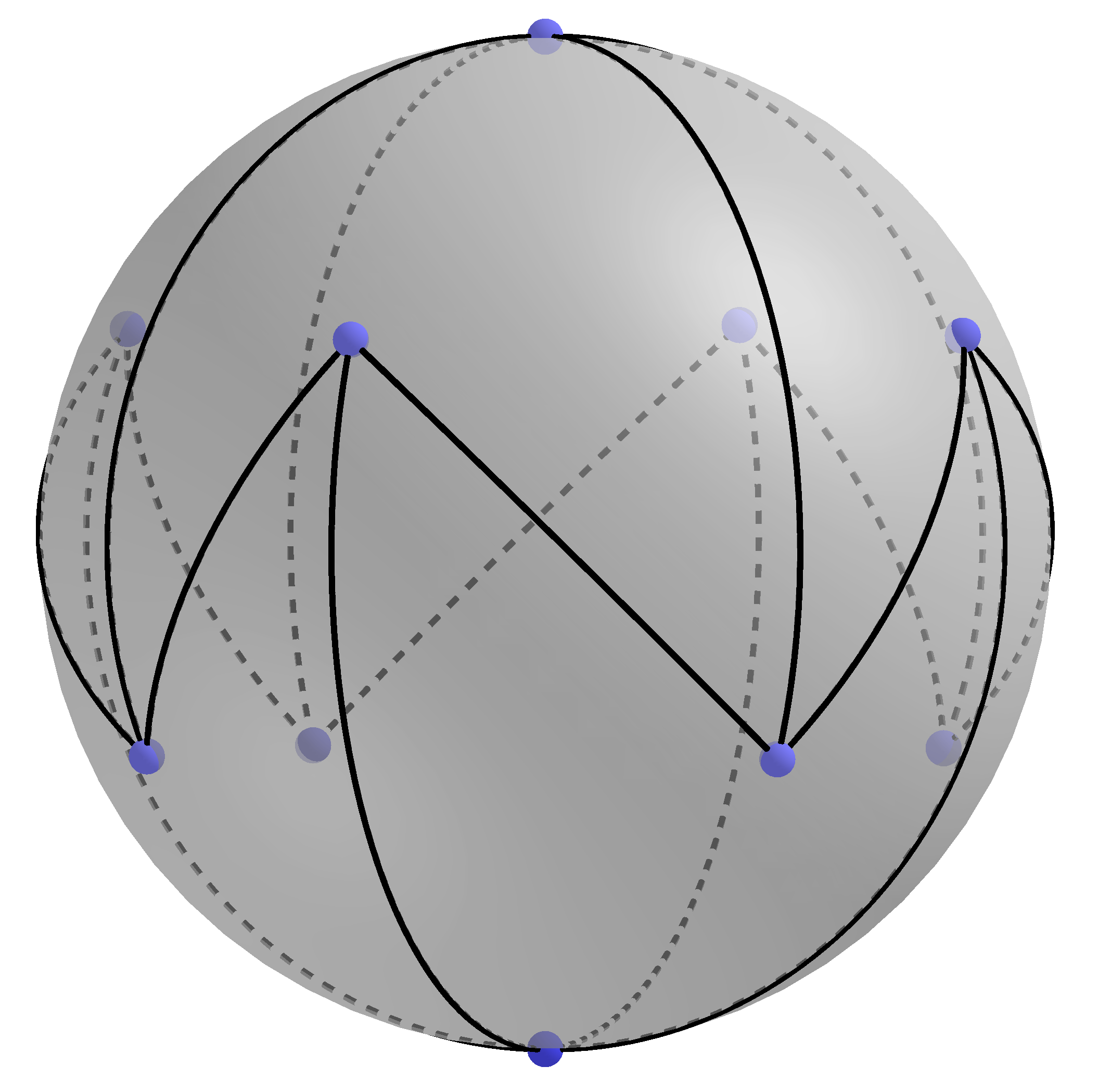}
	
	\caption{Concave $2$-layer earth map tilings drawn by GeoGebra.}
	\label{real' figure}
\end{figure}


The second picture of Fig. \ref{333dModuli} shows the moduli of $T(6\bbb\ccc\ddd, 2\aaa^3)$,  and the dotted curves inside the moduli represent reductions of the quadrilateral from Type $a^2 bc$ to Type $a^2 b^2$ ($b=c$), Type  $a^3 b$ ($a=b$ or $a=c$), and Type $a^4$ ($a=b=c$). The third and fourth pictures of Fig. \ref{333dModuli} are for $f=8$ and $f\ge 10$ respectively, where the reduction curves have different positions inside the moduli. In the next two papers \cite{lw,lpwx} of this series, it turns out that most Type $a^3 b$ quadrilateral tilings of the sphere come from these $2$-layer earth map tilings on the reduction curves together with certain modifications under  extra conditions. Thus the detailed study of the reduction curves will be shown in \cite{lpwx}. 

We remark that the last picture in Fig. \ref{333dQuad} could produce strange quadrilaterals with $c>\pi$, as shown in the first picture of Fig. \ref{real' figure}.

\section{$334d$-Tile and $335d$-Tile}
\label{3345d}

After Proposition \ref{333d}, we will always assume that $\delta$ never appears in any degree $3$ vertex. 

\begin{lemma}\label{lemd4}
	In an $a^2bc$-tiling, if $\bbb\ccc\ddd$ is not a vertex (i.e. $\ddd$ never appears in degree $3$ vertices), then all degree $3$ vertices are $\aaa^3$, $\aaa\bbb^2$, or $\aaa\ccc^2$. Furthermore, there always exist vertices $\aaa^2\cdots$, $\bbb^2\cdots$, $\ccc^2\cdots$, and $\ddd^2\cdots$. In particular, the quadrilateral is convex with all angles being $<\pi$ and $\bbb\neq\ccc$. 
\end{lemma}
\begin{proof}
	If $\bbb\ccc\ddd$ is not a vertex, Lemma \ref{lem3} implies the first statement directly, and Balance Lemma \ref{balance} implies that $\bbb^2\cdots$, $\ccc^2\cdots$, and $\ddd^2\cdots$ are all vertices. We just need to show $\aaa^2\cdots$ is a vertex. Otherwise all degree $3$ vertices ($v_3\ge8$) could only be one of $\aaa\bbb^2, \aaa\ccc^2$ by Proposition \ref{lembc3}. Take $\aaa\ccc^2$, then Lemma \ref{lemac2} and no $\aaa^2\cdots$ imply that some vertex $\aaa\bbb^{2t}$ ($t\ge2$) must appear. But its AAD $\bbb^\alpha\thin^\alpha\bbb\cdots$ implies a vertex $\aaa^2\cdots$, a contradiction. So all angles are $<\pi$. Then we have $\bbb\neq\ccc$ by Lemma \ref{geometry3}.   
\end{proof}

\begin{figure}[htp]
	\centering
	\begin{tikzpicture}[>=latex,scale=0.6]      
	\draw (0,0) -- (0,2) 
	(0,0) -- (2,0);
	\draw[dashed]  (0,2)--(2,2);
	\draw[line width=1.5] (2,0)--(2,2);
	\node at (0.35,0.35){\small $\aaa$};
	\node at (1.65,0.35){\small $\bbb$};
	\node at (1.65,1.65){\small $\ddd$};
	\node at (0.35,1.65){\small $\ccc$};
	\node at (0,-0.45) {\small $3$};          
	\node at (0,2.45) {\small $3$};
	\node at (2,-0.45) {\small $d$};          
	\node at (2,2.45) {\small $4$};
	\fill (2,0) circle (0.2);
	\node at (1,-1.6) {\small $4\le d\le11$};
	\begin{scope}[xshift=4.5 cm]      
	\draw (0,0) -- (0,2) 
	(0,0) -- (2,0);
	\draw[dashed]  (0,2)--(2,2);
	\draw[line width=1.5] (2,0)--(2,2);
	\node at (0.35,0.35){\small $\aaa$};
	\node at (1.65,0.35){\small $\bbb$};
	\node at (1.65,1.65){\small $\ddd$};
	\node at (0.35,1.65){\small $\ccc$};
	\node at (0,-0.45) {\small $3$};          
	\node at (0,2.45) {\small $3$};
	\node at (2,-0.45) {\small $4$};          
	\node at (2,2.45) {\small $d$};
	\fill (2,2) circle (0.2);
	\node at (1,-1.6) {\small $5\le d\le11$};
	\end{scope}  

	\begin{scope}[xshift=9 cm]       
	\draw (0,0) -- (0,2) 
	(0,0) -- (2,0);
	\draw[dashed]  (0,2)--(2,2);
	\draw[line width=1.5] (2,0)--(2,2);
	\node at (0.35,0.35){\small $\aaa$};
	\node at (1.65,0.35){\small $\bbb$};
	\node at (1.65,1.65){\small $\ddd$};
	\node at (0.35,1.65){\small $\ccc$};
	\node at (0,-0.45) {\small $3$};          
	\node at (0,2.45) {\small $3$};
	\node at (2,-0.45) {\small $d$};          
	\node at (2,2.45) {\small $5$};
	\fill (2,0) circle (0.2);
	\node at (1,-1.6) {\small $d=5,6,7$};
	\end{scope}
	\begin{scope}[xshift=13.5 cm]   
	\draw (0,0) -- (0,2) 
	(0,0) -- (2,0);
	\draw[dashed]  (0,2)--(2,2);
	\draw[line width=1.5] (2,0)--(2,2);
	\node at (0.35,0.35){\small $\aaa$};
	\node at (1.65,0.35){\small $\bbb$};
	\node at (1.65,1.65){\small $\ddd$};
	\node at (0.35,1.65){\small $\ccc$};
	\node at (0,-0.45) {\small $3$};          
	\node at (0,2.45) {\small $3$};
	\node at (2,-0.45) {\small $5$};          
	\node at (2,2.45) {\small $d$};
	\fill (2,2) circle (0.2);
	\node at (1,-1.6) {\small $d=6,7$};
	\end{scope}
    \begin{scope}[xshift=19.8 cm,scale=0.8,yshift=-0.4 cm]       
    		\draw (0,0) -- (0,2) 
    		(0,0) -- (2,0)
    		(0,2)--(-2,4)
    		(-2,4)--(2,4)
    		(2,0)--(2,-2);
    		\draw[dashed]  (0,2)--(2,2)
    		(0,0)--(-2,-2);
    		\draw[line width=1.5] (2,0)--(2,2)
    		(2,2)--(2,4)
    		(-2,4)--(-2,-2)
    		(-2,-2)--(2,-2);

    		\node at (0.35,0.35) {\small $\aaa$};
    		\node at (0.2,-0.45) {\small $\ccc$};
    		\node at (-0.35,0.25) {\small $\ccc$};
    		
    		\node at (1.65,0.45) {\small $\bbb$};
    		\node at (1.65,-0.35) {\small $\aaa$};
    		
    		\node at (1.65,1.55) {\small $\ddd$};
    		\node at (1.65,2.45) {\small $\ddd$};
    		
    		\node at (1.65,3.5) {\small $\bbb$};
    		
    		\node at (0,2.4) {\small $\ccc$};
    		\node at (0.35,1.55) {\small $\ccc$};
    		\node at (-0.4,1.75) {\small $\aaa$};
    		
    		\node at (-1.2,3.7) {\small $\aaa$};
    		\node at (-1.65,3) {\small $\bbb$};
    		\node at (-1.1,-1.6) {\small $\ddd$};
    		\node at (-1.7,-1.2) {\small $\ddd$};
    		\node at (1.65,-1.5) {\small $\bbb$};
    		
    		\node[draw,shape=circle, inner sep=0.5] at (1,1) {\small $1$};
    		\node[draw,shape=circle, inner sep=0.5] at (0.8,3) {\small $2$};
    		\node[draw,shape=circle, inner sep=0.5] at (1,-1) {\small $4$};
    		\node[draw,shape=circle, inner sep=0.5] at (-1,1) {\small $3$};
    	\end{scope}	
    	\end{tikzpicture}	
    	\caption{Special $334d$, $335d$-Tiles and their common partial neighborhood} \label{33--}
    \end{figure}
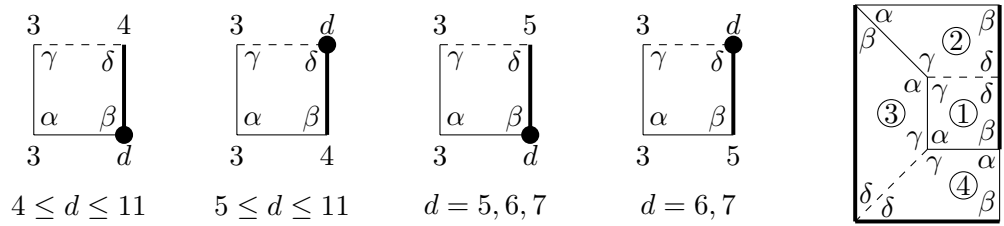

By Proposition \ref{lembc3} and the symmetry of exchanging $\bbb\leftrightarrow\ccc$, there are only $4$ different configurations for degree $3$ vertices and degree $d$ vertex $H$ (indicated by $\bullet$) in $334d$-Tile and $335d$-Tile in Fig. \ref{33--}. Recall that $4\le d \le11$ for Case $334d$, and $ d=5,6,7$ for Case $335d$.

The fifth picture of Fig. \ref{33--} shows the common partial neighborhood in $334d$-Tile and $335d$-Tile:  
Firstly $\ccc_1\cdots=\aaa_3\ccc_1\ccc_2$ determines $T_2$. By $\aaa_3$ and Proposition \ref{lembc3}, we get $E_{34}=c$ and $\aaa_1\cdots=\aaa\ccc^2$. This determines $T_3,T_4$.

\begin{proposition} \label{334dCase1}
 All $a^2bc$-tilings with the 1st special tile in Fig. \ref{33--} are 
 \begin{enumerate}[$(1)$]
 \item
 The flip modification of a unique quadrilateral subdivision of the octahedron   $T(2\aaa^3,6\aaa\ccc^2,6\ddd^4,6\bbb^2\ccc^2,6\aaa^2\bbb^2)$ with $24$ tiles;
 \item
 A sequence of $3$-layer earth map tilings (each has a unique quadrilateral)
 $T(4n\,\aaa\ccc^2,2\bbb^{2n},2n\,\ddd^4,2n\,\aaa^2\bbb^2)$ with $8n$ tiles for any $n\ge2$, among which each odd $n=2m+1$ case admits exactly two flip modifications:
 \begin{itemize}
 	\item $T((8m+4)\aaa\ccc^2,4\aaa\bbb^{2m+2},(4m+2)\ddd^4,4m\,\aaa^2\bbb^2)$;
 	\item $T((8m+2)\aaa\ccc^2,2\aaa\bbb^{2m+2},2\bbb^{2m}\ccc^2,(4m+2)\ddd^4,(4m+2)\aaa^2\bbb^2)$.
 \end{itemize}
 \end{enumerate}
\end{proposition}
\begin{proof}
Let the first of Fig. \ref{33--} be the center tile $T_1$ in the partial neighborhoods in Fig.  \ref{334d1}.  If $E_{56}=a$ in the first picture of Fig.  \ref{334d1}, then $T_5,T_6$ are determined and $H=\aaa\bbb\ddd\cdots$, contradicting  
Lemma \ref{lemac2}. So we have $E_{56}=c$ in the second picture, which determines $T_5,T_6$. 

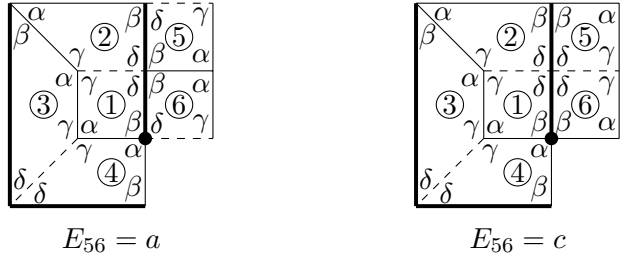
\begin{figure}[htp]
	\centering	
	\begin{tikzpicture}[>=latex,scale=0.45]       
	\draw (0,0) -- (0,2) 
	(0,0) -- (2,0)
	(0,2)--(-2,4)
	(-2,4)--(2,4)
	(2,0)--(2,-2)
	(4,2)--(4,4)
	(2,4)--(4,4)
	(4,2)--(4,0)
	(2,0)--(4,0);
	\draw[dashed]  (0,2)--(2,2)
	(0,0)--(-2,-2)
	(2,2)--(4,2);
	\draw[line width=1.5] (2,0)--(2,2)
	(2,2)--(2,4)
	(-2,4)--(-2,-2)
	(-2,-2)--(2,-2);
	(4,-1)--(3,-2)
	(3,-2)--(2,-2);

	\fill (2,0) circle (0.2);
	
	\node at (0.35,0.35) {\small $\aaa$};
	\node at (0.2,-0.45) {\small $\ccc$};
	\node at (-0.35,0.25) {\small $\ccc$};
	
	\node at (1.65,0.45) {\small $\bbb$};
	\node at (1.65,-0.35) {\small $\aaa$};
	
	\node at (1.65,1.55) {\small $\ddd$};
	\node at (1.65,2.45) {\small $\ddd$};
	
	\node at (1.65,3.5) {\small $\bbb$};
	
	\node at (0,2.4) {\small $\ccc$};
	\node at (0.35,1.55) {\small $\ccc$};
	\node at (-0.4,1.75) {\small $\aaa$};
	
	\node at (-1.2,3.7) {\small $\aaa$};
	\node at (-1.65,3) {\small $\bbb$};
	\node at (-1.1,-1.6) {\small $\ddd$};
	\node at (-1.7,-1.2) {\small $\ddd$};
	\node at (1.65,-1.5) {\small $\bbb$};
	
	\node at (2.35,1.55) {\small $\ddd$};
	\node at (2.35,2.45) {\small $\ddd$};
	\node at (2.35,3.5) {\small $\bbb$};
	
	\node at (3.65,1.55) {\small $\ccc$};
	\node at (3.65,2.45) {\small $\ccc$};
	\node at (3.65,3.6) {\small $\aaa$};
	
	\node at (2.35,0.45) {\small $\bbb$}; \node at (3.65,0.45) {\small $\aaa$}; 
	
	\node[draw,shape=circle, inner sep=0.5] at (1,1) {\small $1$};
	\node[draw,shape=circle, inner sep=0.5] at (0.8,3) {\small $2$};
	\node[draw,shape=circle, inner sep=0.5] at (3,3) {\small $5$};
	\node[draw,shape=circle, inner sep=0.5] at (3,1) {\small $6$};
	\node[draw,shape=circle, inner sep=0.5] at (1,-1) {\small $4$};
	\node[draw,shape=circle, inner sep=0.5] at (-1,1) {\small $3$};
	
	\node at (1,-3) {\small $E_{56}=c$};
	
	\begin{scope}[xshift=-12 cm]
		\draw (0,0) -- (0,2) 
		(0,0) -- (2,0)
		(0,2)--(-2,4)
		(-2,4)--(2,4)
		(2,0)--(2,-2)				
		(4,2)--(4,0)		
		(2,2)--(4,2)
		(4,2)--(4,4);
		\draw[dashed]  (0,2)--(2,2)
		(0,0)--(-2,-2)
		(2,4)--(4,4)
		(2,0)--(4,0);
		\draw[line width=1.5] (2,0)--(2,2)
		(2,2)--(2,4)
		(-2,4)--(-2,-2)
		(-2,-2)--(2,-2);
		(4,-1)--(3,-2)
		(3,-2)--(2,-2);

		\fill (2,0) circle (0.2);
		
		\node at (0.35,0.35) {\small $\aaa$};
		\node at (0.2,-0.45) {\small $\ccc$};
		\node at (-0.35,0.25) {\small $\ccc$};
		
		\node at (1.65,0.45) {\small $\bbb$};
		\node at (1.65,-0.35) {\small $\aaa$};
		
		\node at (1.65,1.55) {\small $\ddd$};
		\node at (1.65,2.45) {\small $\ddd$};
		
		\node at (1.65,3.5) {\small $\bbb$};
		
		\node at (0,2.4) {\small $\ccc$};
		\node at (0.35,1.55) {\small $\ccc$};
		\node at (-0.4,1.75) {\small $\aaa$};
		
		\node at (-1.2,3.7) {\small $\aaa$};
		\node at (-1.65,3) {\small $\bbb$};
		\node at (-1.1,-1.6) {\small $\ddd$};
		\node at (-1.7,-1.2) {\small $\ddd$};
		\node at (1.65,-1.5) {\small $\bbb$};
		
		\node at (2.35,1.55) {\small $\bbb$};
		\node at (2.35,2.45) {\small $\bbb$};
		\node at (2.35,3.5) {\small $\ddd$};
		
		\node at (3.65,1.55) {\small $\aaa$};
		\node at (3.65,2.45) {\small $\aaa$};
		\node at (3.65,3.6) {\small $\ccc$};
		
		\node at (2.35,0.45) {\small $\ddd$}; \node at (3.65,0.45) {\small $\ccc$}; 
		\node at (1,-3) {\small $E_{56}=a$};
		
		\node[draw,shape=circle, inner sep=0.5] at (1,1) {\small $1$};
		\node[draw,shape=circle, inner sep=0.5] at (0.8,3) {\small $2$};
		\node[draw,shape=circle, inner sep=0.5] at (3,3) {\small $5$};
		\node[draw,shape=circle, inner sep=0.5] at (3,1) {\small $6$};
		\node[draw,shape=circle, inner sep=0.5] at (1,-1) {\small $4$};
		\node[draw,shape=circle, inner sep=0.5] at (-1,1) {\small $3$};
	\end{scope}
	\end{tikzpicture}	
	\caption{Partial neighborhoods of the 1st special tile.} \label{334d1}
\end{figure}

By $\aaa\ccc^2$ and Lemma \ref{lemac2}, $H=\aaa_4\bbb_1\bbb_6\cdots=\aaa^k\bbb^{2t}$ with $k\ge1,t\ge1$. If $k \ge 2$ and $t \ge 2$, we have $2\aaa+4\bbb \le 2\pi $. By $\aaa\ccc^2$ and $\ddd^4$, we get $4(\aaa+\bbb+\ccc+\ddd)=2(\aaa+2\ccc)+4\ddd+(2\aaa+4\bbb)\le8\pi$, contradicting Lemma \ref{anglesum}. Therefore $k$ or $t=1$, and $H=\aaa^{d-2}\bbb^2(4\le d\le11)$ or $\aaa\bbb^{d-1}(d=5,7,9,11)$. 

\subsection*{$\mathbf{H=\aaa^{d-2}\bbb^2\,(4\le d\le11)}$}

By $\aaa\ccc^2$, $\ddd^4$ and $\aaa^{d-2}\bbb^2$, we get
\begin{equation*}
	\aaa=\tfrac{\pi}{d-3}-\tfrac{8\pi}{(d-3)f},\,\, \bbb=\tfrac{\pi}{2}-\tfrac{\pi}{2(d-3)}+\tfrac{4(d-2)\pi}{(d-3)f},\,\,
	\ccc=\pi-\tfrac{\pi}{2(d-3)}+\tfrac{4\pi}{(d-3)f},\,\,
	\ddd=\tfrac{\pi}{2}.
\end{equation*}

If $d\ge5$, then $\aaa<\frac\pi2,\bbb>\frac\pi4,\ccc>\frac{3\pi}{4}$. By $d\ge5$, the AAD of $H=\aaa\aaa\aaa\cdots$ induces $\bbb\thin\ccc\cdots$ or $\ccc\thin\ccc\cdots$.  By Lemma $\ref {lembd}'$, we know $\ccc\thin\ccc\cdots$ is not a vertex. By $\bbb+\ccc>\pi$, $\bbb\ccc\cdots$ must be odd. But $R(\bbb\ccc\ddd\cdots)<2\bbb,2\ccc,2\ddd$, so $\bbb\ccc\cdots=\bbb\ccc\ddd$, contradicting Proposition \ref{333d}. 

We conclude $d=4$, $H=\aaa^2\bbb^2$ and $\aaa=(1-\tfrac{8}{f})\pi$, $ \bbb=\tfrac{8\pi}{f}$, $\ccc=(\tfrac{1}{2}+\tfrac{4}{f})\pi$. Now we show $f\ge16$. If $f<16$, we get $\aaa<\frac{\pi}{2},\bbb>\frac{\pi}{2},\ccc>\frac{3\pi}{4},\ddd=\frac{\pi}{2}$. This implies $\bbb^2\cdots=\aaa^2\bbb^2$. Then by the parity lemma, we further get $\aaa\bbb\cdots=\aaa\bbb^2\cdots=\aaa^2\bbb^2$.
Extend the second picture of Fig. \ref{334d1} to Fig. \ref{334d1.1} to show the complete neighborhood of the 1st special tile.
Then $\aaa_2\bbb_3\cdots=\aaa_2\aaa_7\bbb_3\bbb$. By $\aaa_7$,  $R(\bbb_2\bbb_5\cdots)$ has $\bbb$ or $\ccc$, contradicting $\bbb^2\cdots=\aaa^2\bbb^2$.  

\begin{figure}[htp]
	\centering	
	\begin{tikzpicture}[>=latex,scale=0.45]       
		\draw (0,0) -- (0,2) 
		(0,0) -- (2,0)
		(0,2)--(-2,4)
		(-2,4)--(2,4)
		(2,0)--(2,-2)
		(4,2)--(4,4)
		(2,4)--(4,4)
		(4,2)--(4,0)
		(2,0)--(4,0);
		\draw[dashed]  (0,2)--(2,2)
		(0,0)--(-2,-2)
		(2,2)--(4,2);
		\draw[line width=1.5] (2,0)--(2,2)
		(2,2)--(2,4)
		(-2,4)--(-2,-2)
		(-2,-2)--(2,-2);
		\draw[dotted] (-2,4)--(-3,5)
		(-3,5)--(2,5)
		(2,4)--(2,5)
		(2,-2)--(4,-2)
		(4,0)--(4,-2);

		\fill (2,0) circle (0.2);
		
		\node at (0.35,0.35) {\small $\aaa$};
		\node at (0.2,-0.45) {\small $\ccc$};
		\node at (-0.35,0.25) {\small $\ccc$};
		
		\node at (1.65,0.45) {\small $\bbb$};
		\node at (1.65,-0.35) {\small $\aaa$};
		
		\node at (1.65,1.55) {\small $\ddd$};
		\node at (1.65,2.45) {\small $\ddd$};
		
		\node at (1.65,3.5) {\small $\bbb$};
		
		\node at (0,2.4) {\small $\ccc$};
		\node at (0.35,1.55) {\small $\ccc$};
		\node at (-0.4,1.75) {\small $\aaa$};
		
		\node at (-1.2,3.7) {\small $\aaa$};
		\node at (-1.65,3) {\small $\bbb$};
		\node at (-1.1,-1.6) {\small $\ddd$};
		\node at (-1.7,-1.2) {\small $\ddd$};
		\node at (1.65,-1.5) {\small $\bbb$};
		
		\node at (2.35,1.55) {\small $\ddd$};
		\node at (2.35,2.45) {\small $\ddd$};
		\node at (2.35,3.5) {\small $\bbb$};
		
		\node at (3.65,1.55) {\small $\ccc$};
		\node at (3.65,2.45) {\small $\ccc$};
		\node at (3.65,3.6) {\small $\aaa$};
		
		\node at (2.35,0.45) {\small $\bbb$}; \node at (3.65,0.45) {\small $\aaa$}; 
		
		\node at (-1.7,4.35) {\small $\aaa$};
		\node at (-2.4,3.6) {\small $\bbb$};
		
		\node at (2.35,-0.35) {\small $\aaa$};
		
		\node[draw,shape=circle, inner sep=0.5] at (1,1) {\small $1$};
		\node[draw,shape=circle, inner sep=0.5] at (1,3) {\small $2$};
		\node[draw,shape=circle, inner sep=0.5] at (3,3) {\small $5$};
		\node[draw,shape=circle, inner sep=0.5] at (3,1) {\small $6$};
		\node[draw,shape=circle, inner sep=0.5] at (1,-1) {\small $4$};
		\node[draw,shape=circle, inner sep=0.5] at (-1,1) {\small $3$};
		\node[draw,shape=circle, inner sep=0.5] at (0,4.5) {\small $7$};
		
	\end{tikzpicture}	
	\caption{The neighborhood of the 1st special tile with $H=\aaa^2\bbb^2$.} \label{334d1.1}
\end{figure}
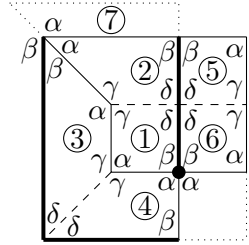

By $f\ge 16$, we have $\aaa\ge\frac{\pi}{2},\bbb \le \frac{\pi}{2},\frac{\pi}{2}<\ccc \le \frac{3\pi}{4},\ddd=\frac{\pi}{2}$. If $\aaa^k\bbb^l\ccc^m\ddd^n$ is a vertex, then we have $k \le 4,m \le 3,n\le4$, and
\[
(1-\tfrac{8}{f})k+\tfrac{8}{f}l+(\tfrac12+\tfrac{4}{f})m+\tfrac{1}{2}n=2
\] 
 We substitute
the finitely many combinations of exponents satisfying the bounds into the equation above and solve for $f$. By the angle values and the edge length consideration, we get all possible AVC in Table \ref{tab-1.1}. Its first row  ``$f=\text{all}$'' means that the angle combinations can be vertices for any $f$; all other rows are mutually exclusive.  All possible tilings based on the AVC of Table \ref{tab-1.1} are deduced as follows. 
\begin{table}[htp]                        
	\centering 
	\caption{The AVC for $H=\aaa^2\bbb^2$ and $f\ge16$ }\label{tab-1.1}                           
	\begin{tabular}{c|c}
		\hline  
		$f$&vertex\\
		\hline\hline
		all&$\aaa\ccc^2,\ddd^4,\aaa^2\bbb^2$ \\
		$16s-4,s=2,3,\cdots$&$\bbb^{2s-1}\ccc\ddd$\\
		$16$&$\aaa^4,\bbb^2\ddd^2,\bbb^4$\\
		$16s,s=2,3,\cdots$&$\bbb^{2s}\ddd^2,\bbb^{4s}$ \\
		$24$&$\aaa^3,\aaa\bbb^4,\bbb^2\ccc^2,\bbb^6$ \\				
		$16s+8,s=2,3,\cdots$&$\aaa\bbb^{2s+2},\bbb^{2s}\ccc^2,\bbb^{4s+2}$ \\
		\hline
	\end{tabular}  
	
\end{table}

If $\text{AVC}=\{\aaa\ccc^2,\ddd^4,\aaa^2\bbb^2$\}, there is no solution satisfying Balance Lemma.

\vspace{9pt}
{\bf{Claim}}: For any $a^2bc$-tiling with the AVC in Table \ref{tab-1.1}, if $\bbb^{2k}\,(k=2,3,\cdots)$ is a vertex, $\ddd\cdots=\ddd^4$, and $\aaa^3,\aaa^4$ are not vertices, then it must be the $3$-layer earth map tiling $T(4k\, \aaa\ccc^2, 2\bbb^{2k},2k\,\ddd^4,2k\,\aaa^2\bbb^2)$ with $8k$ tiles in Fig. \ref{case1jie3}.
\vspace{9pt}

In Fig. \ref{case1jie3}, $\bbb^{2k}=\thick\bbb_1\thin\bbb_2\thick\cdots$ determines $T_1,T_2$. Then $\aaa_1\aaa_2\cdots=\aaa_1\aaa_2\bbb_3\bbb_4$ determines $T_3,T_4$; $\aaa_3\ccc_1\cdots=\aaa_3\ccc_1\ccc_5$ determines $T_5$;  $\aaa_4\ccc_2\cdots=\aaa_4\ccc_2\ccc_6$ determines $T_6$;  $\ddd_3\ddd_4\cdots=\ddd^4$ determines $T_7,T_8$. 
The eight tiles $\{T_i\}_{1\le i \le8}$ together form a time zone. Similarly, we can determine $\{T_{i'}\}_{1\le i \le8}$.
After repeating the process $(k-1)$ times, we get the $3$-layer earth map tiling. 

\begin{figure}[htp]
	\centering
	\begin{tikzpicture}[>=latex,scale=0.6]
		
		\foreach \b in {0,1}
		{
			\begin{scope}[xshift=4*\b cm]	
				\draw (0,0)--(0,-2)--(1,-2)--(1,-4)--(2,-4)--(2,-6)
				(0,-2)--(-1,-2)--(-1,-4)--(-2,-4)--(-2,-6);
				
				\draw[line width=1.5]	
				(2,0)--(2,-4)
				(-2,0)--(-2,-4)
				(0,-2)--(0,-6);
				
				\draw[dashed]	
				(1,-2)--(2,-2)
				(0,-4)--(1,-4)
				(-1,-2)--(-2,-2)
				(0,-4)--(-1,-4);  
				
				\node at (1,-0.25) {\small $\bbb$};
				\node at (1,-1.75) {\small $\ccc$};
				\node at (0.25,-1.75) {\small $\aaa$};
				\node at (1.75,-1.75) {\small $\ddd$};
				\node at (0.75,-2.25) {\small $\aaa$};
				\node at (1.25,-2.3) {\small $\ccc$};
				\node at (0.25,-2.35) {\small $\bbb$};
				\node at (1.75,-2.3) {\small $\ddd$};
				
				\node at (0.75,-3.75) {\small $\ccc$};
				\node at (1.25,-3.75) {\small $\aaa$};
				\node at (0.25,-3.75) {\small $\ddd$};
				\node at (1.75,-3.7) {\small $\bbb$};
				
				\node at (0.25,-4.35) {\small $\ddd$};
				\node at (1.75,-4.35) {\small $\aaa$};
				\node at (1,-4.35) {\small $\ccc$};
				\node at (1,-5.8) {\small $\bbb$};
				
				\node at (-1,-0.25) {\small $\bbb$};
				\node at (-1,-1.75) {\small $\ccc$};
				\node at (-0.25,-1.75) {\small $\aaa$};
				\node at (-1.75,-1.75) {\small $\ddd$};
				\node at (-0.75,-2.25) {\small $\aaa$};
				\node at (-1.25,-2.3) {\small $\ccc$};
				\node at (-0.25,-2.35) {\small $\bbb$};
				\node at (-1.75,-2.3) {\small $\ddd$};
				
				\node at (-0.75,-3.75) {\small $\ccc$};
				\node at (-1.25,-3.75) {\small $\aaa$};
				\node at (-0.25,-3.75) {\small $\ddd$};
				\node at (-1.75,-3.7) {\small $\bbb$};
				
				\node at (-0.25,-4.35) {\small $\ddd$};
				\node at (-1.75,-4.35) {\small $\aaa$};
				\node at (-1,-4.35) {\small $\ccc$};
				\node at (-1,-5.8) {\small $\bbb$};

			\end{scope}
		}

		\node[draw,shape=circle, inner sep=0.5] at (-1,-1) {\small $1$};
		\node[draw,shape=circle, inner sep=0.5] at (1,-1) {\small $2$};
		\node[draw,shape=circle, inner sep=0.5] at (-0.5,-3) {\small $3$};
		\node[draw,shape=circle, inner sep=0.5] at (0.5,-3) {\small $4$};
		\node[draw,shape=circle, inner sep=0.5] at (-1.5,-3) {\small $5$};
		\node[draw,shape=circle, inner sep=0.5] at (1.5,-3) {\small $6$};
		\node[draw,shape=circle, inner sep=0.5] at (-1,-5) {\small $7$};
		\node[draw,shape=circle, inner sep=0.5] at (1,-5) {\small $8$};
		
		\node[draw,shape=circle, inner sep=0.5] at (-1+4,-1) {\small $1'$};
		\node[draw,shape=circle, inner sep=0.5] at (1+4,-1) {\small $2'$};
		\node[draw,shape=circle, inner sep=0.5] at (-0.5+4,-3) {\small $3'$};
		\node[draw,shape=circle, inner sep=0.5] at (0.5+4,-3) {\small $4'$};
		\node[draw,shape=circle, inner sep=0.5] at (-1.5+4,-3) {\small $5'$};
		\node[draw,shape=circle, inner sep=0.5] at (1.5+4,-3) {\small $6'$};
		\node[draw,shape=circle, inner sep=0.5] at (-1+4,-5) {\small $7'$};
		\node[draw,shape=circle, inner sep=0.5] at (1+4,-5) {\small $8'$};
		
		\fill (7,-3) circle (0.05);
		\fill (7.25,-3) circle (0.05);
		\fill (7.5,-3) circle (0.05);

	\end{tikzpicture}
	
	\caption{The $3$-layer earth map tiling $T(4k\, \aaa\ccc^2, 2\bbb^{2k},2k\,\ddd^4,2k\,\aaa^2\bbb^2)$.} \label{case1jie3}
\end{figure}
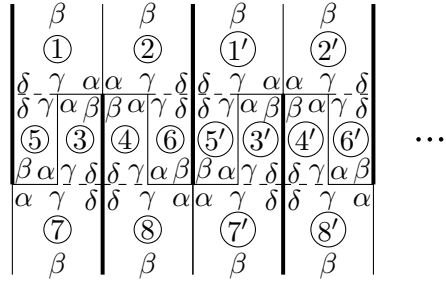

{\bf Remark}: The above discussion actually proved that $\thick\bbb_1\thin\bbb_2\thick$ determines $T_1,\dots,T_8$. This fact will be very useful to deduce other possible tilings.

\subsubsection*{Case $f=16s-4$}
The extra vertex is $\bbb^{2s-1}\ccc\ddd$ in Table \ref{tab-1.1}.
In the left of Fig. \ref{case1taotai1}, $\bbb^{2s-1}\ccc\ddd=\thick\bbb_1\thin\ccc_2\dash\ddd_3\thick\cdots$ determines $T_1,T_2,T_3$. Then $\aaa_1\aaa_2\cdots=\aaa_1\aaa_2\bbb_4\bbb$ determines $T_4$; $\aaa_4\bbb_2\cdots=\aaa_4\aaa\bbb_2\bbb_5$ determines $T_5$. We get $\ccc_3\ddd_2\ddd_5\cdots$, contradicting the AVC.

\begin{figure}[htp]
	\centering
	\begin{tikzpicture}[>=latex,scale=0.5] 
		\draw (0,0)--(-1,-2)--(-3,-2)
		(-1,-2)--(0,-4)--(1,-6)--(2,-4)
		(0,-4)--(-1,-6)
		(2,0)--(3,-2)--(1,-2);
		\draw[dashed](-3,-2)--(-2,0)
		(0,0)--(2,-4)
		(-2,-4)--(-1,-6) ;
		
		\draw[line width=1.5] (-2,0)--(2,0)
		(-2,-4)--(-1,-2)
		(0,-4)--(1,-2);
		

		\node at (-0.6,-0.5){\small $\bbb$};
		\node at (0,-0.8){\small $\ccc$};       
		\node at (0.6,-0.4){\small $\ddd$};
		\node at (-1.8,-0.4){\small $\ddd$};
		\node at (-2.4,-1.7){\small $\ccc$};
		\node at (-1.2,-1.7){\small $\aaa$};
		
		\node at (-0.6,-2){\small $\aaa$};
		\node at (0.55,-2){\small $\ddd$};
		\node at (0.1,-3.1){\small $\bbb$};
		\node at (-1.6,-2.4){\small $\bbb$};
		\node at (-1,-2.8){\small $\bbb$};
		
		\node at (-1.6,-4){\small $\ddd$};
		\node at (-1,-5.4){\small $\ccc$};
		\node at (1,-5.4){\small $\aaa$};
		\node at (-0.3,-4){\small $\aaa$};
		
		\node at (0,-4.7){\small $\aaa$};
		\node at (0.4,-4){\small $\bbb$};
		\node at (1.7,-4){\small $\ccc$};
		
		\node at (1.8,-0.45){\small $\bbb$};
		\node at (2.4,-1.7){\small $\aaa$};
		\node at (1.2,-1.7){\small $\ccc$};
		\node at (1,-2.7){\small $\ddd$};
		
		\node[draw,shape=circle, inner sep=0.5] at (-1.5,-1) {\small $1$};
		\node[draw,shape=circle, inner sep=0.5] at (0,-1.46) {\small $2$};
		\node[draw,shape=circle, inner sep=0.5] at (1.6,-1.1) {\small $3$};
		\node[draw,shape=circle, inner sep=0.5] at (-1,-4) {\small $4$};
		\node[draw,shape=circle, inner sep=0.5] at (1.1,-4) {\small $5$};
		
		
	\end{tikzpicture}\hspace{60pt}	
	\begin{tikzpicture}[>=latex,scale=0.5] 
		\draw (-5,-2)--(-3,0)--(3,0)--(5,-2)
		(-2,-6)--(0,-2)--(2,-6)
		(-2,-2)--(-1,-4)
		(2,-2)--(1,-4);
		\draw[dashed] (0,0)--(0,-2)
		(-5,-2)--(-2,-2)
		(5,-2)--(2,-2)
		(-3,-4)--(-2,-6)
		(3,-4)--(2,-6);
		
		\draw[line width=1.5] (-2,-2)--(0,0)--(2,-2)	
		(-2,-2)--(-3,-4)
		(2,-2)--(3,-4);
		
		\draw[dotted]
		(-1,-4)--(0,-6)
		(0,-6)--(1,-4);
		
		\node at (-1.1,-0.5){\small $\bbb$};
		\node at (1.2,-0.5){\small $\bbb$};
		\node at (-0.3,-0.8){\small $\ddd$};
		\node at (0.3,-0.8){\small $\ddd$};
		\node at (-2.8,-0.35){\small $\aaa$};
		\node at (2.8,-0.35){\small $\aaa$};
		\node at (-4.2,-1.7){\small $\ccc$};
		\node at (4.2,-1.7){\small $\ccc$};
		\node at (-2.2,-1.7){\small $\ddd$};
		\node at (2.2,-1.7){\small $\ddd$};
		\node at (-1.6,-2.1){\small $\bbb$};
		\node at (1.6,-2.3){\small $\bbb$};
		
		\node at (-0.35,-1.9){\small $\ccc$};
		\node at (0.35,-1.9){\small $\ccc$};
		\node at (0,-2.6){\small $\aaa$};
		
		\node at (-2,-2.8){\small $\bbb$};
		\node at (2,-2.9){\small $\bbb$};
		\node at (-2.6,-4){\small $\ddd$};
		\node at (2.60,-04){\small $\ddd$};
		\node at (-1.35,-4){\small $\aaa$};
		\node at (1.35,-4){\small $\aaa$};
		\node at (-2,-5.3){\small $\ccc$};
		\node at (2,-5.3){\small $\ccc$};
		
		\node at (-1,-3.3){\small $\aaa$};
		\node at (1,-3.3){\small $\aaa$};

		\node[draw,shape=circle, inner sep=0.5] at (-2.5,-1) {\small $1$};
		\node[draw,shape=circle, inner sep=0.5] at (-1,-2.2) {\small $2$};
		\node[draw,shape=circle, inner sep=0.5] at (1,-2) {\small $3$};
		\node[draw,shape=circle, inner sep=0.5] at (2.5,-1) {\small $4$};
		\node[draw,shape=circle, inner sep=0.5] at (-2,-4) {\small $5$};
		\node[draw,shape=circle, inner sep=0.5] at (0,-4) {\small $7$};
		\node[draw,shape=circle, inner sep=0.5] at (2,-4) {\small $6$};


	\end{tikzpicture}
	\caption{Vertex $\bbb^{2s-1}\ccc\ddd$ or $\bbb^{2s}\ddd^2$ appears.} \label{case1taotai1}
\end{figure}

\subsubsection*{Case $f=16s$, including $f=16$}
The extra vertices are $\aaa^4\,(f=16),\bbb^{2s}\ddd^2,\bbb^{4s}$ in Table \ref{tab-1.1}.
If $\aaa^4$ appears, then its AAD gives a vertex $\bbb\thin\ccc\cdots$ or $\ccc\thin\ccc\cdots$, contradicting the AVC in Table \ref{tab-1.1}. 

If $\bbb^{2s}\ddd^2$ appears, then $\bbb^{2s}\ddd^2=\thin\bbb_1\thick\ddd_2\dash\ddd_3\thick\bbb_4\thin\cdots$ in the right of Fig. \ref{case1taotai1}, which  determines $T_1,T_2,T_3,T_4$. Then $\dash\ddd_1\thick\bbb_2\thin\cdots=\dash\ddd_1\thick\bbb_2\thin\bbb_5\thick\cdots$ determines $T_5$;  $\dash\ddd_4\thick\bbb_3\thin\cdots=\dash\ddd_4\thick\bbb_3\thin\bbb_6\thick\cdots$ determines $T_6$. By AVC, $\ccc_2\ccc_3\cdots=\aaa_7\ccc_2\ccc_3$. By $\aaa_7$, either $\aaa_2\aaa_5\cdots$ or $\aaa_3\aaa_6\cdots$ is $\aaa^2\ccc\cdots$, contradicting the AVC. 

 If $\bbb^{4s}$ appears, by $\ddd\cdots=\ddd^4$ and the Claim after Table \ref{tab-1.1},  the tiling is the $3$-layer earth map tiling in Fig. \ref{case1jie3}.

\subsubsection*{Case $f=16s+8$, including $f=24$}
	
The extra vertices are $\aaa^3\,(f=24),\aaa\bbb^{2s+2},\bbb^{2s}\ccc^2,\bbb^{4s+2}$ in Table \ref{tab-1.1}. We divide our discussions into four subcases. 

\vspace{9pt}
\noindent{\textbf{Subcase. For $f=24$, $\aaa^3$ is a vertex.}}\label{subcase3a}

In the left of Fig. \ref{case1jie1}, no $\dash\ccc\thin\ccc\dash\cdots$ by AVC implies the unique AAD $\thin^{\ccc}\aaa_1^{\bbb}\thin^{\ccc}\aaa_2^{\bbb}\thin^{\ccc}\aaa_3^{\bbb}\thin$, which determines $T_1,T_2,T_3$. Then $\ddd_1\cdots=\ddd_2\cdots=\ddd_3\cdots=\ddd^4$ determines $T_4,T_5,\cdots,T_{12}$.
So $\aaa_4\aaa_{12}\cdots=\aaa^3$ or $\aaa^2\bbb^2$, shown respectively in Fig. \ref{case1jie1} and \ref{case1jie2}.

In the left of Fig. \ref{case1jie1}, $\aaa_4\aaa_{12}\cdots=\aaa_4\aaa_{12}\aaa_{15}=\thin^{\ccc}\aaa_4^{\bbb}\thin^{\ccc}\aaa_{12}^{\bbb}\thin^{\ccc}\aaa_{15}^{\bbb}\thin$ determines $T_{15}$.
Then $\ddd_{15}\cdots=\ddd^4$ determines $T_{16},T_{17},T_{18}$.
We have $\aaa_9\aaa_{10}\cdots=\aaa^3$ or $\aaa^2\bbb^2$. In the right of Fig. \ref{case1jie1}, $\aaa_9\aaa_{10}\cdots=\aaa_9\aaa_{10}\bbb_{13}\bbb$ determines $T_{13}$. We get $\aaa_{11}\aaa_{16}\ccc_{13}\cdots$, contradicting the AVC. Therefore, $\aaa_9\aaa_{10}\cdots=\aaa_9\aaa_{10}\aaa_{13}=\thin^{\bbb}\aaa_9^{\ccc}\thin^{\bbb}\aaa_{10}^{\ccc}\thin^{\bbb}\aaa_{13}^{\ccc}\thin$
determines $T_{13}$. Then $\ddd_{13}\cdots=\ddd^4$ determines $T_{14},T_{19},T_{24}$. Similarly, we have $\aaa_6\aaa_7\cdots=\aaa_6\aaa_7\aaa_{22}=\thin^{\bbb}\aaa_6^{\ccc}\thin^{\bbb}\aaa_{7}^{\ccc}\thin^{\bbb}\aaa_{22}^{\ccc}\thin$, which determines $T_{22}$. Then $\ddd_{22}\cdots=\ddd^4$ determines $T_{20}$, $T_{21}$, $T_{23}$.
This tiling is the quadrilateral subdivision of the octahedron. Each tile of this tiling is a $3444$-Tile, and it actually belongs to Section \ref{344d}.

\begin{figure}[htp]
	\centering
	\begin{tikzpicture}[>=latex]
		\begin{scope}[scale=2]
		\foreach \a in {0,1,2,3}
		\draw[rotate=90*\a]
		(0.4,0.4)--(0,0.55)--(0,0.9)
		(0,0.55)--(-0.4,0.4)
		(-1.3,1.3)--(0,1.3)--(1.3,1.3)
		(0,0.9)--(0,1.3);
		\foreach \a in {0,2}
		\draw[line width=1.5, rotate=90*\a]	
		(0,0)--(0.4,0.4)
		(-0.4,0.4)--(-0.9,0.9)
		(0,0.9)--(0.9,0.9)
		(0.9,0)--(0.9,0.9)
		(-0.9,0.9)--(-1.3,1.3)
		(1.3,1.3)--(1.6,1.6);
		\foreach \a in {0,2}
		\draw[dashed, rotate=90*\a]	
		(0,0)--(-0.4,0.4)
		(0.4,0.4)--(0.9,0.9)
		(-0.9,0)--(-0.9,0.9)--(0,0.9)
		(0.9,0.9)--(1.3,1.3)
		(-1.3,1.3)--(-1.6,1.6);
		
		\node[draw,shape=circle, inner sep=0.5] at (0,0.3) {\small $1$};
		\node[draw,shape=circle, inner sep=0.5] at (0.3,0.7) {\small $3$};
		\node[draw,shape=circle, inner sep=0.5] at (-0.3,0.7) {\small $2$};
		\node[draw,shape=circle, inner sep=0.5] at (0.3,0) {\small $4$};
		\node[draw,shape=circle, inner sep=0.5] at (0,-0.3) {\small $5$};
		\node[draw,shape=circle, inner sep=0.5] at (-0.3,0) {\small $6$};
		\node[draw,shape=circle, inner sep=0.5] at (-0.7,0.3) {\small $7$};
		\node[draw,shape=circle, inner sep=0.5] at (-1.1,0.5) {\small $8$};
		\node[draw,shape=circle, inner sep=0.5] at (-0.5,1.1) {\small $9$};
		\node[draw,shape=circle, inner sep=0.5] at (0.7,0.3) {\small $12$};
		\node[draw,shape=circle, inner sep=0.5] at (1.1,0.5) {\small $11$};
		\node[draw,shape=circle, inner sep=0.5] at (0.5,1.1) {\small $10$};
		
		\node[draw,shape=circle, inner sep=0.5] at (0,1.57) {\small $13$};
		\node[draw,shape=circle, inner sep=0.5] at (1.57,0) {\small $14$};
		\node[draw,shape=circle, inner sep=0.5] at (0.7,-0.3) {\small $15$};
		\node[draw,shape=circle, inner sep=0.5] at (1.1,-0.5) {\small $16$};
		\node[draw,shape=circle, inner sep=0.5] at (0.3,-0.7) {\small $17$};
		\node[draw,shape=circle, inner sep=0.5] at (0.5,-1.1) {\small $18$};
		\node[draw,shape=circle, inner sep=0.5] at (0,-1.57) {\small $19$};
		
		\node[draw,shape=circle, inner sep=0.5] at (-1.57,0) {\small $24$};
		\node[draw,shape=circle, inner sep=0.5] at (-0.7,-0.3) {\small $22$};
		\node[draw,shape=circle, inner sep=0.5] at (-1.1,-0.5) {\small $23$};
		\node[draw,shape=circle, inner sep=0.5] at (-0.3,-0.7) {\small $20$};
		\node[draw,shape=circle, inner sep=0.5] at (-0.5,-1.1) {\small $21$};
		
		\node at (0,1.37){\small $\aaa$};
		\node at (0,-1.37){\small $\aaa$};
		\node at (-1.37,0){\small $\aaa$};
		\node at (1.37,0){\small $\aaa$};

		\node at (0,0.13){\small $\ddd$};
		\node at (0,-0.13){\small $\ddd$};
		\node at (0.13,0){\small $\ddd$};
		\node at (-0.13,0){\small $\ddd$};
		\node at (0,0.46){\small $\aaa$};
		\node at (0,-0.46){\small $\aaa$};
		\node at (0.46,0){\small $\aaa$};
		\node at (-0.46,0){\small $\aaa$};
		\node at (0.63,0.1){\small $\aaa$};
		\node at (-0.63,0.1){\small $\aaa$};
		\node at (-0.63,-0.1){\small $\aaa$};
		\node at (0.63,-0.1){\small $\aaa$};
		\node at (0.1,0.6){\small $\aaa$};
		\node at (-0.1,0.6){\small $\aaa$};
		\node at (-0.1,-0.6){\small $\aaa$};
		\node at (0.1,-0.6){\small $\aaa$};
		\node at (0.37,0.25){\small $\bbb$};
		\node at (-0.355,-0.25){\small $\bbb$};
		\node at (0.247,0.347){\small $\bbb$};
		\node at (-0.238,-0.343){\small $\bbb$};
		\node at (0.37,-0.25){\small $\ccc$};
		\node at (0.247,-0.365){\small $\ccc$};
		\node at (-0.24,0.345){\small $\ccc$};
		\node at (-0.36,0.23){\small $\ccc$};
		\node at (-0.5,0.33){\small $\bbb$};
		\node at (-0.35,0.5){\small $\bbb$};
		\node at (-0.36,-0.5){\small $\ccc$};
		\node at (-0.5,-0.37){\small $\ccc$};
		\node at (0.1,0.78){\small $\bbb$};
		\node at (-0.1,0.78){\small $\ccc$};
		\node at (0.1,1){\small $\bbb$};
		\node at (-0.1,1){\small $\ccc$};
		\node at (0.1,1.22){\small $\aaa$};
		\node at (-0.1,1.22){\small $\aaa$};
		\node at (0.1,-1.22){\small $\aaa$};
		\node at (-0.1,-1.22){\small $\aaa$};
		\node at (0.1,-0.78){\small $\ccc$};
		\node at (-0.1,-0.79){\small $\bbb$};
		\node at (0.1,-1){\small $\ccc$};
		\node at (-0.1,-1){\small $\bbb$};
		\node at (0.8,0.1){\small $\bbb$};
		\node at (0.8,-0.1){\small $\ccc$};
		\node at (1,-0.1){\small $\ccc$};
		\node at (1,0.1){\small $\bbb$};
		\node at (1.2,0.1){\small $\aaa$};
		\node at (1.2,-0.1){\small $\aaa$};
		\node at (-1.2,0.1){\small $\aaa$};
		\node at (-1.2,-0.1){\small $\aaa$};
		\node at (-0.8,0.1){\small $\ccc$};
		\node at (-0.8,-0.1){\small $\bbb$};
		\node at (-1,-0.1){\small $\bbb$};
		\node at (-1,0.1){\small $\ccc$};
		\node at (-0.83,0.7){\small $\ddd$};
		\node at (-0.65,0.8){\small $\ddd$};
		\node at (-1,0.8){\small $\ddd$};
		\node at (-0.83,1){\small $\ddd$};
		\node at (0.83,0.7){\small $\ddd$};
		\node at (0.65,0.8){\small $\ddd$};
		\node at (1,0.8){\small $\ddd$};
		\node at (0.83,1){\small $\ddd$};
		\node at (-0.83,-0.7){\small $\ddd$};
		\node at (-0.65,-0.8){\small $\ddd$};
		\node at (-1,-0.8){\small $\ddd$};
		\node at (-0.83,-1){\small $\ddd$};
		\node at (0.83,-0.7){\small $\ddd$};
		\node at (0.65,-0.8){\small $\ddd$};
		\node at (1,-0.8){\small $\ddd$};
		\node at (0.83,-1){\small $\ddd$};
		\node at (1.1,1.22){\small $\ccc$};
		\node at (1.242,1.1){\small $\ccc$};
		\node at (-1.1,-1.22){\small $\ccc$};
		\node at (-1.242,-1.1){\small $\ccc$};
		\node at (0.5,-0.33){\small $\bbb$};
		\node at (0.35,-0.5){\small $\bbb$};
		\node at (0.36,0.5){\small $\ccc$};
		\node at (0.5,0.37){\small $\ccc$};
		\node at (-1,1.2){\small $\bbb$};
		\node at (-1.242,1.1){\small $\bbb$};
		\node at (1,-1.2){\small $\bbb$};
		\node at (1.242,-1.1){\small $\bbb$};
		\node at (0,1.78){\small $\ddd$};
		\node at (0,-1.78){\small $\ddd$};
		\node at (1.78,0){\small $\ddd$};
		\node at (-1.78,0){\small $\ddd$};
		\node at (-1.3,1.4){\small $\ccc$};
		\node at (-1.4,1.2){\small $\ccc$};
		\node at (1.3,-1.4){\small $\ccc$};
		\node at (1.4,-1.2){\small $\ccc$};
		\node at (1.3,1.4){\small $\bbb$};
		\node at (1.4,1.2){\small $\bbb$};
		\node at (-1.3,-1.4){\small $\bbb$};
		\node at (-1.4,-1.2){\small $\bbb$};

		\node at (0,-2.2){\small $T(8\aaa^3,6\ddd^4,12\bbb^2\ccc^2)$} ;
	\end{scope}

		\begin{scope}[xshift=6 cm,yshift=0cm,scale=0.8]	
			
			\draw (0,0)--(4,0)
			(1,0)--(1,-2)
			(3,0)--(3,-2);
			\draw[dashed] 
			(0,-2)--(1,-2)
			(2,0)--(2,-2)
			(3,-2)--(4,-2)
			(2,2)--(3,0);
			
			\draw[line width=1.5] (0,0)--(0,-2)
			(1,-2)--(3,-2)
			(1,0)--(2,2)
			(4,0)--(4,-2);

			\node[draw,shape=circle, inner sep=0.5] at (0.5,-1) {\small $9$};
			\node[draw,shape=circle, inner sep=0.5] at (1.5,-1) {\small $10$};
			\node[draw,shape=circle, inner sep=0.5] at (2.5,-1) {\small $11$};
			\node[draw,shape=circle, inner sep=0.5] at (3.5,-1) {\small $16$};
			\node[draw,shape=circle, inner sep=0.5] at (2,1) {\small $13$};
			
			\node at (0.75,-1.75){\small $\ccc$}; \node at (1.25,-1.75){\small $\bbb$};
			\node at (0.75+2,-1.75){\small $\bbb$}; \node at (1.25+2,-1.75){\small $\ccc$};
			\node at (0.25,-1.75){\small $\ddd$}; \node at (1.75,-1.75){\small $\ddd$};
			\node at (2.25,-1.75){\small $\ddd$}; \node at (3.75,-1.75){\small $\ddd$};

			\node at (0.25,-0.25){\small $\bbb$}; \node at (0.25+2,-0.25){\small $\ccc$};
			\node at (0.25+1.5,-0.25){\small $\ccc$}; 
			\node at (0.25+3.5,-0.25){\small $\bbb$}; \node at (0.75,-0.3){\small $\aaa$};
			\node at (1.25,-0.25){\small $\aaa$};  \node at (2.75,-0.25){\small $\aaa$};
			\node at (3.25,-0.25){\small $\aaa$}; 
			
			\node at (0.8,0.25){\small $\bbb$};  \node at (1.4,0.25){\small $\bbb$};  
			\node at (2,0.25){\small $\aaa$};  \node at (2.6,0.25){\small $\ccc$};  
			\node at (2,1.6){\small $\ddd$};  
			
			\node at (2,-2.8){\small $\aaa_9 \aaa_{10} \cdots = \aaa^2 \bbb^2$} ;
		\end{scope}

	\end{tikzpicture}
	\caption{Vertex $\aaa^3$ appears, $\aaa_4\aaa_{12}\cdots=\aaa^3$, and $\aaa_9\aaa_{10}\cdots=\aaa^3$ or  $\aaa^2 \bbb^2$.} \label{case1jie1}
\end{figure}
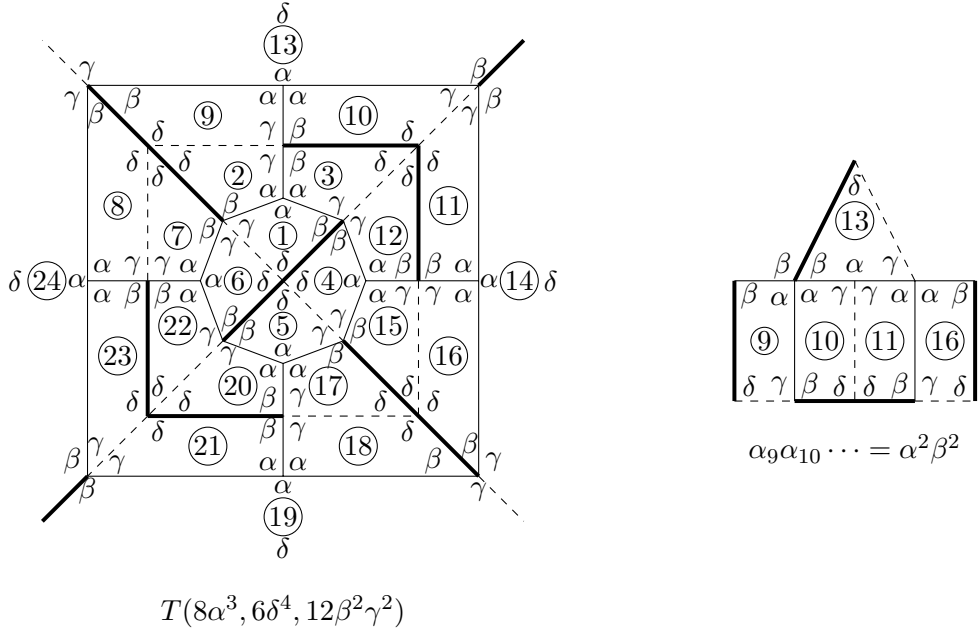

In the left of Fig. \ref{case1jie2},  $\aaa_4\aaa_{12}\cdots=\aaa_4\aaa_{12}\bbb_{13}\bbb_{14}$ determines $T_{13}$, $T_{14}$. Then $\ddd_{13}\ddd_{14}\cdots=\ddd^4$ determines $T_{15},T_{16}$. We have $\aaa_{15}\bbb_{5}\bbb_{6}\cdots=\aaa^2\bbb^2$ or $\aaa\bbb^4$. If $\aaa_{15}\bbb_{5}\bbb_{6}\cdots=\aaa\bbb^4$, then $\aaa_{6}\aaa_{7}\cdots=\aaa_{6}\aaa_{7}\aaa$. We get $\ccc_7\ccc_{8}\cdots=\ccc_7\ccc_{8}\ccc\cdots$, contradicting the AVC. Therefore,  $\aaa_{15}\bbb_{5}\bbb_{6}\cdots=\aaa^2\bbb^2$. By $\aaa_{18}$, we have $\aaa_{6}\aaa_{7}\cdots=\aaa_{6}\aaa_{7}\bbb_{17}\bbb_{18}$, which determines $T_{17},T_{18}$. Then $\ddd_{17}\ddd_{18}\cdots=\ddd^4$ determines $T_{19},T_{20}$. 
Similarly, we can determine $T_{21}$, $T_{22}$, $T_{23}$, $T_{24}$. This tiling turns out to be a flip modification of the quadrilateral subdivision of the octahedron in Fig. \ref{case1jie1}, as explained later using Fig. \ref{quadsubd} \& \ref{quadsubd2}. 
 

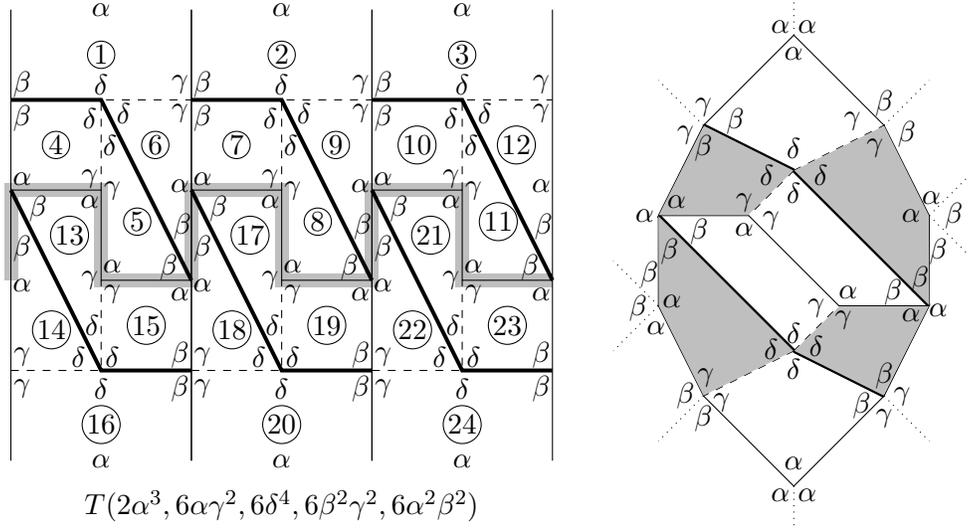
\begin{figure}[htp]
	\centering
	\begin{tikzpicture}[>=latex,scale=0.6]

		\fill[gray!50]
		(4,1.15) -- (4,0.85) -- (6,0.85)-- (6,1.15);
		\fill[gray!50]
		(6.15,1.15) -- (6.15,-1) -- (5.85,-1)-- (5.85,1.15);
		\fill[gray!50]
		(5.85,-0.85) -- (5.85,-1.15) -- (8,-1.15)-- (8,-0.85);
		\fill[gray!50]
		(8.15,1.15) -- (8.15,-1.15) -- (7.85,-1.15)-- (7.85,1.15);
		
		\fill[gray!50]
		(4+4,1.15) -- (4+4,0.85) -- (6+4,0.85)-- (6+4,1.15);
		\fill[gray!50]
		(6.15+4,1.15) -- (6.15+4,-1) -- (5.85+4,-1)-- (5.85+4,1.15);
		\fill[gray!50]
		(5.85+4,-0.85) -- (5.85+4,-1.15) -- (8+4,-1.15)-- (8+4,-0.85);
		\fill[gray!50]
		(8.15+4,1.15) -- (8.15+4,-1.15) -- (7.85+4,-1.15)-- (7.85+4,1.15);
		\fill[gray!50]
		(4+8,1.15) -- (4+8,0.85) -- (6+8,0.85)-- (6+8,1.15);
		\fill[gray!50]
		(6.15+8,1.15) -- (6.15+8,-1) -- (5.85+8,-1)-- (5.85+8,1.15);
		\fill[gray!50]
		(5.85+8,-0.85) -- (5.85+8,-1.15) -- (8+8,-1.15)-- (8+8,-0.85);
		
		\fill[gray!50]
		(3.85,1.15) -- (4.15,1.15) -- (4.15,-1)-- (3.85,-1);

		\foreach \b in {1,2,3}
		{
			\begin{scope}[xshift=4*\b cm]	
				\draw (0,-5)--(0,5)
				(4,-5)--(4,5)
				(0,1)--(2,1)--(2,-1)--(4,-1);
				
				\draw[line width=1.5]	
				(0,3)--(2,3)--(4,-1)
				(0,1)--(2,-3)--(4,-3);
				
				\draw[dashed]	
				(2,1)--(2,3)--(4,3)
				(0,-3)--(2,-3)--(2,-1);  
				
				\node at (2,5) {\small $\aaa$}; \node at (2,-5) {\small $\aaa$};
				\node at (0.25,3.35) {\small $\bbb$}; \node at (2,3.35) {\small $\ddd$};
				\node at (3.75,3.35) {\small $\ccc$}; \node at (3.75,2.7) {\small $\ccc$};
				\node at (0.25,2.6) {\small $\bbb$}; \node at (1.75,2.6) {\small $\ddd$};
				\node at (2.5,2.7) {\small $\ddd$}; \node at (2.2,2) {\small $\ddd$};
				\node at (0.25,1.25) {\small $\aaa$}; \node at (1.75,1.25) {\small $\ccc$};
				\node at (2.25,1) {\small $\ccc$}; \node at (1.75,0.75) {\small $\aaa$};
				\node at (0.6,0.6) {\small $\bbb$}; 
				\node at (1.75,-1) {\small $\ccc$}; \node at (2.25,-1.3) {\small $\ccc$};
				\node at (2.25,-0.7) {\small $\aaa$}; \node at (3.5,-0.7) {\small $\bbb$}; 
				\node at (3.8,0.2) {\small $\bbb$}; 
				\node at (3.75,-1.3) {\small $\aaa$}; \node at (2.25,-2.7) {\small $\ddd$};
				\node at (3.75,-2.65) {\small $\bbb$}; \node at (1.8,-2) {\small $\ddd$};
				\node at (1.5,-2.7) {\small $\ddd$}; \node at (2,-3.4) {\small $\ddd$};
				\node at (0.25,-2.7) {\small $\ccc$};
				\node at (0.25,-1.1) {\small $\aaa$};  \node at (0.25,-0.3) {\small $\bbb$};
				\node at (3.75,1.1) {\small $\aaa$};
				
				\node at (0.25,-3.4) {\small $\ccc$};  \node at (3.75,-3.4) {\small $\bbb$};
				
			\end{scope}
		}
		
		\node[draw,shape=circle, inner sep=0.5] at (6,4) {\small $1$};
		\node[draw,shape=circle, inner sep=0.5] at (10,4) {\small $2$};
		\node[draw,shape=circle, inner sep=0.5] at (14,4) {\small $3$};
		\node[draw,shape=circle, inner sep=0.5] at (5,2) {\small $4$};
		\node[draw,shape=circle, inner sep=0.5] at (6.8,0.3) {\small $5$};
		\node[draw,shape=circle, inner sep=0.5] at (7.2,2) {\small $6$};
		\node[draw,shape=circle, inner sep=0.5] at (9,2) {\small $7$};
		\node[draw,shape=circle, inner sep=0.5] at (10.8,0.3) {\small $8$};
		\node[draw,shape=circle, inner sep=0.5] at (11.2,2) {\small $9$};
		\node[draw,shape=circle, inner sep=0.5] at (13,2) {\small $10$};
		\node[draw,shape=circle, inner sep=0.5] at (14.8,0.3) {\small $11$};
		\node[draw,shape=circle, inner sep=0.5] at (15.2,2) {\small $12$};
		
		\node[draw,shape=circle, inner sep=0.5] at (5.3,0) {\small $13$};
		\node[draw,shape=circle, inner sep=0.5] at (4.9,-2.1) {\small $14$};
		\node[draw,shape=circle, inner sep=0.5] at (7,-2) {\small $15$};
		\node[draw,shape=circle, inner sep=0.5] at (6,-4.2) {\small $16$};
		
		\node[draw,shape=circle, inner sep=0.5] at (5.3+4,0) {\small $17$};
		\node[draw,shape=circle, inner sep=0.5] at (4.9+4,-2.1) {\small $18$};
		\node[draw,shape=circle, inner sep=0.5] at (7+4,-2) {\small $19$};
		\node[draw,shape=circle, inner sep=0.5] at (6+4,-4.2) {\small $20$};
		
		\node[draw,shape=circle, inner sep=0.5] at (5.3+8,0) {\small $21$};
		\node[draw,shape=circle, inner sep=0.5] at (4.9+8,-2.1) {\small $22$};
		\node[draw,shape=circle, inner sep=0.5] at (7+8,-2) {\small $23$};
		\node[draw,shape=circle, inner sep=0.5] at (6+8,-4.2) {\small $24$};

		\node at (10,-6){\small $T(2\aaa^3,6\aaa\ccc^2,6\ddd^4,6\bbb^2\ccc^2,6\aaa^2\bbb^2)$};
	\end{tikzpicture}\hspace{15pt}	
	\begin{tikzpicture}[>=latex,scale=0.6] 
		\draw (0,0)--(2,-2)
		(0,0)--(-2,-2)
		(2,-2)--(3,-4)
		(-2,-2)--(-3,-4)
		(3,-4)--(3,-6)
		(3,-6)--(2,-8)
		(2,-8)--(0,-10)
		(-3,-4)--(-3,-6)
		(-2,-8)--(-3,-6)
		(-2,-8)--(0,-10)
		(-3,-4)--(-1,-4)
		(-1,-4)--(1,-6)
		(1,-6)--(3,-6);
		\draw[dashed] (0,-3)--(2,-2)
		(0,-3)--(-1,-4)
		(0,-7)--(1,-6)
		(0,-7)--(-2,-8);
		
		\draw[line width=1.5] (-2,-2)--(0,-3)
		(0,-3)--(3,-6)
		(-3,-4)--(0,-7)
		(0,-7)--(2,-8);
		
		\draw[dotted] (2,-2)--(3,-1)
		(-2,-2)--(-3,-1)
		(3,-4)--(4,-3)
		(3,-4)--(4,-5)
		(-3,-6)--(-4,-5)
		(-3,-6)--(-4,-7)
		(-2,-8)--(-3,-9)
		(2,-8)--(3,-9)
		(0,0)--(0,1)
		(0,-10)--(0,-11);

		\fill[gray!50]
		(-2.99,-3.99) -- (-1.99,-2.01) -- (-0.01,-3)-- (-1.01,-3.99);
		
		\fill[gray!50]
		(2.99,-6.01) -- (1.99,-7.99) -- (0.01,-7)-- (1.01,-6.01);

		\fill[gray!50]
		(0.01,-3) -- (2,-2.01) -- (2.99,-4)-- (3,-5.99);
		
		\fill[gray!50]
		(-0.01,-7) -- (-2,-7.99) -- (-2.99,-6)-- (-3,-4.01);

		\node at (0,-0.4){\small $\aaa$};
		\node at (-1.3,-1.9){\small $\bbb$};
		\node at (1.3,-1.9){\small $\ccc$};
		\node at (-2,-2.5){\small $\bbb$};
		\node at (1.8,-2.5){\small $\ccc$};
		\node at (-2,-1.6){\small $\ccc$};
		\node at (2,-1.5){\small $\bbb$};
		\node at (-2.4,-2.2){\small $\ccc$};
		\node at (2.5,-2.2){\small $\bbb$};
		\node at (0,-2.5){\small $\ddd$};
		\node at (0,-3.5){\small $\ddd$};
		\node at (0.6,-3.1){\small $\ddd$};
		\node at (-0.6,-3.1){\small $\ddd$};
		\node at (-1.1,-3.7){\small $\ccc$};
		\node at (-0.5,-4.1){\small $\ccc$};
		\node at (-1.1,-4.3){\small $\aaa$};
		
		\node at (-2.6,-3.7){\small $\aaa$};
		\node at (-2.1,-4.35){\small $\bbb$};
		\node at (-2.75,-4.8){\small $\bbb$};
		\node at (-3.3,-4){\small $\aaa$};
		
		\node at (2.6,-4){\small $\aaa$};
		\node at (3.05,-3.6){\small $\aaa$};
		\node at (3.2,-4.8){\small $\bbb$};
		\node at (3.55,-4){\small $\bbb$};
		
		\node at (1.2,-5.7){\small $\aaa$};
		\node at (1.1,-6.3){\small $\ccc$};
		\node at (0.5,-6.05){\small $\ccc$};
		
		\node at (2.8,-5.4){\small $\bbb$};
		\node at (2.2,-5.7){\small $\bbb$};
		\node at (2.6,-6.2){\small $\aaa$};
		\node at (3.2,-6.1){\small $\aaa$};
		
		\node at (0,-6.6){\small $\ddd$};
		\node at (0,-7.4){\small $\ddd$};
		\node at (0.5,-6.9){\small $\ddd$};
		\node at (-0.5,-6.95){\small $\ddd$};
		
		\node at (-2.7,-5.9){\small $\aaa$};
		\node at (-3.2,-5.3){\small $\bbb$};
		\node at (-3.45,-6.05){\small $\bbb$};
		\node at (-3.05,-6.5){\small $\aaa$};
		
		\node at (-1.95,-7.6){\small $\ccc$};
		\node at (2,-7.55){\small $\bbb$};
		\node at (-1.6,-8.15){\small $\ccc$};
		\node at (1.5,-8.2){\small $\bbb$};
		\node at (-2.4,-8){\small $\bbb$};
		\node at (2.4,-8){\small $\ccc$};
		\node at (-2,-8.5){\small $\bbb$};
		\node at (2,-8.5){\small $\ccc$};
		\node at (0,-9.5){\small $\aaa$};
		
		\node at (0.3,0.2){\small $\aaa$};
		\node at (-0.3,0.2){\small $\aaa$};
		
		\node at (0.3,-10.2){\small $\aaa$};
		\node at (-0.3,-10.2){\small $\aaa$};

	\end{tikzpicture}
	\caption{Vertex $\aaa^3$ appears, and $\aaa_4\aaa_{12}\cdots=\aaa^2\bbb^2$.} \label{case1jie2}
\end{figure}

\vspace{9pt}
\noindent{\textbf{Subcase. $\beta^{4s+2}$ appears and $\aaa^3$ is not vertex}}

If $\bbb^{4s+2}$ appears, by $\ddd\cdots=\ddd^4$ and the Claim after Table \ref{tab-1.1}, the tiling is the $3$-layer earth map tiling in Fig. \ref{case1jie3}.

\vspace{9pt}
\noindent{\textbf{Subcase. $\aaa\bbb^{2s+2}$ appears and $\aaa^3,\beta^{4s+2}$ are not vertices}}

In the left of Fig. \ref{case1jie4}, the unique AAD of $\aaa\bbb^{2s+2}=\thick^{\ddd}\bbb_{2}^{\aaa}\thin^{\bbb}\aaa_{1}^{\ccc}\thin^{\aaa}\bbb_{3}^{\ddd}\thick\cdots$ determines $T_1,T_2,T_3$. Then $R(\thick\bbb_2\thin\aaa_1\thin\bbb_3\thick\cdots)=\bbb^{2s}$ and this $\bbb^{2s}$ determines $8s$ tiles by the Remark after the Claim. So we get
$\bbb_{11}\bbb_{12}\cdots=\bbb^{2s}\cdots=\aaa\bbb^{2s+2}$ or $\bbb^{2s}\ccc^2$, shown in  Fig. \ref{case1jie4} and \ref{case1jie5} respectively.

\vspace{9pt}
In Fig. \ref{case1jie4}, $\bbb_{11}\bbb_{12}\cdots=\aaa\bbb^{2s+2}$.
Then $\aaa_3\ccc_1\cdots=\aaa_3\ccc_1\ccc_4$ determines $T_{4}$; $\aaa_4\ccc_3\cdots=\aaa_4\ccc_3\ccc_5$ determines $T_{5}$.
By $\bbb_{11}\bbb_{12}\cdots=\aaa\bbb^{2s+2}$, we get $E_{67}=a$ or $b$. 
If $E_{67}=b$ in the right of Fig. \ref{case1jie4}, this determines $T_{6}$. So we have $\aaa_5\bbb_4\ccc_7\cdots$, contradicting the AVC. Therefore, $E_{67}=a$. 

By $\bbb_{11}\bbb_{12}\cdots=\aaa\bbb^{2s+2}$, we get $E_{78}=b$. This determines $T_7,T_{8}$. By $\aaa_6$, we get $\aaa_{11}\bbb\bbb_5\cdots=\aaa_{11}\bbb\bbb_5\bbb_6\cdots=\aaa_{11}\bbb^{2s+2}$, which determines $T_6$. Then $\aaa_8\aaa\bbb\cdots=\aaa_8\aaa\bbb\bbb_{9}$ determines $T_9$; $\aaa_9\ccc_8\cdots=\aaa_9\ccc_8\ccc_{10}$ determines $T_{10}$. 

By $\aaa\bbb\bbb_5\bbb_6\cdots=\aaa\bbb^{2s+2}$, this determines $8s$ tiles (including $T_1,T_3,T_4,T_5$) similarly and $\aaa_2\bbb_1\cdots=\bbb^{2s}\cdots=\aaa\bbb^{2s+2}$, which determines another $8s$ tiles (including $T_6,T_7$). So we obtain a new tiling $T((8s+4)\aaa\ccc^2,4\aaa\bbb^{2s+2},(4s+2)\ddd^4,4s\,\aaa^2\bbb^2)$.
It turns out to be the first flip modification of the $3$-layer earth map tiling in Fig. \ref{case1jie3}, as explained later using Fig. \ref{fliptiling1} \& \ref{flip1}. The flipped part is a hemisphere between two shaded vertical lines in Fig. \ref{case1jie4}.

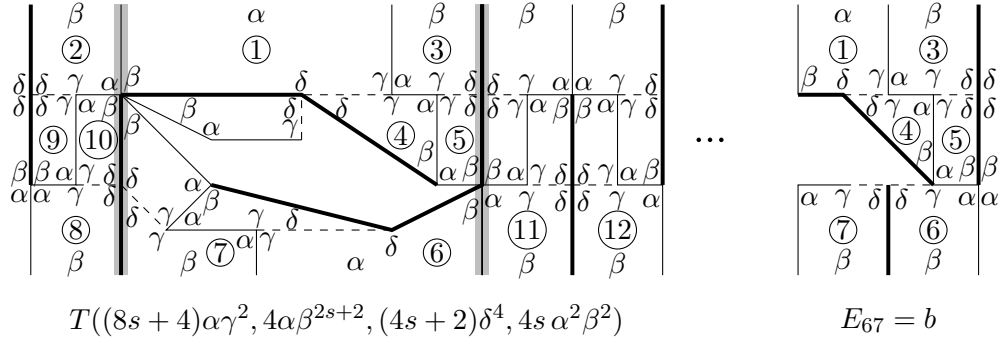
\begin{figure}[htp]
	\centering
	\begin{tikzpicture}[>=latex,scale=0.6]
		
	\begin{scope}[xshift=-1cm]
				
		\fill[gray!50]
		(3.85,0) -- (3.85,-6) -- (4.15,-6)-- (4.15,0);
		\fill[gray!50]
		(11.85,0) -- (12.15,0) -- (12.15,-6)-- (11.85,-6);

		\draw (0+14,0)--(0+14,-2)--(1+14,-2)--(1+14,-4)--(2+14,-4)--(2+14,-6)
		(0+14,-2)--(-1+14,-2)--(-1+14,-4)--(-2+14,-4)--(-2+14,-6)
		(2,-4)--(3,-4)--(3,-2)--(4,-2)--(4,0)
		(10,0)--(10,-2)--(11,-2)--(11,-4)--(12,-4)--(12,-6)
		(4,-2)--(6,-3)--(8,-3)
		(4,-2)--(6,-4)--(5,-5)--(7,-5)--(7,-6)
		(2,-4)--(2,-6);
		
		\draw[line width=1.5]	
		(2+14,0)--(2+14,-4)
		(-2+14,0)--(-2+14,-4)
		(0+14,-2)--(0+14,-6)
		(4,-2)--(4,-6)
		(4,-2)--(8,-2)--(11,-4)
		(12,0)--(12,-4)--(10,-5)--(6,-4)
		(2,0)--(2,-4);
		
		\draw[dashed]	
		(1+14,-2)--(2+14,-2)
		(0+14,-4)--(1+14,-4)
		(-1+14,-2)--(-2+14,-2)
		(0+14,-4)--(-1+14,-4)
		(2,-2)--(3,-2)
		(3,-4)--(4,-4)
		(8,-3)--(8,-2)--(10,-2)
		(11,-2)--(12,-2)
		(4,-4)--(5,-5)
		(7,-5)--(10,-5);  
		
		\node at (1+14,-0.25) {\small $\bbb$};
		\node at (1+14,-1.75) {\small $\ccc$};
		\node at (0.25+14,-1.75) {\small $\aaa$};
		\node at (1.75+14,-1.75) {\small $\ddd$};
		\node at (0.75+14,-2.25) {\small $\aaa$};
		\node at (1.25+14,-2.3) {\small $\ccc$};
		\node at (0.25+14,-2.35) {\small $\bbb$};
		\node at (1.75+14,-2.3) {\small $\ddd$};
		
		\node at (0.75+14,-3.75) {\small $\ccc$};
		\node at (1.25+14,-3.75) {\small $\aaa$};
		\node at (0.25+14,-3.75) {\small $\ddd$};
		\node at (1.75+14,-3.7) {\small $\bbb$};
		
		\node at (0.25+14,-4.35) {\small $\ddd$};
		\node at (1.75+14,-4.35) {\small $\aaa$};
		\node at (1+14,-4.35) {\small $\ccc$};
		\node at (1+14,-5.8) {\small $\bbb$};
		
		\node at (-1+14,-0.25) {\small $\bbb$};
		\node at (-1+14,-1.75) {\small $\ccc$};
		\node at (-0.25+14,-1.75) {\small $\aaa$};
		\node at (-1.75+14,-1.75) {\small $\ddd$};
		\node at (-0.75+14,-2.25) {\small $\aaa$};
		\node at (-1.25+14,-2.3) {\small $\ccc$};
		\node at (-0.25+14,-2.35) {\small $\bbb$};
		\node at (-1.75+14,-2.3) {\small $\ddd$};
		
		\node at (-0.75+14,-3.75) {\small $\ccc$};
		\node at (-1.25+14,-3.75) {\small $\aaa$};
		\node at (-0.25+14,-3.75) {\small $\ddd$};
		\node at (-1.75+14,-3.7) {\small $\bbb$};
		
		\node at (-0.25+14,-4.35) {\small $\ddd$};
		\node at (-1.75+14,-4.35) {\small $\aaa$};
		\node at (-1+14,-4.35) {\small $\ccc$};
		\node at (-1+14,-5.8) {\small $\bbb$};

		\node[draw,shape=circle, inner sep=0.5] at (7,-1) {\small $1$};
		\node[draw,shape=circle, inner sep=0.5] at (3,-1) {\small $2$};
		\node[draw,shape=circle, inner sep=0.5] at (11,-1) {\small $3$};
		\node[draw,shape=circle, inner sep=0.5] at (10.2,-2.9) {\small $4$};
		\node[draw,shape=circle, inner sep=0.5] at (11.5,-3) {\small $5$};
		\node[draw,shape=circle, inner sep=0.5] at (11,-5.5) {\small $6$};
		\node[draw,shape=circle, inner sep=0.5] at (6.2,-5.5) {\small $7$};			
		
		\node[draw,shape=circle, inner sep=0.5] at (2.5,-3) {\small $9$};
		\node[draw,shape=circle, inner sep=0.5] at (3.5,-3) {\small $10$};
		\node[draw,shape=circle, inner sep=0.5] at (3,-5) {\small $8$};
		
		\node[draw,shape=circle, inner sep=0.5] at (13,-5) {\small $11$};
		\node[draw,shape=circle, inner sep=0.5] at (15,-5) {\small $12$};

		\node at (1.75,-1.75) {\small $\ddd$};
		\node at (1.75,-2.25) {\small $\ddd$};
		\node at (1.75,-3.75) {\small $\bbb$};
		\node at (1.75,-4.25) {\small $\aaa$};

		\node at (2.25,-1.75) {\small $\ddd$};
		\node at (2.25,-2.25) {\small $\ddd$};
		\node at (2.25,-3.75) {\small $\bbb$};
		\node at (2.25,-4.25) {\small $\aaa$};
		\node at (3,-0.25) {\small $\bbb$};
		\node at (3,-1.75) {\small $\ccc$};
		\node at (3.75,-1.75) {\small $\aaa$};
		\node at (2.75,-2.25) {\small $\ccc$};
		\node at (3.25,-2.25) {\small $\aaa$};
		\node at (3.75,-2.3) {\small $\bbb$};		
		\node at (2.75,-3.75) {\small $\aaa$};
		\node at (3.25,-3.75) {\small $\ccc$};
		\node at (3.75,-3.75) {\small $\ddd$};
		\node at (3,-4.25) {\small $\ccc$};
		\node at (3.75,-4.25) {\small $\ddd$};
		\node at (3,-5.75) {\small $\bbb$};
		
		\node at (7,-0.25) {\small $\aaa$}; \node at (5.5,-5.75) {\small $\bbb$};
		\node at (11,-0.25) {\small $\bbb$}; \node at (9.2,-5.75) {\small $\aaa$};
		
		\node at (4.25,-1.65) {\small $\bbb$};\node at (5.5,-2.35) {\small $\bbb$};
		\node at (4.25,-2.7) {\small $\bbb$};
		\node at (4.25,-3.85) {\small $\ddd$};\node at (4.25,-4.7) {\small $\ddd$};
		\node at (5,-4.7) {\small $\ccc$};\node at (4.8,-5.25) {\small $\ccc$};
		\node at (5.6,-4.75) {\small $\aaa$}; \node at (6,-4.4) {\small $\bbb$};
		\node at (7,-4.75) {\small $\ccc$};\node at (7.8,-4.75) {\small $\ddd$};
		\node at (6.75,-5.25) {\small $\aaa$};\node at (7.25,-5.25) {\small $\ccc$};
		
		\node at (5.6,-4) {\small $\aaa$};
		
		\node at (6,-2.75) {\small $\aaa$};\node at (7.75,-2.75) {\small $\ccc$};
		\node at (8,-1.7) {\small $\ddd$};\node at (9.75,-1.7) {\small $\ccc$};
		\node at (10.25,-1.7) {\small $\aaa$}; \node at (10,-2.25) {\small $\ccc$};
		\node at (8.9,-2.3) {\small $\ddd$}; \node at (7.75,-2.3) {\small $\ddd$};
		\node at (11,-1.7) {\small $\ccc$};  \node at (11.75,-1.7) {\small $\ddd$}; 
		\node at (10.75,-2.25) {\small $\aaa$}; \node at (11.25,-2.25) {\small $\ccc$}; 
		\node at (11.75,-2.25) {\small $\ddd$}; \node at (11.25,-3.75) {\small $\aaa$}; 
		\node at (11.75,-3.65) {\small $\bbb$}; \node at (10.75,-3.3) {\small $\bbb$}; 
		\node at (10,-5.35) {\small $\ddd$}; \node at (11.75,-4.6) {\small $\bbb$};

		\fill (13+3.8,-3) circle (0.05);
		\fill (13.25+3.8,-3) circle (0.05);
		\fill (13.5+3.8,-3) circle (0.05);
		
		
		\node at (9,-7){\small $T((8s+4)\aaa\ccc^2,4\aaa\bbb^{2s+2},(4s+2)\ddd^4,4s\,\aaa^2\bbb^2)$};
	\end{scope}
		\begin{scope}[xshift=18cm]
			
			\draw (0,0)--(0,-2)
			(2,0)--(2,-2)--(3,-2)--(3,-4)--(4,-4)--(4,-6)
			(1,-4)--(0,-4)--(0,-6);
			
			\draw[line width=1.5]	
			(4,0)--(4,-4)
			(0,-2)--(1,-2)--(3,-4)
			(2,-4)--(2,-6);
			
			\draw[dashed]	
			(1,-2)--(2,-2)
			(3,-2)--(4,-2)
			(1,-4)--(3,-4);  
			
			\node[draw,shape=circle, inner sep=0.5] at (1,-1) {\small $1$};
			\node[draw,shape=circle, inner sep=0.5] at (3,-1) {\small $3$};
			\node[draw,shape=circle, inner sep=0.5] at (2.4,-2.8) {\small $4$};
			\node[draw,shape=circle, inner sep=0.5] at (3.5,-3) {\small $5$};
			\node[draw,shape=circle, inner sep=0.5] at (1,-5) {\small $7$};
			\node[draw,shape=circle, inner sep=0.5] at (3,-5) {\small $6$};
			
			\node at (1,-0.25) {\small $\aaa$};\node at (1,-1.65) {\small $\ddd$};
			\node at (3,-0.25) {\small $\bbb$};\node at (3,-1.65) {\small $\ccc$};
			\node at (0.25,-1.65) {\small $\bbb$};\node at (1.75,-1.65) {\small $\ccc$};
			\node at (2.25,-1.65) {\small $\aaa$};\node at (3.75,-1.65) {\small $\ddd$};
			\node at (2.1,-2.3) {\small $\ccc$}; \node at (1.65,-2.25) {\small $\ddd$}; \node at (2.75,-2.3) {\small $\aaa$};
			\node at (2.8,-3.3) {\small $\bbb$};
			\node at (3.25,-2.3) {\small $\ccc$}; \node at (3.75,-2.3) {\small $\ddd$};
			\node at (3.25,-3.7) {\small $\aaa$}; \node at (3.75,-3.65) {\small $\bbb$};
			
			\node at (3.75,-4.3) {\small $\aaa$}; \node at (3,-4.3) {\small $\ccc$};
			\node at (2.3,-4.3) {\small $\ddd$}; \node at (3,-5.7) {\small $\bbb$};
			
			\node at (0.25,-4.3) {\small $\aaa$}; \node at (1,-4.3) {\small $\ccc$};
			\node at (1.7,-4.3) {\small $\ddd$}; \node at (1,-5.7) {\small $\bbb$};
			
			\node at (4.25,-4.3) {\small $\aaa$}; \node at (4.25,-3.65) {\small $\bbb$};
			\node at (4.25,-2.3) {\small $\ddd$};\node at (4.25,-1.65) {\small $\ddd$};
			
			\node at (2,-7) {\small $E_{67}=b$};
		\end{scope}

	\end{tikzpicture}

	\caption{Vertex $\aaa\bbb^{2s+2}$ appears, $\bbb_{11}\bbb_{12}\cdots=\aaa\bbb^{2s+2}$, and   $E_{67} = a$ or $b$. } \label{case1jie4}
\end{figure}

\vspace{9pt}
In Fig. \ref{case1jie5}, $\bbb_{11}\bbb_{12}\cdots=\bbb^{2s}\ccc_6\ccc_7$ determines $T_6,T_7$. We get $T_4,T_5$ similarly. Then  $\aaa_7\aaa\bbb\cdots=\aaa^2\bbb^2$ determines $T_8$. So $\thick\bbb_4\thin\aaa_5\thin\bbb_6\thick\cdots=\aaa\bbb^{2s+2}$. This determines $8s$ tiles and $\ccc_2\ccc_8\cdots=\bbb^{2s}\cdots=\bbb^{2s}\ccc^2$. So we obtain a new tiling $T((8s+2)\aaa\ccc^2,2\aaa\bbb^{2s+2},2\bbb^{2s}\ccc^2,(4s+2)\ddd^4,(4s+2)\aaa^2\bbb^2)$.
It turns out to be the second flip modification of the $3$-layer earth map tiling in Fig. \ref{case1jie3}, as explained later using Fig. \ref{fliptiling2} \& \ref{flip2}. The flipped part lies between two shaded line segments in Fig. \ref{case1jie5}.

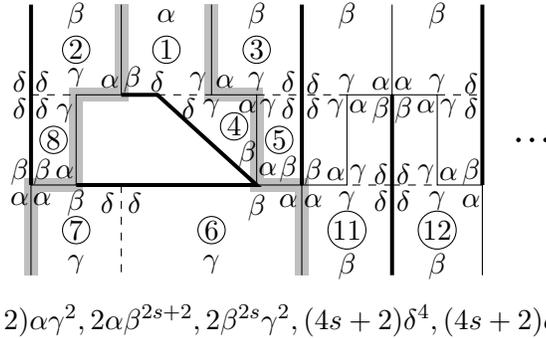
\begin{figure}[htp]
	\centering
	
	\begin{tikzpicture}[>=latex,scale=0.6]

		\fill[gray!50]
		(-0.15,-6) -- (0.15,-6) -- (0.15,-4)-- (-0.15,-4);
		
		\fill[gray!50]
		(-0.15,-3.85) -- (1,-3.85) -- (1,-4.15)-- (-0.15,-4.15);
		
		\fill[gray!50]
		(-0.15+1,-6+2-0.15) -- (0.15+1,-6+2-0.15) -- (0.15+1,-4+2)-- (-0.15+1,-4+2);
		
		\fill[gray!50]
		(-0.15+1,-3.85+2) -- (1+1,-3.85+2) -- (1+1,-4.15+2)-- (-0.15+1,-4.15+2);
		
		\fill[gray!50]
		(-0.15+2,-6+4-0.15) -- (0.15+2,-6+4-0.15) -- (0.15+2,-4+4)-- (-0.15+2,-4+4);
		
		\fill[gray!50]
		(-0.15+4,-6+4-0.15) -- (0.15+4,-6+4-0.15) -- (0.15+4,-4+4)-- (-0.15+4,-4+4);
		
		\fill[gray!50]
		(-0.15+4,-3.85+2) -- (1+4.15,-3.85+2) -- (1+4.15,-4.15+2)-- (-0.15+4,-4.15+2);
		
		\fill[gray!50]
		(-0.15+5,-6+2-0.15) -- (0.15+5,-6+2-0.15) -- (0.15+5,-4+2)-- (-0.15+5,-4+2);
		
		\fill[gray!50]
		(-0.15+5,-3.85) -- (1+5.15,-3.85) -- (1+5.15,-4.15)-- (-0.15+5,-4.15);
		
		\fill[gray!50]
		(-0.15+6,-6) -- (0.15+6,-6) -- (0.15+6,-4)-- (-0.15+6,-4);

		\draw (0+8,0)--(0+8,-2)--(1+8,-2)--(1+8,-4)--(2+8,-4)--(2+8,-6)
		(0+8,-2)--(-1+8,-2)--(-1+8,-4)--(-2+8,-4)--(-2+8,-6)
		(0,-6)--(0,-4)--(1,-4)--(1,-2)--(2,-2)--(2,0)
		(4,0)--(4,-2)--(5,-2)--(5,-4)--(6,-4);
		
		\draw[line width=1.5]	
		(2+8,0)--(2+8,-4)
		(-2+8,0)--(-2+8,-4)
		(0+8,-2)--(0+8,-6)
		(0,0)--(0,-4)
		(1,-4)--(5,-4)--(2.8,-2)--(2,-2);
		
		\draw[dashed]	
		(1+8,-2)--(2+8,-2)
		(0+8,-4)--(1+8,-4)
		(-1+8,-2)--(-2+8,-2)
		(0+8,-4)--(-1+8,-4)
		(0,-2)--(1,-2)
		(2,-4)--(2,-6)
		(2.8,-2)--(4,-2)
		(5,-2)--(6,-2);  
		
		\node at (1+8,-0.25) {\small $\bbb$};
		\node at (1+8,-1.75) {\small $\ccc$};
		\node at (0.25+8,-1.75) {\small $\aaa$};
		\node at (1.75+8,-1.75) {\small $\ddd$};
		\node at (0.75+8,-2.25) {\small $\aaa$};
		\node at (1.25+8,-2.3) {\small $\ccc$};
		\node at (0.25+8,-2.35) {\small $\bbb$};
		\node at (1.75+8,-2.3) {\small $\ddd$};
		
		\node at (0.75+8,-3.75) {\small $\ccc$};
		\node at (1.25+8,-3.75) {\small $\aaa$};
		\node at (0.25+8,-3.75) {\small $\ddd$};
		\node at (1.75+8,-3.7) {\small $\bbb$};
		
		\node at (0.25+8,-4.35) {\small $\ddd$};
		\node at (1.75+8,-4.35) {\small $\aaa$};
		\node at (1+8,-4.35) {\small $\ccc$};
		\node at (1+8,-5.8) {\small $\bbb$};
		
		\node at (-1+8,-0.25) {\small $\bbb$};
		\node at (-1+8,-1.75) {\small $\ccc$};
		\node at (-0.25+8,-1.75) {\small $\aaa$};
		\node at (-1.75+8,-1.75) {\small $\ddd$};
		\node at (-0.75+8,-2.25) {\small $\aaa$};
		\node at (-1.25+8,-2.3) {\small $\ccc$};
		\node at (-0.25+8,-2.35) {\small $\bbb$};
		\node at (-1.75+8,-2.3) {\small $\ddd$};
		
		\node at (-0.75+8,-3.75) {\small $\ccc$};
		\node at (-1.25+8,-3.75) {\small $\aaa$};
		\node at (-0.25+8,-3.75) {\small $\ddd$};
		\node at (-1.75+8,-3.7) {\small $\bbb$};
		
		\node at (-0.25+8,-4.35) {\small $\ddd$};
		\node at (-1.75+8,-4.35) {\small $\aaa$};
		\node at (-1+8,-4.35) {\small $\ccc$};
		\node at (-1+8,-5.8) {\small $\bbb$};

		\node[draw,shape=circle, inner sep=0.5] at (3,-1) {\small $1$};
		\node[draw,shape=circle, inner sep=0.5] at (1,-1) {\small $2$};
		\node[draw,shape=circle, inner sep=0.5] at (5,-1) {\small $3$};
		\node[draw,shape=circle, inner sep=0.5] at (5.5,-3) {\small $5$};
		\node[draw,shape=circle, inner sep=0.5] at (4.5,-2.7) {\small $4$};
		\node[draw,shape=circle, inner sep=0.5] at (4,-5) {\small $6$};
		\node[draw,shape=circle, inner sep=0.5] at (1,-5) {\small $7$};
		\node[draw,shape=circle, inner sep=0.5] at (0.5,-3) {\small $8$};
		
		\node[draw,shape=circle, inner sep=0.5] at (7,-5) {\small $11$};
		\node[draw,shape=circle, inner sep=0.5] at (9,-5) {\small $12$};
		
		\node at (3,-0.25) {\small $\aaa$};
		\node at (1,-0.25) {\small $\bbb$}; 	\node at (1,-5.75) {\small $\ccc$};
		\node at (5,-0.25) {\small $\bbb$};     \node at (4,-5.75) {\small $\ccc$};
		
		\node at (0.25,-1.7) {\small $\ddd$}; \node at (1,-1.6) {\small $\ccc$}; 
		\node at (1.75,-1.7) {\small $\aaa$}; \node at (0.25,-2.3) {\small $\ddd$};
		\node at (0.75,-2.4) {\small $\ccc$}; \node at (0.25,-3.7) {\small $\bbb$};
		\node at (0.75,-3.7) {\small $\aaa$};  \node at (0.25,-4.3) {\small $\aaa$};
		\node at (1,-4.4) {\small $\bbb$};  \node at (1.7,-4.4) {\small $\ddd$}; 
		\node at (2.3,-4.4) {\small $\ddd$};  \node at (5,-4.5) {\small $\bbb$}; 
		\node at (5.7,-4.35) {\small $\aaa$}; \node at (5.7,-3.65) {\small $\bbb$}; 
		\node at (5.25,-3.65) {\small $\aaa$};  \node at (4.8,-3.3) {\small $\bbb$};
		\node at (5.7,-2.3) {\small $\ddd$}; \node at (5.25,-2.3) {\small $\ccc$}; 
		\node at (5.7,-1.7) {\small $\ddd$}; \node at (5,-1.7) {\small $\ccc$}; 
		\node at (4.8,-2.25) {\small $\aaa$}; \node at (4,-2.3) {\small $\ccc$};
		\node at (3.5,-2.3) {\small $\ddd$}; \node at (3.7,-1.7) {\small $\ccc$};
		\node at (4.3,-1.7) {\small $\aaa$}; \node at (2.8,-1.7) {\small $\ddd$};
		\node at (2.25,-1.65) {\small $\bbb$};
		
		\node at (-0.25,-1.7) {\small $\ddd$}; \node at (-0.25,-2.3) {\small $\ddd$};
		\node at (-0.25,-3.7) {\small $\bbb$}; \node at (-0.25,-4.3) {\small $\aaa$};

		
		\fill (10.8,-3) circle (0.05);
		\fill (11.1,-3) circle (0.05);
		\fill (11.4,-3) circle (0.05);
		
		\node at (5,-7){\small $T((8s+2)\aaa\ccc^2,2\aaa\bbb^{2s+2},2\bbb^{2s}\ccc^2,(4s+2)\ddd^4,(4s+2)\aaa^2\bbb^2)$};
		
	\end{tikzpicture}

	\caption{Vertex $\aaa\bbb^{2s+2}$ appears, and  $\bbb_{11}\bbb_{12}\cdots=\bbb^{2s}\ccc^2$.} \label{case1jie5}
\end{figure}

\vspace{9pt}
\noindent{\textbf{Subcase. $\beta^{2s}\gamma^2$ appears and $\aaa^3,\beta^{4s+2},\aaa\bbb^{2s+2}$ are not vertices}}

In Fig. \ref{case1taotai2}, $\bbb^{2s}\ccc^2=\thick\bbb_1\thin\ccc_3\dash\ccc_4\thin\bbb_2\thick\cdots$ determines $T_1,T_2,T_3,T_4$. Then $\ddd_3\ddd_4\cdots=\ddd^4$ determines $T_5,T_6$; $\aaa_1\aaa_3\cdots=\aaa_1\aaa_3\bbb_7\bbb$ determines $T_7$, and $\aaa_7\bbb_3\bbb_5\cdots=\aaa_7\aaa_8\bbb_3\bbb_5$ (recall that $\aaa\bbb^{2s+2}$ is not a vertex). By $\aaa_8$, we get $\dash\ccc_7\thin\cdots=\dash\ccc_7\thin\bbb_8\thick\cdots$ which determines $T_8$. Then $\aaa_5\ccc_8\cdots=\aaa_5\ccc_8\ccc_9$ determines $T_9$; $\aaa_2\aaa_4\cdots=\aaa_2\aaa_4\bbb_{10}\bbb$ determines $T_{10}$, and  $\aaa_{10}\bbb_4\bbb_6\cdots=\aaa_{10}\aaa_{11}\bbb_4\bbb_6$. By $\aaa_{11}$, we have $\aaa_6\bbb_9\cdots \neq \aaa^2\bbb^2$, contradicting the AVC.

\begin{figure}[htp]
	\centering
	\begin{tikzpicture}[>=latex,scale=0.33] 
		\draw (0,0)--(3,-1)
		(0,0)--(-3,-1)
		(3,-1)--(6,0)
		(-3,-1)--(-6,0)
		(-3,-1)--(-3,-4)
		(3,-1)--(3,-4)
		(-3,-4)--(-7,-6)
		(-3,-4)--(-4,-7)
		(3,-4)--(7,-6)
		(3,-4)--(4,-7)
		(-4,-7)--(0,-6)
		(4,-7)--(0,-6) ;
		
		\draw[dashed] (0,0)--(0,-6)
		(-6,0)--(-3,3)
		(6,0)--(3,3)
		(-6,-3)--(-7,-6)
		(6,-3)--(7,-6)
		(-4,-7)--(-6,-9);
		
		\draw[line width=1.5] (0,0)--(3,3)
		(0,0)--(-3,3)
		(-3,-1)--(-6,-3)
		(3,-1)--(6,-3)
		(-3,-4)--(0,-3)
		(3,-4)--(0,-3)
		(-6,-9)--(4,-7)
		(-7,-6)--(-6,-9);
		
		\draw[dotted](4,-7)--(7,-8)
		(7,-8)--(7,-6);
		


		\node[draw,shape=circle, inner sep=0.5] at (-3,0.8) {\small $1$};
		\node[draw,shape=circle, inner sep=0.5] at (3,0.8) {\small $2$};
		\node[draw,shape=circle, inner sep=0.5] at (-1.5,-2) {\small $3$};
		\node[draw,shape=circle, inner sep=0.5] at (1.5,-2) {\small $4$};
		\node[draw,shape=circle, inner sep=0.5] at (-1.5,-5) {\small $5$};
		\node[draw,shape=circle, inner sep=0.5] at (1.5,-5) {\small $6$};
		\node[draw,shape=circle, inner sep=0.5] at (-4.5,-3.5) {\small $7$};
		\node[draw,shape=circle, inner sep=0.5] at (-5.3,-6.6) {\small $8$};
		\node[draw,shape=circle, inner sep=0.5] at (-2,-7.3) {\small $9$};
		\node[draw,shape=circle, inner sep=0.5] at (4.6,-3.5) {\small $10$};
		\node[draw,shape=circle, inner sep=0.5] at (5.5,-6.5) {\small $11$};

		\node at (-4.5,-1.3){\small $\bbb$};
		\node at (4.5,-1.2){\small $\bbb$};
		
		\node at (-1.3,0.1){\small $\bbb$};
		\node at (1.3,0.3){\small $\bbb$};
		
		\node at (-0.55,-1){\small $\ccc$}; \node at (0.55,-1){\small $\ccc$};
		\node at (-2.8,-0.5){\small $\aaa$}; \node at (2.9,-0.5){\small $\aaa$};
		\node at (-3,2.1){\small $\ddd$}; \node at (3,2.15){\small $\ddd$};
		\node at (-4.8,0.1){\small $\ccc$};  \node at (4.8,0.1){\small $\ccc$};
		\node at (-2.5,-1.4){\small $\aaa$}; \node at (2.5,-1.4){\small $\aaa$};
		\node at (-3.4,-2){\small $\bbb$}; \node at (3.45,-2.25){\small $\bbb$}; 
		\node at (-2.4,-3.2){\small $\bbb$}; \node at (2.5,-3.1){\small $\bbb$};
		\node at (-0.55,-2.6){\small $\ddd$}; \node at (0.55,-2.6){\small $\ddd$};
		\node at (-3.4,-3.7){\small $\aaa$};\node at (3.4,-3.7){\small $\aaa$};
		\node at (-5.6,-3.4){\small $\ddd$};\node at (5.7,-3.6){\small $\ddd$};
		\node at (-6.1,-5){\small $\ccc$};\node at (6.1,-5){\small $\ccc$};
		\node at (-3.7,-4.8){\small $\aaa$};\node at (3.7,-4.8){\small $\aaa$};
		\node at (-2.8,-4.7){\small $\bbb$};\node at (2.7,-4.7){\small $\bbb$};
		\node at (-0.55,-3.8){\small $\ddd$};\node at (0.55,-3.8){\small $\ddd$};
		\node at (-0.55,-5.6){\small $\ccc$}; \node at (0.55,-5.6){\small $\ccc$};
		\node at (-3.2,-6.3){\small $\aaa$}; \node at (3.2,-6.3){\small $\aaa$};
		\node at (0,-6.5){\small $\aaa$};\node at (1.5,-7){\small $\bbb$};
		\node at (-3.8,-7.5){\small $\ccc$};\node at (-4.3,-6.7){\small $\ccc$};
		\node at (-6.2,-6.3){\small $\bbb$};
		\node at (-5.8,-7.9){\small $\ddd$};
		\node at (-4.7,-8.2){\small $\ddd$};
		
	\end{tikzpicture}
	\caption{Vertex $\bbb^{2s}\ccc^2$ appears.} \label{case1taotai2}
\end{figure}

\subsection*{$\mathbf{H=\aaa\bbb^{d-1}\,(d=5,7,9,11)}$}

By $\aaa\ccc^2$, $\ddd^4$ and $\aaa\bbb^{d-1}$, we get
\begin{equation*}	
	\aaa=\pi-\tfrac{2\pi}{d-3}+\tfrac{8(d-1)\pi}{(d-3)f},\,\,
	\bbb=\tfrac{\pi}{d-3}-\tfrac{8\pi}{(d-3)f},\,\,
	\ccc=\tfrac{\pi}{2}+\tfrac{\pi}{d-3}-\tfrac{4(d-1)\pi}{(d-3)f},\,\,
	\ddd=\tfrac{\pi}{2}. 
\end{equation*}

By Lemma \ref{lemd4}, we have $\aaa^2\cdots$ and $\aaa<\pi$, which implies $f>4(d-1)$. Then $\aaa+2\bbb>\pi$, $\frac{\pi}{d-1}<\bbb<\frac{\pi}{d-3}$, $\frac{\pi}{2}<\ccc<\frac{\pi}{2}+\frac{\pi}{d-3}$.

By $\aaa\ccc^2$ and Lemma \ref{lemac2}, $\aaa^2\cdots=\aaa^k\bbb^{2t}$. By $\aaa+2\bbb>\pi$, we deduce that $t=0$ or $1$. Therefore, $\aaa^2\cdots=\aaa^k$ or $\aaa^k\bbb^2$.

Suppose $\ccc^2\cdots=\aaa^k\bbb^l\ccc^m\ddd^n$, by $R(\ccc^2\cdots)<\ccc+\ddd$ and Parity Lemma, we get $\ccc^2\cdots$ is even. By $\ccc>\frac\pi2$ and $2\ccc+2\ddd>2\pi$, we get $m=2$ and $n=0$. By $\aaa+2\ccc=2\pi$, we get $k=0$ or $1$. If $k=1$, then $\ccc^2\cdots=\aaa\ccc^2$. If $k=0$, by $\frac{\pi}{d-1}<\bbb<\frac{\pi}{d-3},\frac{\pi}{2}<\ccc<\frac{\pi}{2}+\frac{\pi}{d-3}$, we deduce that $d-5 < l < d-1$, which forces $l=d-3$. Therefore, $\ccc^2\cdots=\aaa\ccc^2$ or $\bbb^{d-3}\ccc^2$. If $\bbb^{d-3}\ccc^2$ is a vertex, we get $f=8(d-2)$, $\aaa=\frac{(d-3)\pi}{d-2}$, $\bbb=\frac{\pi}{d-2}$, $\ccc=\frac{(d-1)\pi}{2(d-2)}$, for $d=5,7,9,11$. Then $(\aaa,\bbb,\ccc)/\pi=(\frac23,\frac13,\frac23)$, $(\frac45,\frac15,\frac35)$,  $(\frac67,\frac17,\frac47)$ or $(\frac89,\frac19,\frac59)$.  
These belong exactly  to Case $f=16s+8$ in Table \ref{tab-1.1} for $s=1,2,3,4$, which has been classified there. They only admit the second flip modification of the $3$-layer earth map tilings, which have special $3344$-Tile and should not be considered to have special $334d$-Tile for $d=5,7,9,11$. 
So $\bbb^{d-3}\ccc^2$ is not a vertex and $\ccc^2\cdots=\aaa\ccc^2$. 
By Parity Lemma and $\aaa\ccc^2$, we have $\aaa\ccc\cdots=\aaa\ccc^2$. By Parity Lemma and $R(\ddd^3\cdots)=\ddd<\ccc$, we get $\ddd^3\cdots=\ddd^4$. Extend the 2nd picture of Fig. \ref{334d1} to Fig. \ref{334d2_taotai1} using $H=\aaa\bbb^{d-1}$ to show more neighborhood of the 1st special tile. 


\begin{figure}[htp]
	\centering
	\begin{tikzpicture}[>=latex,scale=0.55]
		\draw (0,0) -- (0,2) 
		(0,0) -- (2,0)
		(0,2)--(-2,4)
		(-2,4)--(2,4)
		(2,0)--(2,-2)
		(3,2)--(4,3)
		(2,4)--(4,3)
		(3,2)--(4,0)
		(2,0)--(4,0)
		(2,4)--(6,4)
		(4,0)--(6,-1)
		(6,-1)--(7,1);
		\draw[dashed]  (0,2)--(2,2)
		(0,0)--(-2,-2)
		(2,2)--(3,2)
		(4,3)--(6,2.5)
		(6,2.5)--(7,1)
		(4,-2)--(6,-1);
		\draw[line width=1.5] (2,0)--(2,2)
		(2,2)--(2,4)
		(-2,4)--(-2,-2)
		(-2,-2)--(2,-2)
		(6,4)--(6,2.5)
		(4,0)--(6,2.5)
		(2,0)--(4,-2);
		\draw[dotted] ;

		\fill (2,0) circle (0.2);

		\node at (0.35,0.35) {\small $\aaa$};
		\node at (0.2,-0.45) {\small $\ccc$};
		\node at (-0.35,0.25) {\small $\ccc$};
		
		\node at (1.65,0.45) {\small $\bbb$};
		\node at (1.65,-0.35) {\small $\aaa$};
		
		\node at (1.65,1.55) {\small $\ddd$};
		\node at (1.65,2.45) {\small $\ddd$};
		
		\node at (1.65,3.5) {\small $\bbb$};
		
		\node at (0,2.4) {\small $\ccc$};
		\node at (0.35,1.55) {\small $\ccc$};
		\node at (-0.4,1.75) {\small $\aaa$};
		
		\node at (-1.2,3.7) {\small $\aaa$};
		\node at (-1.65,3) {\small $\bbb$};
		\node at (-1.1,-1.6) {\small $\ddd$};
		\node at (-1.7,-1.2) {\small $\ddd$};
		\node at (1.65,-1.5) {\small $\bbb$};

		\node at (2.35,0.4) {\small $\bbb$};
		
		\node at (2.3,1.6) {\small $\ddd$};
		\node at (2.3,2.4) {\small $\ddd$};

		\node at (3.5,0.3) {\small $\aaa$};
		\node at (2.8,1.65) {\small $\ccc$};
		\node at (2.85,2.35) {\small $\ccc$};
		
		\node at (3.5,2.9) {\small $\aaa$};
		\node at (2.25,3.45) {\small $\bbb$};

		\node at (3.1,3.75) {\small $\aaa$};
		\node at (4,3.3) {\small $\ccc$};
		\node at (5.75,3.6) {\small $\bbb$};
		\node at (5.75,2.9) {\small $\ddd$};
		
		\node at (4,2.5) {\small $\ccc$};
		5			\node at (3.45,1.95) {\small $\aaa$};
		\node at (4.1,0.65) {\small $\bbb$};
		\node at (5.4,2.25) {\small $\ddd$};
		
		\node at (5.95,2) {\small $\ddd$};
		\node at (4.65,0.15) {\small $\bbb$};
		\node at (6.55,0.9) {\small $\ccc$};
		\node at (5.9,-0.55) {\small $\aaa$};
		
		\node at (3,-0.40) {\small $\bbb$};
		\node at (4,-0.3) {\small $\aaa$};
		\node at (4.1,-1.6) {\small $\ddd$};
		\node at (5.2,-1.0) {\small $\ccc$};

		\node at (6.3,2.7) {\small $\ddd$};

		\node[draw,shape=circle, inner sep=0.5] at (1,1) {\small $1$};
		\node[draw,shape=circle, inner sep=0.5] at (1,3) {\small $2$};
		\node[draw,shape=circle, inner sep=0.5] at (2.8,3) {\small $5$};
		\node[draw,shape=circle, inner sep=0.5] at (3,1) {\small $6$};
		\node[draw,shape=circle, inner sep=0.5] at (1,-1) {\small $4$};
		\node[draw,shape=circle, inner sep=0.5] at (-1,1) {\small $3$};
		\node[draw,shape=circle, inner sep=0.5] at (4,-0.9) {\small $7$};
		\node[draw,shape=circle, inner sep=0.5] at (4.4,1.8) {\small $8$};
		\node[draw,shape=circle, inner sep=0.5] at (5.7,0.8) {\small $9$};
		\node[draw,shape=circle, inner sep=0.5] at (4.9,3.5) {\small $10$};

		\begin{scope}[xshift=13 cm, yshift=0.5 cm]

		\draw (0,0) -- (0,2) 
		(0,0) -- (2,0)
		(0,2)--(-2,4)
		(-2,4)--(2,4)
		(2,0)--(2,-2)
		(3,2)--(4,3)
		(2,4)--(4,3)
		(3,2)--(4,0)
		(2,0)--(4,0)
		(2,0)--(3,-3)
		(3,-3)--(5,-3)
		(4,0)--(6,-1)
		(6,-1)--(8,-1)
		(8,-1)--(8,-3)
		(4,0)--(8,0.5);
		\draw[dashed]  (0,2)--(2,2)
		(0,0)--(-2,-2)
		(2,2)--(3,2)
		(4,3)--(6,2.5)
		(5,-2)--(5,-3)
		(5,-2)--(6,-1);
		\draw[line width=1.5] (2,0)--(2,2)
		(2,2)--(2,4)
		(-2,4)--(-2,-2)
		(-2,-2)--(2,-2)
		(4,0)--(6,2.5)	
		(2,0)--(5,-2)
		(5,-2)--(8,-3);
		\draw[dotted] 
		(8,0.5)--(8,-0.6)
		(8,-0.6)--(6,-1);

		\fill (2,0) circle (0.2);
		
		\node at (0.35,0.35) {\small $\aaa$};
		\node at (0.2,-0.45) {\small $\ccc$};
		\node at (-0.35,0.25) {\small $\ccc$};
		
		\node at (1.65,0.45) {\small $\bbb$};
		\node at (1.65,-0.35) {\small $\aaa$};
		
		\node at (1.65,1.55) {\small $\ddd$};
		\node at (1.65,2.45) {\small $\ddd$};
		
		\node at (1.65,3.5) {\small $\bbb$};
		
		\node at (0,2.4) {\small $\ccc$};
		\node at (0.35,1.55) {\small $\ccc$};
		\node at (-0.4,1.75) {\small $\aaa$};
		
		\node at (-1.2,3.7) {\small $\aaa$};
		\node at (-1.65,3) {\small $\bbb$};
		\node at (-1.1,-1.6) {\small $\ddd$};
		\node at (-1.7,-1.2) {\small $\ddd$};
		\node at (1.65,-1.5) {\small $\bbb$};
		\node at (2.35,0.4) {\small $\bbb$};

		\node at (2.3,1.6) {\small $\ddd$};
		\node at (2.3,2.4) {\small $\ddd$};		
		
		\node at (3.5,0.25) {\small $\aaa$};
		\node at (2.8,1.6) {\small $\ccc$};
		\node at (2.85,2.35) {\small $\ccc$};
		
		\node at (3.5,2.9) {\small $\aaa$};
		\node at (2.25,3.4) {\small $\bbb$};

		\node at (4,2.5) {\small $\ccc$};
		\node at (3.4,2) {\small $\aaa$};
		\node at (4.1,0.7) {\small $\bbb$};
		\node at (5.4,2.25) {\small $\ddd$};
		
		
		\node at (3.3,-0.38) {\small $\bbb$};
		\node at (4,-0.3) {\small $\aaa$};
		\node at (5.4,-1.1) {\small $\ccc$};
		\node at (5,-1.55) {\small $\ddd$};
		
		\node at (2.8,-1.1) {\small $\bbb$};
		\node at (3.25,-2.6) {\small $\aaa$};
		\node at (4.6,-2.2) {\small $\ddd$};
		\node at (4.6,-2.8) {\small $\ccc$};
		
		\node at (4.9,-0.15) {\small $\aaa$};
		\node at (6.1,-1.4) {\small $\ccc$};
		\node at (5.6,-1.8) {\small $\ddd$};
		\node at (7.6,-1.3) {\small $\aaa$};
		\node at (7.6,-2.4) {\small $\bbb$};
		
		\node at (5.3,-2.55) {\small $\ddd$};

		\node[draw,shape=circle, inner sep=0.5] at (1,1) {\small $1$};
		\node[draw,shape=circle, inner sep=0.5] at (1,3) {\small $2$};
		\node[draw,shape=circle, inner sep=0.5] at (2.8,3) {\small $5$};
		\node[draw,shape=circle, inner sep=0.5] at (3,1) {\small $6$};
		\node[draw,shape=circle, inner sep=0.5] at (1,-1) {\small $4$};
		\node[draw,shape=circle, inner sep=0.5] at (-1,1) {\small $3$};
		\node[draw,shape=circle, inner sep=0.5] at (4.4,-0.85) {\small $7$};
		\node[draw,shape=circle, inner sep=0.5] at (4.4,1.8) {\small $8$};
		\node[draw,shape=circle, inner sep=0.5] at (3.7,-1.9) {\small $9$};
		\node[draw,shape=circle, inner sep=0.5] at (6.75,-1.8) {\small $10$};
		\node[draw,shape=circle, inner sep=0.5] at (6.85,-0.3) {\small $11$};
		
	\end{scope}	
	\end{tikzpicture}	
	\caption{$H=\aaa\bbb^{d-1}$, $\aaa_6\aaa_7\cdots=\aaa^2\bbb^2$ or $\aaa^m\bbb^2(m\ge3)$.} \label{334d2_taotai1}
\end{figure}

In Fig. \ref{334d2_taotai1}, $H=\thin\aaa_4\thin\bbb_1\thick\bbb_6\thin\bbb_7\thick\cdots$  determines $T_7$. Then $\ccc_5\ccc_6\cdots=\aaa_8\ccc_5\ccc_6$. By $\aaa_8$, we deduce  $\aaa_6\aaa_7\cdots=\aaa^m\bbb^2(m\ge2)$. We discuss two cases $m=2$ and $m\ge3$ in Fig. \ref{334d2_taotai1}. 

In the left of Fig. \ref{334d2_taotai1}, $\aaa_6\aaa_7\cdots=\aaa^2\bbb^2$. This determines $T_8,T_9$. Then $\aaa_5\ccc_8\cdots=\aaa_5\ccc_8\ccc_{10}$ determines $T_{10}$ and $\ddd_8\ddd_9\ddd_{10}\cdots=\ddd^4$. Note that $T_8$ is a special $3344$-Tile which  has been  handled in Case $H=\aaa^{d-2}\bbb^2$ with $d=4$.

In the right of Fig. \ref{334d2_taotai1}, $\aaa_6\aaa_7\cdots=\aaa^m\bbb^2\,(m\ge3)$. We  determine $T_7,T_8,T_{9}$ similarly. So $\thick\bbb_8\thin\aaa_6\thin\aaa_7\thin\cdots=\thick\bbb_8\thin\aaa_6\thin\aaa_7\thin\aaa_{11}\thin\cdots$. Then by $\ccc>\ddd=\frac\pi2$, we get $\dash\ddd_{7}\thick\ddd_{9}\dash\cdots=\ddd^4$, which determines $T_{10}$. By $\aaa_{11}$, we have $\ccc_7\ccc_{10}\cdots\neq \aaa\ccc^2$, a contradiction. 
\end{proof}

\subsubsection*{Calculate the quadrilaterals in $3$-layer earth map tilings}

By $\aaa\ccc^2,\aaa^2\bbb^2,\bbb^{2n},\ddd^4$ in a $3$-layer earth map tiling with $f=8n$ tiles ($n\ge2$), we get
\[
\aaa=\tfrac{n-1}{n}\pi, \quad \bbb=\tfrac{1}{n}\pi,  \quad \ccc=\tfrac{n+1}{2n}\pi, \quad  \ddd=\tfrac{1}{2}\pi,
\]
as shown in Fig. \ref{jiejisuan}. By the sine and cosine law, we get
\begin{align*}
	\cos x&=\cos b \cos c =\cos^2 a -\sin^2 a  \cos \frac{\pi}{n},\,\,\,\,\,
	\cos a=\cot A \cot\frac{(n-1)\pi}{2n}\\
	\frac{\sin A}{\sin a}&=\frac{\sin (\pi-\frac{\pi}{n})}{\sin x},\quad \,\,\,\,
	\frac{\sin B}{\sin b}=\frac{\sin C}{\sin c}=\frac{1}{\sin x}, \quad\,\,\,\, a+2b=\pi.
\end{align*}
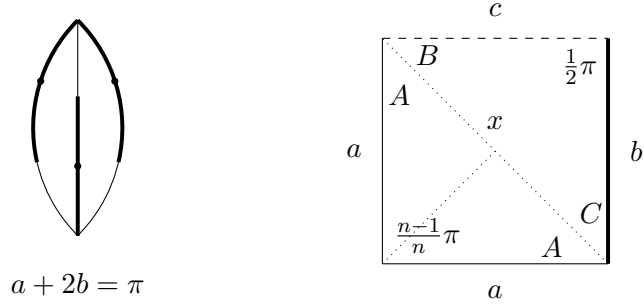
\begin{figure}[htp]
	\centering
	\begin{tikzpicture}[>=latex,scale=1.5] 
		\draw (0,0)--(2,0)
		(0,0)--(0,2);
		\draw[dashed]  (0,2)--(2,2);
		\draw[line width=1.5] (2,0)--(2,2);
		\draw[dotted] (2,0)--(0,2);
		\draw[dotted] (0,0)--(1,1);
		
		\node at (0.4,0.25){\small $\frac{n-1}{n}\pi$};
		
		\node at (1.5,0.15){\small $A$};
		\node at (1.85,0.45){\small $C$};
		
		\node at (1.75,1.75){\small $\frac{1}{2}\pi$};
		
		\node at (0.15,1.5){\small $A$};
		\node at (0.4,1.85){\small $B$};
		
		\node at (1,1.25){\small $x$};
		
		\node at (-0.25,1){\small $a$};
		\node at (1,-0.25){\small $a$};
		\node at (2.25,1){\small $b$};
		\node at (1,2.25){\small $c$};
		
		\begin{scope}[scale=0.45, xshift=-12 cm, yshift=1.8 cm] 
			\draw (6,-1)--(6,3);
			\draw[line width=1.5] (6,1.5)--(6,-1.25);
			\draw[line width=1.5] (6,3) arc (45:-13:3); 
			\draw[line width=1.5] (6,3) arc (-45+180:13+180:3); 
			
			\draw (6,3) arc (45:-45:3); \draw (6,3) arc (135:225:3);
			
			\fill (6.73,1.8) circle (0.07); \fill (5.27,1.8) circle (0.07);
			\fill (6,0.125) circle (0.07);
			\node at (6,-2.2){\small $a+2b=\pi$};
		\end{scope}		
	\end{tikzpicture}
	\caption{The quadrilateral in a $3$-layer earth map tiling with $f=8n$.} \label{jiejisuan}
\end{figure} 
Then
\begin{align*}
	\frac{\sin A}{\sin a}=&\frac{\sin(\pi-\frac{\pi}{n})}{\sin x}=\frac{\sin\frac{\pi}{n} \sin B}{\sin b}= \frac{\sin\frac{\pi}{n}\sin(\frac{n+1}{2n}\pi-A)}{\sin(\frac{\pi-a}{2})}=\frac{\sin\frac\pi n \cos(A-\frac{\pi}{2n})}{\cos\frac a2}\\
	&=\frac{\sin\frac\pi n (\cos A\cos\frac{\pi}{2n}+\sin A \sin \frac{\pi}{2n})}{\cos\frac a2}\\
	&=\frac{\sin\frac\pi n (\sin A\cos a \tan\frac{(n-1)\pi}{2n}\cos\frac{\pi}{2n}+\sin A \sin \frac{\pi}{2n})}{\cos\frac a2}
\end{align*}

After division by $\sin A$, we deduce a cubic equation of $\sin\frac{a}{2}$: 
\begin{equation*}
	8 (\cos \frac{\pi}{2n}\sin \frac{a}{2} )^3-4 (\cos \frac{\pi}{2n}\sin\frac{a}{2} )+ 1=0.
\end{equation*}
Thus $\cos \frac{\pi}{2n}\sin \frac{a}{2} = \frac{1}{2}$, $\frac{\pm\sqrt{5}-1}{4}$. 
Note that $0<a<1$ implies $\sin\frac{a}{2}>0$. If $\sin\frac{a}{2}=\frac{1}{2} \sec\frac{\pi}{2n}$, then $A=\frac\pi n$ and $C=0$, a contradiction. So we get a unique solution $\sin\frac{a}{2}=\frac{(\sqrt{5}-1)}{4} \sec\frac{\pi}{2n}$, and we conclude that 
\begin{equation*}
	a=2\arcsin \frac{\sqrt{5}-1}{4\cos \frac{\pi}{2n}},  \quad
	b=\frac{\pi-a}{2},  \quad
	c=\arccos \frac{(3-\sqrt{5})\cos^2\frac{\pi}{2n}+\sqrt{5}-2}{\cos \frac{\pi}{2n}}.
\end{equation*}

For $f=16$, we get $\aaa=\frac{\pi}{2}$, $\bbb=\frac{\pi}{2}$, $\ccc=\frac{3\pi}{4}$, $a\approx0.2879\pi$, $b\approx0.3560\pi$, $c\approx0.1615\pi$. 

For $f=24$, we get $\aaa=\frac{2\pi}{3}$, $\bbb=\frac{\pi}{3}$, $\ccc=\frac{2\pi}{3}$, $a\approx0.2323\pi$, $b\approx0.3839\pi$, $c\approx0.1161\pi$. This quadrilateral also gives the first tiling in Proposition \ref{334dCase1}.

As $f=8n \to \infty$, we get $\aaa \nearrow \pi$, $\bbb \searrow 0$, $\ccc \searrow \frac\pi2$, $a \searrow \frac{\pi}{5}$, $b \nearrow \frac{2\pi}{5}$, $c \searrow 0$. In summary $a,b,c$ are distinct for all $n\ge 2$ and the quadrilateral is indeed of Type $a^2bc$.

\subsubsection*{$3$-layer earth map tilings and their flip modifications}

The previous computation shows that there is a unique quadrilateral admitting the $3$-layer earth map tiling $T(4n\, \aaa\ccc^2, 2\bbb^{2n},2n\,\ddd^4,2n\,\aaa^2\bbb^2)$ with $n\ge2$ time zones ($8n$ tiles). Furthermore, there are two more tilings in Fig. \ref{case1jie4} \& \ref{case1jie5} when $n$ is odd, i.e. $n=2m+1$, which implies $\ccc=(m+1)\bbb$, $\aaa=2m\bbb$. These two tilings can be explained as two flip modifications of the $3$-layer earth map tiling as follows. 

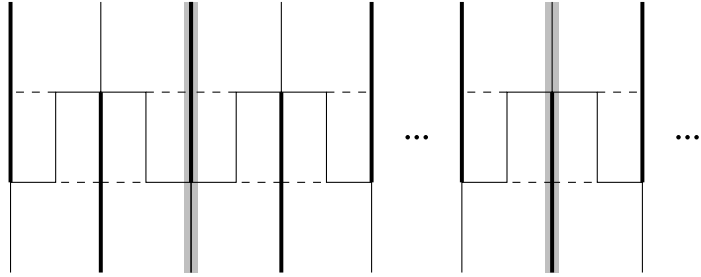
\begin{figure}[htp]
	\centering
	\begin{tikzpicture}[>=latex,scale=0.6]
		
		\begin{scope}[xshift=2 cm]			
			\fill[gray!50]
			(3.85,0) -- (3.85,-6) -- (4.15,-6)-- (4.15,0);
			\fill[gray!50]
			(11.85,0) -- (12.15,0) -- (12.15,-6)-- (11.85,-6);
		\end{scope}
		\foreach \b in {2,4,7}
		{
			\begin{scope}[xshift=2*\b cm]	
				\draw (0,0)--(0,-2)--(1,-2)--(1,-4)--(2,-4)--(2,-6)
				(0,-2)--(-1,-2)--(-1,-4)--(-2,-4)--(-2,-6);
				
				\draw[line width=1.5]	
				(2,0)--(2,-4)
				(-2,0)--(-2,-4)
				(0,-2)--(0,-6);
				
				\draw[dashed]	
				(1,-2)--(2,-2)
				(0,-4)--(1,-4)
				(-1,-2)--(-2,-2)
				(0,-4)--(-1,-4);  
			\end{scope}
		}
		\fill (10.8,-3) circle (0.05); \fill (11,-3) circle (0.05); \fill (11.2,-3) circle (0.05); 
		\fill (16.8,-3) circle (0.05); \fill (17,-3) circle (0.05); \fill (17.2,-3) circle (0.05);

	\end{tikzpicture}

	\caption{A hemisphere between two thick grey longitudes.} \label{fliptiling1}
\end{figure}

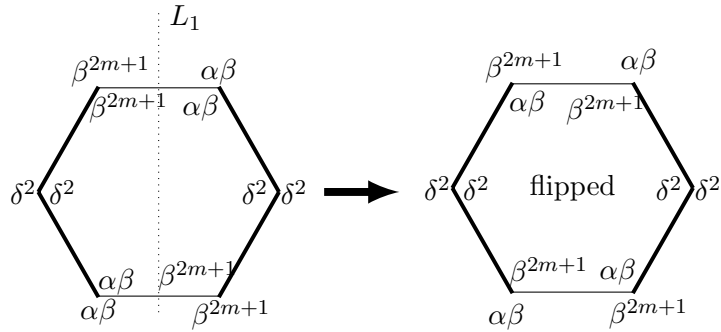
\begin{figure}[htp]
	\centering

	\begin{tikzpicture}[>=latex,scale=0.4]
		\coordinate (O) at (0,0);
		\def\m{6}
		\pgfmathsetmacro\i{\m-1}
		\foreach \x in {0,...,\i}
		{
			\def\pointname{\x}
			\coordinate (\pointname) at ($(0,0) +(\x*360/\m:4cm)$)  ;
			
		}
		
		\draw (0)
		\foreach \x in {0,...,\i}
		{ -- (\x) } -- cycle;

		\draw[line width=1.5] (2,3.5)--(4,0)
		(2,-3.5)--(4,0)
		(-2,3.5)--(-4,0)
		(-2,-3.5)--(-4,0);
		
		\draw[dotted] (0,-4)--(0,6);
		
		\node at (1.4,2.8){\small $\aaa\bbb$};
		\node at (2,4){\small $\aaa\bbb$};
		
		\node at (-1.4,-3){\small $\aaa\bbb$};
		\node at (-2,-4){\small $\aaa\bbb$};
		
		\node at (3.2,0){\small $\ddd^2$};
		\node at (4.5,0){\small $\ddd^2$};
		\node at (-3.2,0){\small $\ddd^2$};
		\node at (-4.5,0){\small $\ddd^2$};
		
		\node at (-1,2.8){\small $\bbb^{2m+1}$};
		\node at (-1.6,4){\small $\bbb^{2m+1}$};
		
		\node at (1.3,-2.8){\small $\bbb^{2m+1}$};
		\node at (2.4,-4.1){\small $\bbb^{2m+1}$};	
		
		\node at (0.9,5.8){\small $L_1$};
		
		\draw[line width=3pt, ->]
		(5.5,0) -- (8,0);
		
		
	\end{tikzpicture}
	\begin{tikzpicture}[>=latex,scale=0.4]
		\coordinate (O) at (0,0);
		\def\m{6}
		\pgfmathsetmacro\i{\m-1}
		\foreach \x in {0,...,\i}
		{
			\def\pointname{\x}
			\coordinate (\pointname) at ($(0,0) +(\x*360/\m:4cm)$)  ;
			
		}
		
		\draw (0)
		\foreach \x in {0,...,\i}
		{ -- (\x) } -- cycle;
		
		
		\draw[line width=1.5] (2,3.5)--(4,0)
		(2,-3.5)--(4,0)
		(-2,3.5)--(-4,0)
		(-2,-3.5)--(-4,0);

		\node at (-1.5,2.8){\small $\aaa\bbb$};
		\node at (2.2,4.2){\small $\aaa\bbb$};
		
		\node at (1.5,-2.8){\small $\aaa\bbb$};
		\node at (-2.2,-4.3){\small $\aaa\bbb$};
		
		\node at (3.2,0){\small $\ddd^2$};
		\node at (4.5,0){\small $\ddd^2$};
		\node at (-3.2,0){\small $\ddd^2$};
		\node at (-4.5,0){\small $\ddd^2$};
		
		\node at (1.1,2.6){\small $\bbb^{2m+1}$};
		\node at (-1.6,4){\small $\bbb^{2m+1}$};
		
		\node at (-0.8,-2.8){\small $\bbb^{2m+1}$};
		\node at (2.520,-4.2){\small $\bbb^{2m+1}$};	
		
		\node at (0,0){\small $\text{flipped}$};
		
	\end{tikzpicture}
	\caption{Flip a hemisphere in special $3$-layer earth map tilings.} \label{flip1}
\end{figure}

The angles and edges along the thick grey boundary in Fig. \ref{fliptiling1} are indicated in Fig. \ref{flip1}. Note that this boundary is a full great circle. Then the flip of the enclosed hemisphere with respect to the line $L_1$ keeps the angle sums of all vertices and transforms the $3$-layer earth map tiling to a new tiling $T((8m+4)\aaa\ccc^2,4\aaa\bbb^{2m+2},(4m+2)\ddd^4,4m\,\aaa^2\bbb^2)$ in Fig. \ref{case1jie4}. 

The angles and edges along the thick grey boundary in Fig. \ref{fliptiling2} are  indicated in Fig. \ref{flip2}. Note that there are $2m+1$ time zones in total. The $10$ $a$-edges in thick grey boundary cross over $m+1$ time zones with $4$ tiles left out, so they divide the tiling into two identical halves (but they are not hemispheres). The flip with respect to the line $L_2$ keeps the angle sums of all vertices and transforms the $3$-layer earth map tiling to a new tiling $T((8m+2)\aaa\ccc^2,2\aaa\bbb^{2m+2},2\bbb^{2m}\ccc^2,(4m+2)\ddd^4,(4m+2)\aaa^2\bbb^2)$  in Fig. \ref{case1jie5}. Equivalently, this flip modification can also be viewed as the  clockwise rotation of the inner half by $\frac{2\pi}{5}$.

\begin{figure}[htp]
	\centering	
	\begin{tikzpicture}[>=latex,scale=0.6]
		
		\begin{scope}[xshift=2 cm]			
			\draw [line width=5pt, gray!50]
			(6,0)--(6,-2) -- (5,-2)--(5,-4) -- (4,-4)-- (4,-6);
			\draw [line width=5pt, gray!50]
			(12,0)--(12,-2) -- (13,-2)--(13,-4) -- (14,-4)-- (14,-6);
		\end{scope}
		\foreach \b in {2,4,7}
		{
			\begin{scope}[xshift=2*\b cm]	
				\draw (0,0)--(0,-2)--(1,-2)--(1,-4)--(2,-4)--(2,-6)
				(0,-2)--(-1,-2)--(-1,-4)--(-2,-4)--(-2,-6);
				
				\draw[line width=1.5]	
				(2,0)--(2,-4)
				(-2,0)--(-2,-4)
				(0,-2)--(0,-6);
				
				\draw[dashed]	
				(1,-2)--(2,-2)
				(0,-4)--(1,-4)
				(-1,-2)--(-2,-2)
				(0,-4)--(-1,-4);  
			\end{scope}
		}
		\fill (10.8,-3) circle (0.05); \fill (11,-3) circle (0.05); \fill (11.2,-3) circle (0.05); 
		\fill (16.8,-3) circle (0.05); \fill (17,-3) circle (0.05); \fill (17.2,-3) circle (0.05);

	\end{tikzpicture}	
	
	\caption{Two identical halves with thick grey boundary.} \label{fliptiling2}
\end{figure}	

\begin{figure}[htp]
	\centering
	\begin{tikzpicture}[>=latex,scale=0.45]
		\coordinate (O) at (0,0);
		\def\m{10}
		\pgfmathsetmacro\i{\m-1}
		\foreach \x in {0,...,\i}
		{
			\def\pointname{\x}
			\coordinate (\pointname) at ($(0,0) +(\x*360/\m:4cm)$)  ;
			
		}
		
		\draw (0)
		\foreach \x in {0,...,\i}
		{ -- (\x) } -- cycle;

		\node at (1,3){\small $\aaa\bbb^2$};
		\node at (2.8,2.2){\small $\aaa$};
		\node at (3.3,0){\small $\ccc^2$};
		\node at (2.8,-2.2){\small $\aaa$};
		\node at (1.1,-3.1){\small $\bbb^{2m+2}$};

		\node at (-1,3.2){\small $\bbb^{2m}$};
		\node at (-2.5,2){\small $\aaa\bbb^2$};
		\node at (-3.4,0){\small $\aaa$};
		\node at (-2.8,-2){\small $\ccc^2$};
		\node at (-1.1,-3.1){\small $\aaa$};

		\node at (1.2,4.2){\small $\aaa$};
		\node at (3.5,2.8){\small $\ccc^2$};
		\node at (4.4,0){\small $\aaa$};
		\node at (4,-2.5){\small $\aaa\bbb^2$};
		\node at (1.6,-4.3){\small $\bbb^{2m}$};

		\node at (-1.2,4.2){\small $\bbb^{2m+2}$};
		\node at (-3.5,2.8){\small $\aaa$};
		\node at (-4.4,0){\small $\ccc^2$};
		\node at (-3.6,-2.5){\small $\aaa$};
		\node at (-1.6,-4.4){\small $\aaa\bbb^2$};

		\node at (2.2,5){\small $L_2$};

		\draw[dotted] (1.95,6)--(-2,-6);
		
		\draw[line width=3pt, ->]
		(5.5,0) -- (8,0);
		
		
	\end{tikzpicture}
	\begin{tikzpicture}[>=latex,scale=0.45]
		\coordinate (O) at (0,0);
		\def\m{10}
		\pgfmathsetmacro\i{\m-1}
		\foreach \x in {0,...,\i}
		{
			\def\pointname{\x}
			\coordinate (\pointname) at ($(0,0) +(\x*360/\m:4cm)$)  ;
			
		}
		
		\draw (0)
		\foreach \x in {0,...,\i}
		{ -- (\x) } -- cycle;

		\node at (1,3){\small $\aaa\bbb^2$};
		\node at (2.4,1.9){\small $\bbb^{2m}$};
		\node at (3.1,0){\small $\aaa\bbb^2$};
		\node at (2.8,-2.2){\small $\aaa$};
		\node at (1.1,-3.1){\small $\ccc^2$};

		\node at (-1,3.2){\small $\aaa$};
		\node at (-2.5,2){\small $\ccc^2$};
		\node at (-3.4,0){\small $\aaa$};
		\node at (-2,-1.7){\small $\bbb^{2m+2}$};
		\node at (-1.1,-3.1){\small $\aaa$};

		\node at (1.2,4.2){\small $\aaa$};
		\node at (3.7,2.9){\small $\ccc^2$};
		\node at (4.4,0){\small $\aaa$};
		\node at (4,-2.5){\small $\aaa\bbb^2$};
		\node at (1.6,-4.3){\small $\bbb^{2m}$};

		\node at (-1.2,4.2){\small $\bbb^{2m+2}$};
		\node at (-3.5,2.8){\small $\aaa$};
		\node at (-4.4,0){\small $\ccc^2$};
		\node at (-3.6,-2.6){\small $\aaa$};
		\node at (-1.6,-4.4){\small $\aaa\bbb^2$};



		\node at (0,0){\small $\text{flipped}$};
		\node at (0,-6){\small };	
	\end{tikzpicture}	
	\caption{Flip half of the sphere in special $3$-layer earth map tilings.} \label{flip2}
\end{figure}  

\subsubsection*{Quadrilateral subdivision tilings and a special flip modification}

The quadrilateral subdivision was introduced in \cite[Section 3.2]{wy1}, obtained by combining any tiling of a closed surface with its dual. Such subdivisions of Platonic solids produce some tilings of the sphere by congruent quadrilaterals of Type $a^4$ with $12$ tiles (tetrahedron), Type $a^2b^2$ with $24$ tiles (cube or octahedron) and $60$ tiles (dodecahedron or icosahedron)  respectively. Furthermore the octahedron's subdivision admits $1$-parameter deformations of Type $a^2bc$, which will reduce to Type $a^3b$ for a particular parameter. The left of Fig. \ref{quadsubd} is such a quadrilateral  subdivision of one triangular face of the regular octahedron. Replacing all triangular faces by this subdivision, we get a tiling $T(8\aaa^3,6\ddd^4,12\bbb^2\ccc^2)$ of the sphere as shown in the middle picture of Fig. \ref{quadsubd}.

\begin{figure}[htp]
	\centering
	
	\begin{tikzpicture}[>=latex,scale=0.6]
		
		\begin{scope}[xshift=-5 cm,yshift=-1.2 cm, scale=0.8]	
			\draw[dashed]
			(-2,0)--(1,0)
			(-1.5,0.87)--(0,3.46)
			(0.5,2.59)--(2,0);
			\draw[line width=1.5]
			(1,0)--(2,0)
			(-2,0)--(-1.5,0.87)
			(0,3.46)--(0.5,2.59);
			\draw
			(0,1.15)--(0.5,2.59)
			(0,1.15)--(1,0)
			(0,1.15)--(-1.5,0.87);
			
			\node at (0,-0.8){\small };
		\end{scope}
		\begin{scope}[xshift=2 cm,yshift=0.4 cm,scale=0.4]
			
			\foreach \a in {0,1,2}
			{
				\begin{scope}[rotate=120*\a ]
					\draw [line width=3pt, gray!50] (-3.8,0)--({-2*sqrt(3)},-2)--(-2,-3.2) --(0,-4)--(2,-3.2)--({2*sqrt(3)},-2);
				\end{scope}				
			}

			\foreach \i in {0,1,2}
			{
				\begin{scope}[rotate=120*\i ]
					
					\draw (0,0)--({-2*sqrt(3)},-2)
					(0,0)--({2*sqrt(3)},-2)
					(0,-6)--(0,-8)
					(0,-4) --(0,-6)
					({-2*sqrt(3)},-2)--(-2,-3.2) --(0,-4)--(2,-3.2) --({2*sqrt(3)},-2);
					
					\draw[line width=1.5] 
					({-3*sqrt(3)},-3)--(-2,-3.2)
					({-sqrt(3)},-1)--(0,-2)--(2,-3.2)
					({-3*sqrt(3)},-3)--({-3*sqrt(3)},3);

					\draw[dashed] 
					({-3*sqrt(3)},-3)--(0,-6)
					({sqrt(3)},-1)--(0,-2)--(-2,-3.2)
					({-3*sqrt(3)},-3)--(-3.8,0);
				\end{scope}
			}

		\end{scope}
		\begin{scope}[xshift=10 cm,yshift=0.4 cm,scale=0.4]

			\foreach \i in {0,1,2}
			{
				\begin{scope}[rotate=120*\i ]
					\draw [line width=3pt, gray!50]({-2*sqrt(3)},-2)--(-2,-3.2) --(0,-4)--(2,-3.2) --({2*sqrt(3)},-2);
					
					\draw (0,0)--({-2*sqrt(3)},-2)
					(0,0)--({2*sqrt(3)},-2)
					(0,-6)--(0,-8)
					(2,-3.2) --(0,-6)
					({-2*sqrt(3)},-2)--(-2,-3.2) --(0,-4)--(2,-3.2) --({2*sqrt(3)},-2);
					
					\draw[line width=1.5] ({-2*sqrt(3)},-2)--({-3*sqrt(3)},-3)--(0,-6)
					({2*sqrt(3)},-2)--({3*sqrt(3)},-3)
					({-sqrt(3)},-1)--(0,-2)--(2,-3.2);

					\draw[dashed] ({3*sqrt(3)},-3)--(0,-6)
					({sqrt(3)},-1)--(0,-2)--(-2,-3.2) 
					({-3*sqrt(3)},-3)--(0,-4);
				\end{scope}
			}
		\end{scope}
	\end{tikzpicture}  	
	
	\caption{Quadrilateral subdivisions of the octahedron and a flip if $\bbb=\frac{\pi}{3}$. }
	\label{quadsubd}
\end{figure}
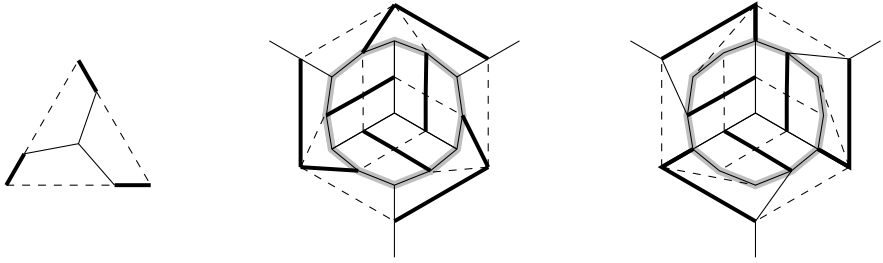

\begin{figure}[htp]
	\centering	
	
	\begin{tikzpicture}[>=latex,scale=0.45]
		\coordinate (O) at (0,0);
		\def\n{12}
		\pgfmathsetmacro\i{\n-1}
		\foreach \x in {0,...,\i}
		{
			\def\pointname{\x}
			\coordinate (\pointname) at ($(0,0) +(\x*360/\n:4cm)$)  ;
			
		}
		
		\draw (0)
		\foreach \x in {0,...,\i}
		{ -- (\x) } -- cycle;

		\node at (0,3.4){\small $\aaa^2$};
		\node at (1.6,3){\small $\ccc^2$};
		\node at (3,1.8){\small $\aaa$};
		\node at (3.4,0){\small $\bbb^2$};
		\node at (3,-1.8){\small $\aaa^2$};
		\node at (1.6,-3){\small $\ccc^2$};
		\node at (0,-3.4){\small $\aaa$};
		
		\node at (-1.5,-2.9){\small $\bbb^2$};
		\node at (-2.8,-1.8){\small $\aaa^2$};
		\node at (-3.4,0){\small $\ccc^2$};
		\node at (-3,1.8){\small $\aaa$};
		\node at (-1.6,3){\small $\bbb^2$};
		
		\node at (0,4.5){\small $\aaa$};
		\node at (2.5,3.9){\small $\bbb^2$};
		\node at (4,2.5){\small $\aaa^2$};
		\node at (4.6,0){\small $\ccc^2$};
		\node at (3.9,-2.4){\small $\aaa$};
		\node at (2.4,-3.9){\small $\bbb^2$};
		\node at (0,-4.5){\small $\aaa^2$};
		
		\node at (-2.5,3.9){\small $\ccc^2$};
		\node at (-4,2.5){\small $\aaa^2$};
		\node at (-4.6,0){\small $\bbb^2$};
		\node at (-3.9,-2.4){\small $\aaa$};
		\node at (-2.4,-3.9){\small $\ccc^2$};

		\node at (4,5.2){\small $L_{\aaa=\ccc=2\bbb}$};

		\draw[dotted] (1.8,6)--(-1.8,-6);
		
		\draw[line width=3pt, ->]
		(6.5,0) -- (11,0);
		
		\node at (8.5,1){\small $\aaa=\ccc=2\bbb$};
		
	\end{tikzpicture}\hspace{5pt}
	\begin{tikzpicture}[>=latex,scale=0.45]
		\coordinate (O) at (0,0);
		\def\n{12}
		\pgfmathsetmacro\i{\n-1}
		\foreach \x in {0,...,\i}
		{
			\def\pointname{\x}
			\coordinate (\pointname) at ($(0,0) +(\x*360/\n:4cm)$)  ;
			
		}
		
		\draw (0)
		\foreach \x in {0,...,\i}
		{ -- (\x) } -- cycle;

		\node at (0,3.4){\small $\ccc^2$};
		\node at (1.6,3){\small $\aaa^2$};
		\node at (3,1.6){\small $\bbb^2$};
		\node at (3.4,0){\small $\aaa$};
		\node at (3,-1.8){\small $\ccc^2$};
		\node at (1.6,-3){\small $\aaa^2$};
		\node at (0,-3.4){\small $\bbb^2$};
		
		\node at (-1.5,-2.9){\small $\aaa$};
		\node at (-2.8,-1.8){\small $\ccc^2$};
		\node at (-3.4,0){\small $\aaa^2$};
		\node at (-3,1.8){\small $\bbb^2$};
		\node at (-1.6,3){\small $\aaa$};
		
		\node at (0,4.5){\small $\aaa$};
		\node at (2.5,3.9){\small $\bbb^2$};
		\node at (4,2.5){\small $\aaa^2$};
		\node at (4.6,0){\small $\ccc^2$};
		\node at (3.9,-2.4){\small $\aaa$};
		\node at (2.4,-3.9){\small $\bbb^2$};
		\node at (0,-4.5){\small $\aaa^2$};
		
		\node at (-2.5,3.9){\small $\ccc^2$};
		\node at (-4,2.5){\small $\aaa^2$};
		\node at (-4.6,0){\small $\bbb^2$};
		\node at (-3.9,-2.4){\small $\aaa$};
		\node at (-2.4,-3.9){\small $\ccc^2$};

		\node at (0,0){\small $\text{flipped}$};
		\node at (0,-5.5){\small };
	\end{tikzpicture}
	\caption{Flip modification of a special $\bbb=\frac{\pi}{3}$ quadrilateral subdivision. }
	\label{quadsubd2} 
\end{figure}
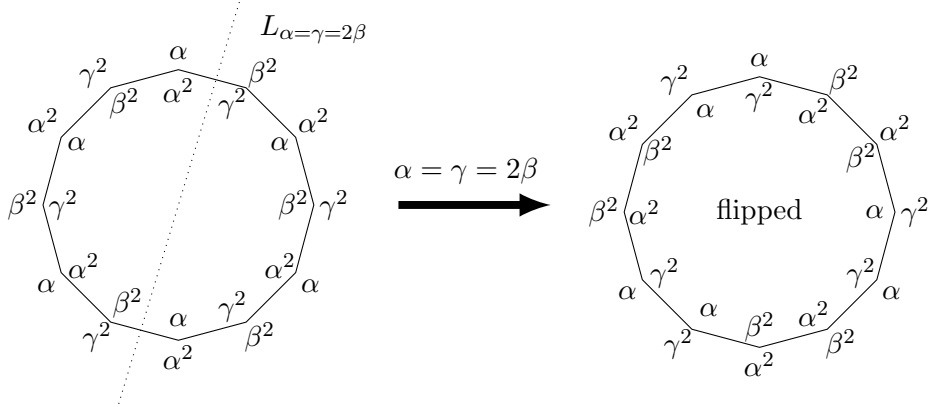	

Then we explain the flip modification for the special $\bbb=\frac{\pi}{3}$ case. In the middle of Fig. \ref{quadsubd}, we may use the thick grey lines to divide the  tiling  into two identical halves (but they are not hemispheres). The boundary between these two halves is illustrated by the left of Fig. \ref{quadsubd2}. Since $\bbb=\frac{\pi}{3}$ implies $\aaa=\ccc=2\bbb$, we may flip the inner half with respect to the line $L$ in Fig. \ref{quadsubd2}, and still keep the angle sums of all the vertices along the boundary to be $2\pi$. So the flip modification gives a new tiling  $T(2\aaa^3,6\aaa\ccc^2,6\ddd^4,6\bbb^2\ccc^2,6\aaa^2\bbb^2)$, as shown in the right of Fig. \ref{quadsubd}.  Its authentic $3$D picture is given in Fig. \ref{real figure}.

\vspace{9pt}
\begin{proposition} \label{334dCase2}
 There is no $a^2bc$-tiling with the 2nd special tile in Fig. \ref{33--}.
\end{proposition}

\begin{proof}
Let the second of Fig. \ref{33--} be the center tile $T_1$ in the partial neighborhood in Fig.  \ref{334dF3}. 
By Lemma \ref{lem4}, we get the degree $4$ vertex $\aaa_4\bbb_1\cdots=\aaa^2\bbb^2=\aaa_4\aaa_5\bbb_1\bbb_6$. This determines $T_6$. 

\begin{figure}[htp]
	\centering
	\begin{tikzpicture}[>=latex,scale=0.45]       
		\draw (0,0) -- (0,2) 
		(0,0)--(2,0)
		(0,2)--(-2,4)
		(-2,4)--(2,4)
		(2,0)--(2,-2)
		(4,0)--(5,2)
		(2,0)--(4,0)
		(-2,4)--(-3,6)
		(2,4)--(5,6)
		(5,6)--(4,7)
		(5,6)--(4,4)
		(4,4)--(7,6)
		(7,6)--(6,4)
		(6,4)--(7,3)
		(7,3)--(5,2)
		(4,7)--(6,9)
		(10,6)--(6,4);
		\draw[dashed]  (0,2)--(2,2)
		(0,0)--(-2,-2)
		(2,2)--(5,2)
		(-3,6)--(2,6)
		(2,6)--(4,7)
		(2,2)--(4,4)
		(6,9)--(8,8)
		(8,8)--(7,6)
		(4,0)--(4,-2);
		\draw[line width=1.5] (2,0)--(2,2)
		(2,2)--(2,4)
		(-2,-2)--(-2,4)
		(-2,-2)--(2,-2)
		(2,6)--(2,4)
		(2,2)--(6,4)
		(4,-2)--(7,3)
		(5,6)--(8,8)
		(8,8)--(10,6)
		(2,-2)--(4,-2);
		\draw[dotted] 
		(7,3)--(9,4)
		(9,4)--(10,6);
		
		\fill (2,2) circle (0.2);
		
		\node[draw,shape=circle, inner sep=0.5] at (1,1) {\small $1$};
		\node[draw,shape=circle, inner sep=0.5] at (0.8,3) {\small $2$};
		\node[draw,shape=circle, inner sep=0.5] at (3.1,1) {\small $6$};
		\node[draw,shape=circle, inner sep=0.5] at (3,-1) {\small $5$};
		\node[draw,shape=circle, inner sep=0.5] at (1,-1) {\small $4$};
		\node[draw,shape=circle, inner sep=0.5] at (-1,1) {\small $3$};
		\node[draw,shape=circle, inner sep=0.5] at (3.1,3.8) {\small $7$};
		\node[draw,shape=circle, inner sep=0.5] at (5,4.05) {\small $8$};
		\node[draw,shape=circle, inner sep=0.5] at (5.1,3) {\small $9$};
		\node[draw,shape=circle, inner sep=0.5] at (0,5) {\small $10$};
		\node[draw,shape=circle, inner sep=0.5] at (3.3,5.6) {\small $11$};
		\node[draw,shape=circle, inner sep=0.5] at (6,7.5) {\small $12$};
		\node[draw,shape=circle, inner sep=0.5] at (6,6) {\small $13$};
		\node[draw,shape=circle, inner sep=0.5] at (8.2,6) {\small $14$};
		\node[draw,shape=circle, inner sep=0.5] at (7.7,4.1) {\small $15$};
		\node[draw,shape=circle, inner sep=0.5] at (4.9,0.55) {\small $16$};

		\node at (0.35,0.35){\small $\aaa$};
		\node at (-0.35,0.25){\small $\ccc$};
		\node at (0.1,-0.45){\small $\ccc$};
		
		\node at (1.65,0.4){\small $\bbb$};
		\node at (1.65,-0.4){\small $\aaa$};
		\node at (2.35,0.4){\small $\bbb$};
		\node at (2.35,-0.4){\small $\aaa$};
		
		\node at (1.65,1.6){\small $\ddd$};
		\node at (1.65,2.4){\small $\ddd$};
		\node at (2.35,1.6){\small $\ddd$};
		
		\node at (0.35,1.6){\small $\ccc$};
		\node at (0,2.4){\small $\ccc$};
		\node at (-0.35,1.7){\small $\aaa$};    
		
		\node at (3.75,0.4){\small $\aaa$};
		\node at (4.5,1.6){\small $\ccc$};
		
		\node at (1.65,-1.55){\small $\bbb$};
		\node at (1.65,3.45){\small $\bbb$};
		
		\node at (-1.2,3.7) {\small $\aaa$};
		\node at (-1.65,3) {\small $\bbb$};
		\node at (-1.1,-1.6) {\small $\ddd$};
		\node at (-1.7,-1.2) {\small $\ddd$};
		
		\node at (-1.75,4.4) {\small $\aaa$};
		\node at (-2.4,5.5) {\small $\ccc$};
		\node at (1.65,5.5) {\small $\ddd$};
		\node at (1.65,4.45) {\small $\bbb$};
		\node at (2.35,4.6) {\small $\bbb$};
		\node at (2.35,5.6) {\small $\ddd$};
		\node at (3.85,6.4) {\small $\ccc$};
		\node at (4.45,6.05) {\small $\aaa$};
		
		\node at (4.5,7) {\small $\aaa$};
		\node at (5.15,6.6) {\small $\bbb$};
		\node at (6,8.4) {\small $\ccc$};
		\node at (7.2,7.9) {\small $\ddd$};
		\node at (2.35,3.6) {\small $\bbb$};
		\node at (2.3,2.8) {\small $\ddd$};
		\node at (3.4,3.05) {\small $\ddd$};
		\node at (3.4,2.33) {\small $\ddd$};
		
		\node at (3.8,4.25) {\small $\ccc$};
		\node at (4.3,5.25) {\small $\aaa$};
		\node at (5.15,5.6) {\small $\bbb$};
		\node at (4.7,4.8) {\small $\aaa$};
		\node at (4.1,3.6) {\small $\ccc$};

		\node at (8,7.2) {\small $\ddd$};
		\node at (7.3,7.1) {\small $\ddd$};
		\node at (6.8,6.3) {\small $\ccc$};
		\node at (6.3,5.2) {\small $\aaa$};
		\node at (7.3,5.8) {\small $\ccc$};
		\node at (9.2,6) {\small $\bbb$};
		
		\node at (5.85,4.35) {\small $\bbb$};
		\node at (6.7,4.7) {\small $\aaa$};
		\node at (6.6,3.9) {\small $\aaa$};
		\node at (5.9,3.4) {\small $\bbb$};
		
		\node at (6.5,3) {\small $\aaa$};
		\node at (6.1,2.1) {\small $\bbb$};
		\node at (5,2.3) {\small $\ccc$};
		\node at (5.25,1.75) {\small $\aaa$};
		
		\node at (3.65,-0.4) {\small $\ccc$};
		\node at (4.3,-0.25) {\small $\ccc$};
		\node at (2.35,-1.6) {\small $\bbb$};
		\node at (3.65,-1.6) {\small $\ddd$};
		\node at (4.25,-1.1) {\small $\ddd$};
		
		\node at (-2.4,3.8) {\small $\bbb$};
		
		\node at (1.65,-2.5) {\small $\bbb$};
		\node at (2.35,-2.5) {\small $\bbb$};

	\end{tikzpicture}
	\caption{Partial neighborhood of the 2nd special tile.}
	\label{334dF3}
\end{figure}
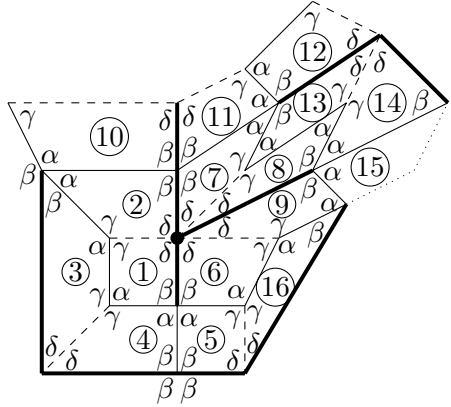

Suppose $H=\ddd_1\ddd_2\ddd_6\cdots=\aaa^k\bbb^l\ccc^m\ddd^n$. Since $\aaa^2\bbb^2$ implies $\aaa+\bbb=\pi$ and $\ccc+\ddd=(1+\frac4f)\pi$, we get $m\le1$.  If $m=1$, $H$ is odd and $l \ge 1$. So we have $\bbb+\ccc+3\ddd \le 2\pi$. By $\aaa^2\bbb^2$ and $\aaa\ccc^2$, we get $3(\aaa+\bbb+\ccc+\ddd)=(2\aaa+2\bbb)+(\aaa+2\ccc)+(\bbb+\ccc+3\ddd)\le6\pi$, contradicting Lemma \ref{anglesum}. 
Therefore, we have $m=0$, $H$ is even and $n \ge 4$. If $l \ge 2$, we have $2\bbb+4\ddd \le 2\pi$. By $\aaa^2\bbb^2$ and $\aaa\ccc^2$, we get $4(\aaa+\bbb+\ccc+\ddd)=(2\aaa+2\bbb)+2(\aaa+2\ccc)+(2\bbb+4\ddd)\le8\pi$, contradicting Lemma \ref{anglesum}. Therefore, we have $l=0$ and $H=\ddd^d$, $d=6,8,10$ by Lemma \ref{ad}. This determines $T_7,T_9$.

The extra angle sum at $H=\ddd^d$, $d=6,8,10$ implies
\[
\aaa=\tfrac{4\pi}{d}-\tfrac{8\pi}{f} ,\quad\bbb=\pi-\tfrac{4\pi}{d}+\tfrac{8\pi}{f},\quad\ccc=\pi-\tfrac{2\pi}{d}+\tfrac{4\pi}{f},\quad\ddd=\tfrac{2\pi}{d}.
\]

Then we have  $\bbb>\frac{\pi}{3},\ccc>\frac{2\pi}{3},\frac{\pi}{5}\le\ddd \le \frac{\pi}{3}$. By the angle values, the edge length consideration and Lemma \ref{lemac2}, we get
\[\text{AVC}\subset\{\aaa\ccc^2,\aaa^2\bbb^2,\ddd^d,\bbb^4,\bbb^2\ddd^2,\aaa^m\}. \]

Then $\aaa_2\bbb_3\cdots=\aaa_2\aaa_{10}\bbb_3\bbb$. By $\aaa_{10}$, 
$\bbb_2\bbb_7\cdots=\bbb_2\bbb_7\bbb_{10}\bbb_{11}$ determines $T_{10},T_{11}$.  By $\aaa\ccc^2$, $\aaa^2\bbb^2$, $\bbb^4$ and $\ddd^{d}$, we get $\aaa=\bbb=\frac{\pi}{2}$, $\ccc=\frac{3\pi}{4}$ and $\ddd=\frac{2\pi}{d}=(\frac{1}{4}+\frac4f)\pi>\frac{\pi}{4}$. So we get $d=6$, $\ddd=\frac{\pi}{3}$ and $f=48$. Therefore 
\[\text{AVC}\sub\{\aaa\ccc^2,\aaa^2\bbb^2,\ddd^6,\bbb^4,\aaa^4\}.\]
Then $H=\ddd^6$ determines $T_8$, and  $\ccc_7\ccc_8\cdots=\aaa_{13}\ccc_7\ccc_8$. By $\aaa_{13}$, 
$\aaa_7\aaa_{11}\cdots=\aaa_7\aaa_{11}\bbb_{12}\bbb_{13}$ determines $T_{12},T_{13}$; $\aaa_8\ccc_{13}\cdots=\aaa_8\ccc_{13}\ccc_{14}$ determines $T_{14}$. So $\aaa_{14}\bbb_8\bbb_9\cdots=\aaa_{14}\aaa_{15}\bbb_8\bbb_9$. By $\aaa_5$, we get $\thick\bbb_4\thin\cdots=\thick\bbb_4\thin\bbb_5\thick\cdots=\bbb^4$, which determines $T_{5}$. Then $\aaa_6\ccc_5\cdots=\aaa_6\ccc_5\ccc_{16}$ and  $\aaa_{16}\ccc_6\ccc_9\cdots=\aaa_{16}\ccc_6\ccc_9$ determine $T_{16}$. By $\aaa_{15}$, we have $\aaa_9\bbb_{16}\cdots \neq \aaa^2\bbb^2$,  contradicting the AVC.
\end{proof}

\begin{proposition} \label{335dCase1}
	There is no $a^2bc$-tiling with the 3rd special tile in Fig. \ref{33--}.
\end{proposition}

\begin{proof}
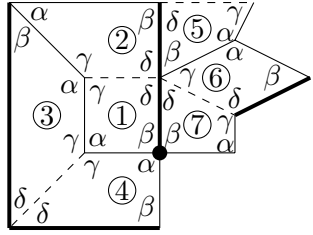
\begin{figure}[htp]
	\centering
	\begin{tikzpicture}[>=latex,scale=0.5]       
		\draw (0,0) -- (0,2) 
		(0,0) -- (2,0)
		(0,2)--(-2,4)
		(-2,4)--(2,4)
		(2,0)--(2,-2)
		(4,0)--(4,1)
		(4,3)--(4.5,4)
		(2,2)--(4,3)
		(6,2)--(4,3)
		(2,0)--(4,0);
		\draw[dashed]  (0,2)--(2,2)
		(2,2)--(4,1)
		(2,4)--(4.5,4)
		(0,0)--(-2,-2);
		\draw[line width=1.5] (2,0)--(2,2)
		(2,2)--(2,4)
		(6,2)--(4,1)
		(-2,-2)--(-2,4)
		(-2,-2)--(2,-2);

		\fill (2,0) circle (0.2);
		
		\node[draw,shape=circle, inner sep=0.5] at (1,1) {\small $1$};
		\node[draw,shape=circle, inner sep=0.5] at (1,3) {\small $2$};
		\node[draw,shape=circle, inner sep=0.5] at (3,3.3) {\small $5$};
		\node[draw,shape=circle, inner sep=0.5] at (3.5,2) {\small $6$};
		\node[draw,shape=circle, inner sep=0.5] at (3,0.75) {\small $7$};
		\node[draw,shape=circle, inner sep=0.5] at (1,-1) {\small $4$};
		\node[draw,shape=circle, inner sep=0.5] at (-1,1) {\small $3$};
		
		\node at (0.35,0.35) {\small $\aaa$};
		\node at (0.2,-0.45) {\small $\ccc$};
		\node at (-0.35,0.25) {\small $\ccc$};
		
		\node at (1.65,0.45) {\small $\bbb$};
		\node at (1.65,-0.35) {\small $\aaa$};
		
		\node at (1.65,1.55) {\small $\ddd$};
		\node at (1.65,2.45) {\small $\ddd$};
		\node at (2.35,2.65) {\small $\bbb$};
		\node at (2.25,1.4) {\small $\ddd$};
		\node at (2.75,2) {\small $\ccc$};
		
		\node at (2.35,3.5) {\small $\ddd$};
		\node at (1.65,3.5) {\small $\bbb$};
		\node at (2.35,0.45) {\small $\bbb$};
		
		\node at (4,3.65) {\small $\ccc$};
		\node at (3.75,3.15) {\small $\aaa$};
		\node at (4,2.65) {\small $\aaa$};
		\node at (5,2) {\small $\bbb$};
		\node at (4,1.45) {\small $\ddd$};
		\node at (3.75,0.75) {\small $\ccc$};
		\node at (3.75,0.2) {\small $\aaa$};
		
		\node at (0,2.4) {\small $\ccc$};
		\node at (0.35,1.55) {\small $\ccc$};
		\node at (-0.4,1.75) {\small $\aaa$};
		
		\node at (-1.2,3.7) {\small $\aaa$};
		\node at (-1.65,3) {\small $\bbb$};
		\node at (-1.1,-1.6) {\small $\ddd$};
		\node at (-1.7,-1.2) {\small $\ddd$};
		\node at (1.65,-1.5) {\small $\bbb$};    
	\end{tikzpicture}
	\caption{Partial neighborhood of the 3rd special tile.} \label{335dsubcase4}
\end{figure}
Let the third of Fig. \ref{33--} be the center tile $T_1$ in the partial neighborhoods in Fig.  \ref{335dsubcase4}.
By the edge length consideration, we get $\ddd_1\ddd_2\cdots=\aaa\bbb^2\ddd^2$ or $\bbb\ccc\ddd^3$. If $\aaa\bbb^2\ddd^2$ is a vertex, by $\aaa\ccc^2$, then $2(\aaa+\bbb+\ccc+\ddd)=4\pi$, contradicting Lemma \ref{anglesum}. So we have $\ddd_1\ddd_2\cdots=\bbb\ccc\ddd^3$. By Lemma \ref{lemac2} and $\aaa\ccc^2$, we get $\aaa_4\bbb_1\cdots=\aaa^k\bbb^{2t}$. This determines $T_7,T_6,T_5$. If $k\ge2$, then $2\aaa+2\bbb\le2\pi$.  By $\aaa\ccc^2$ and $\bbb\ccc\ddd^3$, we get $3(\aaa+\bbb+\ccc+\ddd)=(\aaa+2\ccc)+(\bbb+\ccc+3\ddd)+(2\aaa+2\bbb)\le6\pi$, contradicting Lemma \ref{anglesum}. Therefore, $H=\aaa\bbb^{d-1}\,(d=5,7)$.

By Lemma \ref{lemd4}, $\aaa^2\cdots$ is a vertex and must be even.
By $\aaa\ccc^2,\bbb\ccc\ddd^3$ and Lemma \ref{anglesum}, we get $\aaa+\bbb=(1+\frac{6}{f})\pi>\pi$ and $\aaa+\ccc>\pi$.  By $R(\aaa^2\cdots)<2\bbb,2\ccc$ and Lemma \ref{ad}, we get $\aaa^2\cdots=\aaa^k$. The angle sums at $\aaa^k,\aaa\bbb^{d-1},\aaa\ccc^2,\bbb\ccc\ddd^3$ imply $\aaa=\frac{2}{k}\pi,\bbb=\frac{2k-2}{k(d-1)}\pi,\ccc=(1-\frac1k)\pi,\ddd=(\frac13+\frac{1}{3k}-\frac{2k-2}{3k(d-1)})\pi$. Then $\aaa+\bbb+\ccc+\ddd=\frac{4}{3}(1+\frac1k+\frac{k-1}{k(d-1)})\pi\le2\pi$ by $k\ge3,d=5,7$, contradicting Lemma \ref{anglesum}. 
\end{proof}

\begin{proposition} \label{335dCase2}
	There is no $a^2bc$-tiling with the 4th special tile in Fig. \ref{33--}.
\end{proposition}
\begin{proof}

Let the fourth of Fig. \ref{33--} be the center tile $T_1$ in the partial neighborhoods in Fig.  \ref {335dsubcase3.4.1}.
By the edge length consideration and Lemma \ref{lemac2}, we get $\aaa_4\bbb_1\cdots=\aaa^3\bbb^2$ or $\aaa\bbb^4$, shown in the first and second pictures of  Fig.  \ref {335dsubcase3.4.1}. 

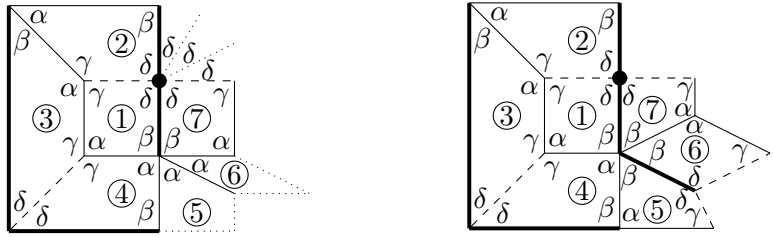
\begin{figure}[htp]
	\centering
	\begin{tikzpicture}[>=latex,scale=0.5]       
		\draw (0,0) -- (0,2) 
		(0,0)--(2,0)
		(0,2)--(-2,4)
		(-2,4)--(2,4)
		(2,0)--(2,-2)
		(2,0)--(4,0)
		(2,0)--(4,-1)
		(4,0)--(4,2);
		\draw[dashed]  (0,2)--(2,2)
		(2,2)--(4,2)
		(0,0)--(-2,-2);
		\draw[line width=1.5] (2,0)--(2,2)
		(-2,-2)--(-2,4)
		(-2,-2)--(2,-2)
		(2,2)--(2,4);
		\draw[dotted] (2,-2)--(4,-2)
		(4,-1)--(6,-1)
		(4,-2)--(4,-1)
		(4,0)--(6,-1)
		(2,2)--(3.1,3.8)
		(2,2)--(4,3);

		\fill (2,2) circle (0.2);
		
		\node[draw,shape=circle, inner sep=0.5] at (1,1) {\small $1$};
		\node[draw,shape=circle, inner sep=0.5] at (1,3) {\small $2$};
		\node[draw,shape=circle, inner sep=0.5] at (3,1) {\small $7$};
		\node[draw,shape=circle, inner sep=0.5] at (4,-0.5) {\small $6$};
		\node[draw,shape=circle, inner sep=0.5] at (3,-1.5) {\small $5$};
		\node[draw,shape=circle, inner sep=0.5] at (1,-1) {\small $4$};
		\node[draw,shape=circle, inner sep=0.5] at (-1,1) {\small $3$};
		
		\node at (0.35,0.35) {\small $\aaa$};
		\node at (0.2,-0.45) {\small $\ccc$};
		\node at (-0.35,0.25) {\small $\ccc$};
		
		\node at (1.65,0.45) {\small $\bbb$};
		\node at (1.65,-0.35) {\small $\aaa$};
		
		\node at (1.65,1.55) {\small $\ddd$};
		\node at (1.65,2.45) {\small $\ddd$};
		
		\node at (1.65,3.5) {\small $\bbb$};
		
		\node at (0,2.4) {\small $\ccc$};
		\node at (0.35,1.55) {\small $\ccc$};
		\node at (-0.4,1.75) {\small $\aaa$};
		
		\node at (-1.2,3.7) {\small $\aaa$};
		\node at (-1.65,3) {\small $\bbb$};
		\node at (-1.1,-1.6) {\small $\ddd$};
		\node at (-1.7,-1.2) {\small $\ddd$};
		\node at (1.65,-1.5) {\small $\bbb$};
		
		\node at (2.35,0.4) {\small $\bbb$}; \node at (3.65,0.35) {\small $\aaa$};
		\node at (2.35,1.55) {\small $\ddd$}; \node at (3.65,1.55) {\small $\ccc$};
		
		\node at (3.1,-0.25) {\small $\aaa$}; \node at (2.35,-0.5) {\small $\aaa$}; 
		\node at (3.3,2.3) {\small $\ddd$}; \node at (2.8,2.75) {\small $\ddd$};
		\node at (2.25,2.9) {\small $\ddd$};

	\end{tikzpicture}	\hspace{50pt}
	\begin{tikzpicture}[>=latex,scale=0.5]       
		\draw (0,0) -- (0,2) 
		(0,0)--(2,0)
		(0,2)--(-2,4)
		(-2,4)--(2,4)
		(2,0)--(2,-2)
		(2,0)--(4,1)
		(4,1)--(6,0)
		(4,1)--(4,2)
		(2,-2)--(4.5,-2);
		\draw[dashed]  (0,2)--(2,2)
		(0,0)--(-2,-2)
		(2,2)--(4,2)
		(4,-1)--(6,0)
		(4.5,-2)--(4,-1);
		\draw[line width=1.5] (2,0)--(2,2)
		(2,2)--(2,4)
		(-2,-2)--(-2,4)
		(-2,-2)--(2,-2)
		(2,0)--(4,-1);

		\fill (2,2) circle (0.2);
		
		\node[draw,shape=circle, inner sep=0.5] at (1,1) {\small $1$};
		\node[draw,shape=circle, inner sep=0.5] at (1,3) {\small $2$};
		\node[draw,shape=circle, inner sep=0.5] at (2.87,1.15) {\small $7$};
		\node[draw,shape=circle, inner sep=0.5] at (4,0.1) {\small $6$};
		\node[draw,shape=circle, inner sep=0.5] at (3,-1.5) {\small $5$};
		\node[draw,shape=circle, inner sep=0.5] at (1,-1) {\small $4$};
		\node[draw,shape=circle, inner sep=0.5] at (-1,1) {\small $3$};

		\node at (0.35,0.35) {\small $\aaa$};
		\node at (0.2,-0.45) {\small $\ccc$};
		\node at (-0.35,0.25) {\small $\ccc$};
		
		\node at (1.65,0.45) {\small $\bbb$};
		\node at (1.65,-0.35) {\small $\aaa$};
		
		\node at (1.65,1.55) {\small $\ddd$};
		\node at (1.65,2.45) {\small $\ddd$};
		\node at (2.25,1.55) {\small $\ddd$};
		
		\node at (1.65,3.5) {\small $\bbb$};
		
		\node at (0,2.4) {\small $\ccc$};
		\node at (0.35,1.55) {\small $\ccc$};
		\node at (-0.4,1.75) {\small $\aaa$};
		
		\node at (-1.2,3.7) {\small $\aaa$};
		\node at (-1.65,3) {\small $\bbb$};
		\node at (-1.1,-1.6) {\small $\ddd$};
		\node at (-1.7,-1.2) {\small $\ddd$};
		\node at (1.65,-1.5) {\small $\bbb$};

		\node at (2.35,0.5){\small $\bbb$};
		\node at (2.25,-0.65){\small $\bbb$};
		\node at (3,0){\small $\bbb$};
		
		\node at (3.7,1.15) {\small $\aaa$}; \node at (3.75,1.65) {\small $\ccc$};
		\node at (4,0.7) {\small $\aaa$};
		
		\node at (3.65,-1.23) {\small $\ddd$};   \node at (5.2,0) {\small $\ccc$};
		\node at (4,-0.6) {\small $\ddd$};
		
		\node at (2.3,-1.65){\small $\aaa$}; \node at (3.95,-1.7){\small $\ccc$};

	\end{tikzpicture}
	\caption{Partial neighborhoods of the 4th special tile.}\label{335dsubcase3.4.1}
\end{figure}
\subsubsection*{Case $\aaa_4\bbb_1\cdots=\aaa^3\bbb^2$}

This $\thin\aaa_4\thin\bbb_1\thick\cdots=\thin\bbb_1\thick\bbb_7\thin\aaa_6\thin\aaa_5\thin\aaa_4\thin$ determines $T_7$.
Suppose $H=\ddd_1\ddd_2\ddd_7\cdots=\aaa^k\bbb^l\ccc^m\ddd^n$. If $m\ge1$, by Parity Lemma, we get $H=\bbb\ccc\ddd^3\cdots$ or $\ccc^2\ddd^4\cdots$. But by $\aaa\ccc^2,\aaa^3\bbb^2$,  
\begin{align*}
3(\aaa+\bbb+\ccc+\ddd)&<(\aaa+2\ccc)+(3\aaa+2\bbb)+(\bbb+\ccc+3\ddd)\le 6\pi, \\ 2(\aaa+\bbb+\ccc+\ddd)&<(3\aaa+2\bbb)+(2\ccc+4\ddd)\le 4\pi.
\end{align*} 
Both $H$ contradict Lemma \ref{anglesum}. Therefore, $m=0$ and $H$ must be even. If $l>0$, we get $H=\bbb^2\ddd^4\cdots$. But  \[4(\aaa+\bbb+\ccc+\ddd)<2(\aaa+2\ccc)+(3\aaa+2\bbb)+(2\bbb+4\ddd)\le8\pi,\]
which contradicts Lemma \ref{anglesum}. Therefore, $l=0$. By Lemma \ref{ad}, we have $H=\ddd^6$. By $\aaa\ccc^2,\aaa^3\bbb^2,\ddd^6$, we get
$\aaa=(\frac{1}{3}-\frac{4}{f})\pi ,\bbb=(\frac{1}{2}+\frac{6}{f})\pi,\ccc=(\frac{5}{6}+\frac{2}{f})\pi,\ddd=\frac{\pi}{3}$. 

The AAD of $\thin\aaa_6\thin\aaa_5\thin\aaa_4\thin\cdots$ induces $\ccc_5\thin\ccc_6\cdots$ or $\bbb\thin\ccc_5\cdots$.	
But $R(\dash\ccc\thin\ccc\dash\cdots)<\frac{\pi}{3}\le\ccc$, $\ddd$, a contradiction. Similarly $R(\thick\bbb\thin\ccc\dash\cdots)<\frac{2\pi}{3}<\ccc$ implies $\thick\bbb\thin\ccc\dash\cdots=\thick\bbb\thin\ccc\dash\ddd\thick\cdots$, whose remainder is $(\frac{1}{3}-\frac{8}{f})\pi<$ all angles, a contradiction.

\subsubsection*{Case $\aaa_4\bbb_1\cdots=\aaa\bbb^4$}

This $\thin\aaa_4\thin\bbb_1\thick\cdots=\aaa\bbb^4$ determines $T_5,T_6,T_7$.
By $\aaa\ccc^2$, $\aaa\bbb^4$ and Lemma \ref{lemd4}, we get
$\pi>\ccc=2\bbb=\pi-\frac{\aaa}{2}>\frac{\pi}{2}$, $\ddd=\frac{\pi}{2}-\frac{\aaa}{4}+\frac{4\pi}{f}>\frac{\pi}{4}$. If $\text{deg} H=7$, then $H$ is odd by Lemma \ref{lemac2}, and $H=\ddd_1\ddd_2\ddd_7\cdots=\bbb^3\ccc\ddd^3,\bbb\ccc^3\ddd^3$ or $\bbb\ccc\ddd^5$. But by $\bbb,\ddd>\frac{\pi}{4}$ and $\ccc>\frac{\pi}{2}$, all of them are $>2\pi$, a contradiction. So $\text{deg} H=6$ and we get $H=\bbb^2\ddd^4,\ccc^2\ddd^4$ or $\ddd^6$ similarly. If $H=\ccc^2\ddd^4$, it contradicts $2\ccc+4\ddd>2\pi$.
If $H=\ddd^6$, we get $\aaa=(\frac{2}{3}+\frac{16}{f})\pi,\bbb=(\frac{1}{3}-\frac{4}{f})\pi,\ccc=(\frac{2}{3}-\frac{8}{f})\pi,\ddd=\frac{\pi}{3}$. 
If $H=\bbb^2\ddd^4$, we get  $\aaa=(\frac{2}{3}+\frac{32}{3f})\pi,\bbb=(\frac{1}{3}-\frac{8}{3f})\pi,\ccc=(\frac{2}{3}-\frac{16}{3f})\pi,\ddd=(\frac{1}{3}+\frac{4}{3f})\pi$. Both imply $R(\aaa_6\aaa_7\cdots)<\aaa,2\bbb$, contradicting Lemma \ref{lemac2}.
\end{proof}

In summary, there is no $a^2bc$-tiling with any $335d$-Tile.

\section{$344d$-Tile} \label{344d}
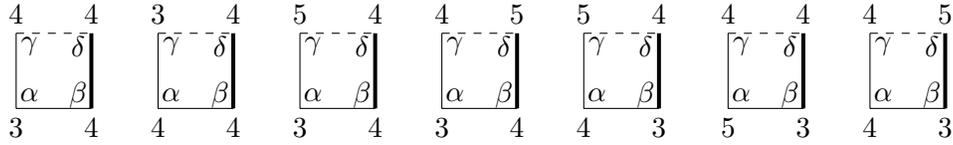
\begin{figure}[htp]
	\centering
	\begin{tikzpicture}[>=latex,scale=0.5]
		\draw
		(0,0)--(0,2)
		(0,0)--(2,0);
		\draw[line width=1.5]
		(2,0)--(2,2);
		\draw[dashed]
		(0,2)--(2,2);
		\node at (0.35,0.35){\small $\aaa$};
		\node at (1.65,0.35){\small $\bbb$};
		\node at (1.65,1.65){\small $\ddd$};
		\node at (0.35,1.65){\small $\ccc$};
		\node at (0,-0.5) {\small $3$};
		\node at (0,2.5) {\small $4$};
		\node at (2,-0.5) {\small $4$};
		\node at (2,2.5) {\small $4$};
	\end{tikzpicture}\hspace{8pt}
	\begin{tikzpicture}[>=latex,scale=0.5]
		\draw
		(0,0)--(0,2)
		(0,0)--(2,0);
		\draw[line width=1.5]
		(2,0)--(2,2);
		\draw[dashed]
		(0,2)--(2,2);
		\node at (0.35,0.35){\small $\aaa$};
		\node at (1.65,0.35){\small $\bbb$};
		\node at (1.65,1.65){\small $\ddd$};
		\node at (0.35,1.65){\small $\ccc$};
		\node at (0,-0.5) {\small $4$};
		\node at (0,2.5) {\small $3$};
		\node at (2,-0.5) {\small $4$};
		\node at (2,2.5) {\small $4$};
	\end{tikzpicture}\hspace{8pt}
	\begin{tikzpicture}[>=latex,scale=0.5]      
		\draw (0,0) -- (0,2) 
		(0,0) -- (2,0);
		\draw[dashed]  (0,2)--(2,2);
		\draw[line width=1.5] (2,0)--(2,2);
		\node at (0.35,0.35){\small $\aaa$};
		\node at (1.65,0.35){\small $\bbb$};
		\node at (1.65,1.65){\small $\ddd$};
		\node at (0.35,1.65){\small $\ccc$};
		\node at (0,-0.5) {\small $3$};          
		\node at (0,2.5) {\small $5$};
		\node at (2,-0.5) {\small $4$};          
		\node at (2,2.5) {\small $4$};
		
	\end{tikzpicture}\hspace{8pt}
	\begin{tikzpicture}[>=latex,scale=0.5]      
		\draw (0,0) -- (0,2) 
		(0,0) -- (2,0);
		\draw[dashed]  (0,2)--(2,2);
		\draw[line width=1.5] (2,0)--(2,2);
		\node at (0.35,0.35){\small $\aaa$};
		\node at (1.65,0.35){\small $\bbb$};
		\node at (1.65,1.65){\small $\ddd$};
		\node at (0.35,1.65){\small $\ccc$};
		\node at (0,-0.5) {\small $3$};          
		\node at (0,2.5) {\small $4$};
		\node at (2,-0.5) {\small $4$};          
		\node at (2,2.5) {\small $5$};
		
	\end{tikzpicture}\hspace{8pt}
	\begin{tikzpicture}[>=latex,scale=0.5]      
		\draw (0,0) -- (0,2) 
		(0,0) -- (2,0);
		\draw[dashed]  (0,2)--(2,2);
		\draw[line width=1.5] (2,0)--(2,2);
		\node at (0.35,0.35){\small $\aaa$};
		\node at (1.65,0.35){\small $\bbb$};
		\node at (1.65,1.65){\small $\ddd$};
		\node at (0.35,1.65){\small $\ccc$};
		\node at (0,-0.5) {\small $4$};          
		\node at (0,2.5) {\small $5$};
		\node at (2,-0.5) {\small $3$};          
		\node at (2,2.5) {\small $4$};
	\end{tikzpicture} \hspace{5pt}
	\begin{tikzpicture}[>=latex,scale=0.5]      
		\draw (0,0) -- (0,2) 
		(0,0) -- (2,0);
		\draw[dashed]  (0,2)--(2,2);
		\draw[line width=1.5] (2,0)--(2,2);
		\node at (0.35,0.35){\small $\aaa$};
		\node at (1.65,0.35){\small $\bbb$};
		\node at (1.65,1.65){\small $\ddd$};
		\node at (0.35,1.65){\small $\ccc$};
		\node at (0,-0.5) {\small $5$};          
		\node at (0,2.5) {\small $4$};
		\node at (2,-0.5) {\small $3$};          
		\node at (2,2.5) {\small $4$};
		
	\end{tikzpicture}\hspace{8pt}
	\begin{tikzpicture}[>=latex,scale=0.5]      
		\draw (0,0) -- (0,2) 
		(0,0) -- (2,0);
		\draw[dashed]  (0,2)--(2,2);
		\draw[line width=1.5] (2,0)--(2,2);
		\node at (0.35,0.35){\small $\aaa$};
		\node at (1.65,0.35){\small $\bbb$};
		\node at (1.65,1.65){\small $\ddd$};
		\node at (0.35,1.65){\small $\ccc$};
		\node at (0,-0.5) {\small $4$};          
		\node at (0,2.5) {\small $4$};
		\node at (2,-0.5) {\small $3$};          
		\node at (2,2.5) {\small $5$};
		
	\end{tikzpicture}

	\caption{Special tiles with vertex degrees $344d$, $d=4,5$.}
	\label{3444-5}
\end{figure}

By Proposition \ref{333d}, we can assume $\delta$ never appears at any degree $3$ vertices. 
By the symmetry of exchanging $\bbb\leftrightarrow\ccc$, there are only $7$ different configurations for degree $3,4,5$ vertices in $3444$-Tile and $3445$-Tile in Fig. \ref{3444-5}. We first prove two useful propositions before studying each special tile.

\begin{proposition} \label{proposition1}
	There is no $a^2bc$-tiling with vertices $\aaa^3,\bbb^4,\ccc^2\ddd^2$.	
\end{proposition}

\begin{proof}
	If $\aaa^3,\bbb^4,\ccc^2\ddd^2$ are vertices, we have $\aaa=\frac{2\pi}{3},\bbb=\frac{\pi}{2},\ccc+\ddd=\pi$. By Lemma \ref{geometry4}, we have $\ccc>\frac{\pi}{6}$. The AAD of $\aaa^3=\thin\aaa^{\bbb}\thin^{\ccc}\aaa\thin\cdots$ gives a vertex $\bbb\ccc\cdots=\aaa^k\bbb^l\ccc^m\ddd^n$. If $k\ge1$, we have $\bbb\ccc\cdots=\aaa\bbb\ccc\cdots$ must be even. So $l\ge2,m\ge2$. Then we have $\aaa+2\bbb+2\ccc>2\pi$, a contradiction. Thus $k=0$. 
	
	If $n=0$, it is even. By $\bbb=\frac12\pi,\ccc>\frac16\pi$, we get $\bbb\ccc\cdots=\bbb^2\ccc^2$ or $\bbb^2\ccc^4$. But $\bbb^2\ccc^2$ and $\bbb^4$ imply $\bbb=\ccc$, contradicting Lemma \ref{lemd4}. So $\bbb\ccc\cdots=\bbb^2\ccc^4=\thin^{\aaa}\ccc^{\ddd}\dash^{\ddd}\ccc^{\aaa}\thin\cdots$. This gives a vertex $\thick\ddd\dash\ddd\thick\cdots$. By $\bbb^2\ccc^4$, we have $\ddd=\frac{3\pi}{4},\bbb+\ddd>\pi$, contradicting Lemma \ref{lembd}.
	
	If $n>0$, by $\bbb^4,\ccc^2\ddd^2$, Lemma \ref{lemd4} and Parity Lemma, we get $\bbb\ccc\cdots=\bbb\ccc^m\ddd$ or $\bbb\ccc\ddd^n$ with $m,n$ being odd integers $\ge3$. 
	
	If $\bbb\ccc^m\ddd$ is a vertex,
	by $\bbb^4,\ccc^2\ddd^2,\bbb^m\ddd$, we get $\ccc=\frac{\pi}{2(m-1)}\le \frac\pi4$ and $\ddd\ge \frac34\pi$. Its AAD $\bbb\ccc^m\ddd=\ccc^{\ddd}\dash^{\ddd}\ccc\cdots$  gives a vertex $\ddd\dash\ddd\cdots$, contradicting Lemma \ref{lembd}. 
	
	If $\bbb\ccc\ddd^n$ is a vertex, then similarly we get $\ccc\ge\frac{3\pi}{4},\ddd\le\frac{\pi}{4}$. Its AAD $\bbb\ccc\ddd^n=\dash^{\ccc}\ddd^{\bbb}\thick^{\bbb}\ddd^{\ccc}\dash^{\ccc}\ddd^{\bbb}\thick\cdots$ gives a vertex $\thin\ccc\dash\ccc\thin\cdots=\theta\thin\ccc\dash\ccc\thin\rho\cdots$, where $\theta,\rho=\aaa,\bbb$ or $\ccc$. By $\aaa=\frac{2\pi}{3},\bbb=\frac{\pi}{2},\ccc\ge\frac{3\pi}{4}$, we get a contradiction. 
\end{proof}

\begin{proposition} \label{proposition2}
	There is no $a^2bc$-tiling with vertices $\aaa^3,\bbb^2\ccc^2,\ccc^2\ddd^2,\bbb\thin\bbb\cdots$. 
\end{proposition}

\begin{proof}
	If $\aaa^3,\bbb^2\ccc^2$ and $\ccc^2\ddd^2$ are vertices, by Lemma \ref{anglesum}, we have $\aaa=\frac{2\pi}{3},\ccc=(\frac{2}{3}-\frac{4}{f})\pi,\bbb=\ddd=(\frac{1}{3}+\frac{4}{f})\pi$. Let $\bbb\thin\bbb\cdots=\aaa^k\bbb^l\ccc^m\ddd^n$.  By $\bbb^2\ccc^2=\thin\bbb\thick\bbb\thin\ccc\dash\ccc\thin$, we get $m\le1$. If $m=1$, then it is odd and $l\ge3,n\ge1$. However, we have $3\bbb+\ccc+\ddd=2\pi+\frac{12\pi}{f}>2\pi$, a contradiction. Therefore $m=0$. By edge length consideration (Parity Lemma is not enough), we have $l+n\ge4$. By $\aaa+4\bbb=\aaa+2\bbb+2\ddd=(2+\frac{16}{f})\pi>2\pi$, we get $k=0$. So we have $\bbb\thin\bbb\cdots=\bbb^l\ddd^n$, where $l,n$ are even. If $l+n\ge6$, we get $l\bbb+n\ddd=(l+n)(\frac{\pi}{3}+\frac{4\pi}{f})>2\pi$, a contradiction. Therefore $l+n=4$. By $\bbb=\ddd$ and $\bbb^2\ccc^2,\bbb^l\ddd^n$, we get $\bbb=\ccc=\ddd$,  contradicting Lemma \ref{lemd4}. 
\end{proof}

\begin{proposition}\label{case1111}
	All $a^2bc$-tilings with the $1st$ special tile in Fig. \ref{3444-5} are quadrilateral subdivisions of the octahedron $T(8\aaa^3, 6\ddd^4, 12\bbb^2\ccc^2)$ with $24$ tiles.
\end{proposition}
\begin{proof}
	
	\begin{figure}[htp]
		\centering
		\begin{tikzpicture}[>=latex,scale=0.45]
			\draw (0,0)--(2,0)
			(0,0)--(0,2)
			
			(0,2)--(-2,2)
			
			;
			\draw[line width=1.5]
			(2,0)--(2,2)
			
			(-2,-2)--(-2,2);
			\draw[dashed]
			(0,2)--(2,2)
			(0,0)--(-2,-2);
			\draw[dotted]
			(-2,-2)--(2,-2)
			(2,0)--(2,-2)
			(0,4)--(2,4)
			(2,2)--(2,4)
			(0,2)--(0,4)
			(0,4)--(-2,4)--(-2,2);
			\node at (0.3,0.3) {\small $\aaa$};
			\node at (0.3,1.6) {\small $\ccc$};
			\node at (-0.3,1.65) {\small $\aaa$};
			\node at (1.6,0.4) {\small $\bbb$};
			\node at (1.65,1.6) {\small $\ddd$};
			\node at (-0.35,0.15) {\small $\ccc$};
			\node at (-1.65,-1.2) {\small $\ddd$};
			\node at (0.2,-0.35) {\small $\ccc$};
			\node at (-1.65,1.55) {\small $\bbb$};
			\node[draw,shape=circle, inner sep=0.5] at (1,1) {\small $1$};
			\node[draw,shape=circle, inner sep=0.5] at (-1,1) {\small $2$};
			\node[draw,shape=circle, inner sep=0.5] at (-1,3) {\small $3$};
			\node[draw,shape=circle, inner sep=0.5] at (1,3) {\small $4$};
			\node[draw,shape=circle,inner sep=0.5] at (1,-1){\small $8$};
			\node at (0,-3){$E_{28}= c$};
		\end{tikzpicture}\hspace{13pt}
		\begin{tikzpicture}[>=latex,scale=0.45]
			\draw (0,0)--(2,0)
			(0,0)--(0,2)
			(0,0)--(-2,-2)
			(2,4)--(4,4)
			(4,2)--(4,4)
			(4,0)--(4,2)
			(2,0)--(4,0)
			(-2,4)--(2,4)		
			(0,2)--(0,4);
			\draw[line width=1.5]
			(-2,2)--(-2,4)
			(-2,-2)--(-2,2)
			(2,0)--(2,2)			
			(2,2)--(2,4);
			
			\draw[dashed]
			(-2,2)--(2,2)				
			(2,2)--(4,2);
			\draw[dotted]
			(-2,-2)--(2,-2)
			(2,0)--(2,-2)
			(2,-2)--(4,-2)
			(4,0)--(4,-2);
			\node[draw,shape=circle, inner sep=0.5] at (1,1) {\small $1$};
			\node[draw,shape=circle, inner sep=0.5] at (-1,3) {\small $3$};
			\node[draw,shape=circle, inner sep=0.5] at (-1,1) {\small $2$};
			\node[draw,shape=circle, inner sep=0.5] at (1,-1) {\small $8$};
			\node[draw,shape=circle, inner sep=0.5] at (3,1) {\small $6$};
			\node[draw,shape=circle, inner sep=0.5] at (3,-1) {\small $7$};
			\node[draw,shape=circle, inner sep=0.5] at (3,3) {\small $5$};
			\node[draw,shape=circle, inner sep=0.5] at (1,3) {\small $4$};
			\node at (0.35,0.3) {\small $\aaa$};
			\node at (0.35,1.6) {\small $\ccc$};
			\node at (0.35,2.4) {\small $\ccc$};
			\node at (0.35,3.65) {\small $\aaa$};
			\node at (0.2,-0.3) {\small $\aaa$};
			\node at (-0.35,0.15) {\small $\aaa$};
			\node at (-0.35,1.6) {\small $\ccc$};
			\node at (-0.35,2.4) {\small $\ccc$};
			\node at (-0.35,3.65) {\small $\aaa$};
			\node at (-1.65,-1.2) {\small $\bbb$};
			\node at (-1.65,1.6) {\small $\ddd$};
			\node at (-1.65,2.4) {\small $\ddd$};
			\node at (-1.65,3.5) {\small $\bbb$};
			\node at (1.65,0.4) {\small $\bbb$};
			\node at (1.65,1.6) {\small $\ddd$};
			\node at (1.65,2.4) {\small $\ddd$};
			\node at (1.65,3.5) {\small $\bbb$};
			\node at (2.35,0.4) {\small $\bbb$};
			\node at (2.35,1.6) {\small $\ddd$};
			\node at (2.35,2.4) {\small $\ddd$};
			\node at (2.35,3.5) {\small $\bbb$};
			\node at (3.7,0.3) {\small $\aaa$};
			\node at (3.7,1.65) {\small $\ccc$};
			\node at (3.7,2.35) {\small $\ccc$};
			\node at (3.7,3.55) {\small $\aaa$};
			\node at (1,-3) {\small $E_{23}=c,E_{34}=a$};
		\end{tikzpicture}\hspace{10pt}
		\begin{tikzpicture}[>=latex,scale=0.45]
			\draw (0,0)--(2,0)
			(0,0)--(0,2)
			(0,0)--(-2,-2)
			(0,0)--(4,0)
			(4,0)--(4,2)
			(-2,2)--(-2,4)
			(-2,4)--(4,4)
			(2,2)--(2,4)
			(4,0)--(4,-2);
			\draw[line width=1.5]
			(0,2)--(0,4)
			(-2,-2)--(-2,2)
			(2,0)--(2,2)
			(4,2)--(4,4)
			(-2,-2)--(4,-2);
			\draw[dashed]
			(-2,2)--(4,2)
			(2,0)--(2,-2);

			\node[draw,shape=circle, inner sep=0.5] at (1,1) {\small $1$};
			\node[draw,shape=circle, inner sep=0.5] at (-1,3) {\small $3$};
			\node[draw,shape=circle, inner sep=0.5] at (-1,1) {\small $2$};
			\node[draw,shape=circle, inner sep=0.5] at (1,-1) {\small $8$};
			\node[draw,shape=circle, inner sep=0.5] at (3,1) {\small $6$};
			\node[draw,shape=circle, inner sep=0.5] at (3,-1) {\small $7$};
			\node[draw,shape=circle, inner sep=0.5] at (3,3) {\small $5$};
			\node[draw,shape=circle, inner sep=0.5] at (1,3) {\small $4$};
			\node at (1.65,-0.4) {\small $\ccc$};
			\node at (2.3,-0.4) {\small $\ccc$};
			\node at (1.65,-1.55) {\small $\ddd$};
			\node at (2.3,-1.55) {\small $\ddd$};
			\node at (-1.05,-1.5) {\small $\bbb$};
			\node at (3.65,-1.55) {\small $\bbb$};
			\node at (3.65,-0.3) {\small $\aaa$};
			\node at (0.35,0.3) {\small $\aaa$};
			\node at (0.35,1.6) {\small $\ccc$};
			\node at (0.35,2.4) {\small $\ddd$};
			\node at (0.35,3.5) {\small $\bbb$};
			\node at (0.2,-0.3) {\small $\aaa$};
			\node at (-0.35,0.15) {\small $\aaa$};
			\node at (-0.35,1.6) {\small $\ccc$};
			\node at (-0.35,2.4) {\small $\ddd$};
			\node at (-0.35,3.5) {\small $\bbb$};
			\node at (-1.65,-1.2) {\small $\bbb$};
			\node at (-1.65,1.6) {\small $\ddd$};
			\node at (-1.65,2.4) {\small $\ccc$};
			\node at (-1.65,3.65) {\small $\aaa$};
			\node at (1.65,0.4) {\small $\bbb$};
			\node at (1.65,1.6) {\small $\ddd$};
			\node at (1.65,2.4) {\small $\ccc$};
			\node at (1.65,3.65) {\small $\aaa$};
			\node at (2.3,0.4) {\small $\bbb$};
			\node at (2.3,1.6) {\small $\ddd$};
			\node at (2.3,2.4) {\small $\ccc$};
			\node at (2.3,3.65) {\small $\aaa$};
			\node at (3.7,0.3) {\small $\aaa$};
			\node at (3.7,1.6) {\small $\ccc$};
			\node at (3.7,2.4) {\small $\ddd$};
			\node at (3.7,3.5) {\small $\bbb$};
			\node at (1,-3) { $E_{23}=c,E_{34}=b$};
		\end{tikzpicture}\hspace{10pt}
		\begin{tikzpicture}[>=latex,scale=0.45]
			\draw (0,0)--(2,0)
			(0,0)--(0,2)
			(0,0)--(-2,-2)
			(0,2)--(0,4)
			(0,4)--(-2,4)
			(0,4)--(2,4)
			(4,2)--(4,4)
			(4,2)--(4,0)
			(4,4)--(2,4)
			(4,0)--(2,0)
			
			(4,-2)--(4,0);
			\draw[line width=1.5]
			(0,2)--(-2,2)
			(2,0)--(2,2)
			(2,2)--(2,4)
			(2,-2)--(-2,-2)
			(2,-2)--(4,-2);
			\draw[dashed]
			(-2,2)--(-2,-2)
			(0,2)--(2,2)
			(-2,2)--(-2,4)
			(2,2)--(4,2)
			(2,0)--(2,-2);

			\node at (1.65,-0.4) {\small $\ccc$};
			\node at (2.3,-0.4) {\small $\ccc$};
			\node at (1.65,-1.55) {\small $\ddd$};
			\node at (2.3,-1.55) {\small $\ddd$};
			\node at (-1.05,-1.5) {\small $\bbb$};
			\node at (3.65,-1.55) {\small $\bbb$};
			\node at (3.65,-0.3) {\small $\aaa$};
			
			\node at (0.3,3.7) {\small $\aaa$};
			\node at (0.3,0.25) {\small $\aaa$};
			\node at (0.2,-0.3) {\small $\aaa$};
			\node at (0.3,2.3) {\small $\ccc$};
			\node at (0.3,1.6) {\small $\ccc$};
			\node at (-0.3,0.1) {\small $\aaa$};
			\node at (-1.7,1.55) {\small $\ddd$};
			\node at (-1.7,2.4) {\small $\ddd$};
			\node at (-1.7,-1.3) {\small $\ccc$};
			\node at (-1.7,3.6) {\small $\ccc$};
			\node at (1.7,0.4) {\small $\bbb$};
			\node at (1.7,1.6) {\small $\ddd$};
			\node at (1.7,2.4) {\small $\ddd$};
			\node at (1.7,3.5) {\small $\bbb$};
			\node at (2.35,0.4) {\small $\bbb$};
			\node at (2.35,1.55) {\small $\ddd$};
			\node at (2.35,2.45) {\small $\ddd$};
			\node at (2.35,3.55) {\small $\bbb$};
			\node at (3.7,0.3) {\small $\aaa$};
			\node at (3.7,1.6) {\small $\ccc$};
			\node at (3.7,2.3) {\small $\ccc$};
			\node at (3.7,3.7) {\small $\aaa$};
			\node at (-0.3,3.7) {\small $\aaa$};
			\node at (-0.3,1.5) {\small $\bbb$};
			\node at (-0.3,2.5) {\small $\bbb$};
			\node[draw,shape=circle, inner sep=0.5] at (1,1) {\small $1$};
			\node[draw,shape=circle, inner sep=0.5] at (-1,3) {\small $3$};
			\node[draw,shape=circle, inner sep=0.5] at (-1,1) {\small $2$};
			\node[draw,shape=circle, inner sep=0.5] at (1,-1) {\small $8$};
			\node[draw,shape=circle, inner sep=0.5] at (3,1) {\small $6$};
			\node[draw,shape=circle, inner sep=0.5] at (3,-1) {\small $7$};
			\node[draw,shape=circle, inner sep=0.5] at (3,3) {\small $5$};
			\node[draw,shape=circle, inner sep=0.5] at (1,3) {\small $4$};
			
			\node at (1,-3) {$E_{23}=b$};
		\end{tikzpicture}
		
		\caption{Partial neighborhoods of the 1st special tile.}
		\label{applec}
	\end{figure}
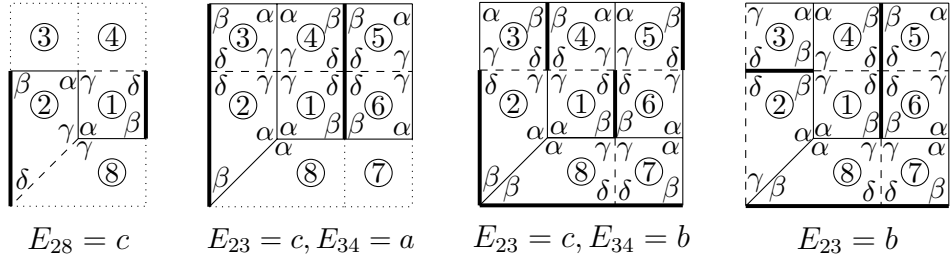
	Let the first of Fig. \ref{3444-5} be the center tile $T_1$ in the partial neighborhoods in Fig. \ref{applec}. If $E_{28}=c$ in the first picture, then $\aaa_1\cdots=\aaa_1\ccc_2\ccc_8$. This determines $T_2$. By Lemma \ref{lem4}, $\aaa_2\ccc_1\cdots=\aaa^2\ccc^2$, contradicting $\aaa_1\ccc_2\ccc_8$. If $E_{28}=b$, we get similar contradiction. So we have $E_{28}=a,\aaa_1\cdots=\aaa^3$ in all other pictures. There are three possibilities shown in the 2nd, 3rd and 4th pictures, according to $E_{23}=b$ or $c$ and $E_{34}=a$ or $b$.
	
	\subsubsection*{Case $E_{23}=c,E_{34}=a$}
	These two edges determine $T_2,T_3,T_4$. By Lemma \ref{lem4}, $\thick\ddd_1\dash\ddd_4\thick\cdots=\bbb^2\ddd^2$ or $\ddd^4$. By Proposition $\ref{proposition1}'$, we know $\bbb^2\ddd^2$ cannot be a vertex. So $\ddd_1\ddd_4\cdots=\ddd^4$, which determines $T_5,T_6$. By $\aaa_8$, we get $\bbb_1\bbb_6\cdots=\bbb^4$ or $\bbb^2\ccc^2$. By $\ccc^4$, we get $\bbb=\ccc$, contradicting Lemma \ref{lemd4}. 	
	
	\subsubsection*{Case $E_{23}=c,E_{34}=b$}
	These two edges determine $T_2,T_3,T_4$. By Lemma \ref{lem4},  $\ccc_4\ddd_1\cdots=\ccc^2\ddd^2$ determines $T_5,T_6$. By $\aaa_8$, $\bbb_1\bbb_6\cdots=\bbb^4$ or $\bbb^2\ccc^2$. By Proposition \ref{proposition1}, we know $\bbb^4$ cannot be vertex. So $\bbb_1\bbb_6\cdots=\bbb^2\ccc^2$, which determine $T_7,T_8$. So  $\aaa^3,\bbb^2\ccc^2,\ccc^2\ddd^2,\bbb_2\thin\bbb_8\cdots$ are vertices. By Proposition \ref{proposition2}, we get a contradiction.
	
	\subsubsection*{Case $E_{23}=b$}
	This edge determines $T_2$. By Lemma \ref{lem4}, $\bbb_2\ccc_1\cdots=\bbb^2\ccc^2$  determines $T_3,T_4$. 
	The proof so far shows we cannot have $\ccc_1\ccc_2\cdots$. This means we must have $\ccc_1\bbb_2\cdots$. By the symmetry, we also must have $\bbb_1\ccc_8\cdots$. Then $E_{78}=c$ determines $T_7,T_8$. And we can further determine $T_5,T_6$. 
	
	The angle sums at $\aaa^3,\bbb^2\ccc^2,\ddd^4$ imply $\aaa=\frac{2\pi}{3},\bbb+\ccc=\pi,\ddd=\frac{\pi}{2},f=24$. Since $\aaa=\frac{2\pi}{3},\bbb,\ccc,\ddd<\pi$, by Lemma \ref{geometry4}, we have $\bbb,\ccc>\frac{\pi}{6}$.
	By symmetry, we may assume $\bbb>\ccc$. Then the AVC is derived as shown in Table \ref{tab-2}.
	\begin{table}[htp]  
		\centering
		\caption{AVC for $\aaa = \frac{2\pi}{3},\ddd=\frac{\pi}{2},\bbb +\ccc=\pi$ with $\bbb>\ccc>\frac{\pi}{6}$.}\label{tab-2}
		\begin{minipage}{0.3\textwidth}
			\begin{tabular}{c|c}
				\hline
				$\text{vertex}$ &$\bbb$ \\
				\hline
				\hline
				$\aaa^3,\bbb^2\ccc^2,\ddd^4$&all\\
				$\ccc^{10}$&{\scriptsize $4\pi/5$}\\
				$\aaa\ccc^6$&{\scriptsize $7\pi/9$}\\
				\hline
			\end{tabular}
		\end{minipage}
		\raisebox{0.58em}{\begin{minipage}{0.3\textwidth}
				\begin{tabular}{c|c}					
					\hline
					$\text{vertex}$ &$\bbb$ \\
					\hline
					\hline					
					$\bbb\ccc^3\ddd,\ccc^4\ddd^2,\ccc^8$&{\scriptsize $3\pi/4$}\\
					$\aaa\bbb^2,\aaa\ccc^4,\aaa^2\ccc^2,\ccc^6$&{\scriptsize $2\pi/3$}\\
					\hline
				\end{tabular}
		\end{minipage}}
	\end{table}
	
	For $\bbb=\frac{3\pi}{4}$, we have 
	\begin{align*}
		\#\bbb&=2\#\bbb^2\ccc^2+\#\bbb\ccc^3\ddd,\\
		\#\ccc&=2\#\bbb^2\ccc^2+3\#\bbb\ccc^3\ddd+4\#\ccc^4\ddd^2+\#\ccc^{8}. 
	\end{align*}
	By $\#\bbb=f=\#\ccc$, the equalities above imply $\#\bbb\ccc^3\ddd=\#\ccc^4\ddd^2=\#\ccc^8=0$. This means $\bbb\ccc^3\ddd,\ccc^4\ddd^2,\ccc^8$ are not vertices.
	
	By the similar argument, we know $\ccc^{10}$ is not a vertex for $\bbb=\frac{4\pi}{5}$ and $\aaa\ccc^6$ is not a vertex for $\bbb=\frac{7\pi}{9}$. Moreover, for $\bbb=\frac{2\pi}{3}$, we know $\aaa\bbb^2$ is a vertex if and only if one of $\aaa\ccc^4,\aaa^2\ccc^2,\ccc^6$ is a vertex.	
	In the first of Fig. \ref{quadsub}, $\aaa_1\bbb_2\bbb_3=\thin^{\bbb}\aaa_1^{\ccc}\thin^{\aaa}\bbb_2^{\ddd}\thick^{\ddd}\bbb_3^{\aaa}\thin$ determines $T_1,T_2,T_3$. Then  $\aaa_3\bbb_1\cdots=\aaa\bbb^2$, $\ddd_2\ddd_3\cdots=\ddd^4$, $\ccc_3\cdots=\bbb^2\ccc^2$, $\aaa\ccc^4$, $\aaa^2\ccc^2$ or $\ccc^6$ make $T_3$ a $334d$-Tile, which has been handled in Section \ref{3345d}.
	\begin{figure}[htp]
		\centering
		\begin{tikzpicture}[>=latex,scale=1.5]
			\begin{scope}[scale=0.4]
				\draw
				(0,0)--(0,-2)
				(0,0)--(0,2)
				(0,2)--(2,2)
				(0,-2)--(2,-2);
				\draw[line width=1.5]
				(0,0)--(2,0)
				(0,-2)--(-2,0);
				\draw[dashed]
				(-2,0)--(0,2)
				(2,2)--(2,-2);
				\node[draw,shape=circle, inner sep=0.5] at (-1,0) {\small $1$};
				\node[draw,shape=circle, inner sep=0.5] at (1,1) {\small $2$};
				\node[draw,shape=circle, inner sep=0.5] at (1,-1) {\small $3$};
				\node at (0.35,0.35) {\small $\bbb$};
				\node at (0.35,1.75) {\small $\aaa$};
				\node at (-0.07,-2.45) {\small $\bbb$};
				\node at (-0.35,0) {\small $\aaa$};
				\node at (0.35,-1.75) {\small $\aaa$};
				\node at (0.35,-0.45) {\small $\bbb$};
				\node at (-1.55,0.07) {\small $\ddd$};
				\node at (1.75,0.35) {\small $\ddd$};
				\node at (1.75,-0.35) {\small $\ddd$};
				\node at (2.2,0.35) {\small $\ddd$};
				\node at (2.2,-0.35) {\small $\ddd$};
				\node at (1.75,1.7) {\small $\ccc$};
				\node at (1.75,-1.7) {\small $\ccc$};
				\node at (-0.25,1.3) {\small $\ccc$};
				\node at (-0.25,-1.25) {\small $\bbb$};	
			\end{scope}	
			\end{tikzpicture}\hspace{30pt}			 		    
			\includegraphics[scale=0.18]{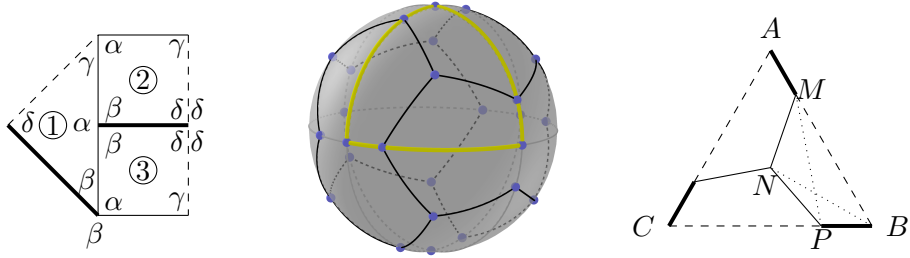} \hspace{15pt}	
	     	\begin{tikzpicture}	[>=latex,scale=1.5]
			\begin{scope}[scale=0.45]
				\draw[dashed]
				(-2,0)--(1,0)
				(-1.5,0.87)--(0,3.46)
				(0.5,2.59)--(2,0);
				\draw[line width=1.5]
				(1,0)--(2,0)
				(-2,0)--(-1.5,0.87)
				(0,3.46)--(0.5,2.59);
				\draw
				(0,1.15)--(0.5,2.59)
				(0,1.15)--(1,0)
				(0,1.15)--(-1.5,0.87);
				\draw[dotted]
				(2,0)--(0,1.15)
				(1,0)--(0.5,2.59);
				
				\node at (0,3.9) {\small $A$};
				\node at (2.5,0) {\small $B$};
				\node at (-2.5,0) {\small $C$};
				\node at (-0.1,0.8) {\small $N$};
				\node at (1,-0.25) {\small $P$};
				\node at (0.8,2.7) {\small $M$};
			\end{scope}
		\end{tikzpicture}  	
		\caption{Pictures for $\aaa\bbb^2$ and quadrilateral subdivisions of the octahedron.}
		\label{quadsub}
	\end{figure}

	In conclusion, we only need to consider $\text{AVC}=\{\aaa^3,\bbb^2\ccc^2,\ddd^4\}$,  which induces a quadrilateral subdivision of the octahedron in Fig. \ref{case1jie1} by exactly the same deduction as ``\textbf{Subcase. For $f=24$, $\aaa^3$ is a vertex}'' on page \pageref{subcase3a} in Section \ref{3345d}.
   \end{proof}
	
	The following computation shows that the quadrilateral in this subdivision tiling could be more general than the one in Section \ref{3345d}.

	\subsubsection*{Moduli of quadrilateral subdivisions of the octahedron}
	
	The second and third pictures of Fig. \ref{quadsub} shows a quadrilateral subdivision of one triangular face of the regular octahedron. We have 
	\[\aaa = \tfrac{2\pi}{3},\ddd=\tfrac{\pi}{2},\bbb+\ccc=\pi, PN= MN=a, BP=b, MB=c, b+c=\tfrac{\pi}{2}.\]
	By symmetry between $\bbb,\ccc$, we may assume that $b\in (0,\frac{\pi}{4}]$. Then 
	\[\cos{MP}=\cos^2 a+\sin^2 a\cos\aaa=\cos b\cos c.\] 
	Solving the above equation, we get   
	\[ a = \arccos{\sqrt{\tfrac{\sin{2b}+1}{3}}}.\]
	Similarly, the cosine laws for $BN$ give	
	\[\cos{BN}=\cos a\cos c+\sin a\sin c\cos \ccc=\cos a\cos b+\sin a\sin b\cos {(\pi-\ccc)},\]
	which implies 
	\[ \ccc =\arccos{\frac{\cos b-\sin b}{(\sin b +\cos b)\tan {(\arccos{\sqrt{\frac{\sin{2b}+1}{3}}})}}}.\]
	
	When $\bbb=\frac{2\pi}{3},\ccc=\frac{\pi}{3}$, we get $b=\arctan(\frac{3-\sqrt{5}}{2})\approx0.116\pi$. This is the special quadrilateral admitting a flip modification in the 3rd picture of Fig. \ref{real figure} in Section \ref{introduction} and Fig. \ref{case1jie2} in Section \ref{3345d} after exchanging $\bbb\leftrightarrow\ccc$.        
	
	Solving $a=b$, we get $b=\arctan(\sqrt 3 -1)\approx0.201\pi$. Solving $a=c$, we get $b=\arctan{\frac{1+\sqrt 3}{2}}\approx0.298\pi > \frac{\pi}{4} $. When $b=c=\frac{\pi}{4}$, we get $\bbb=\ccc=\frac{\pi}{2}$ and it is of Type $a^2b^2$. Therefore the moduli of these quadrilateral subdivisions is given by $(0,\frac{\pi}{4}]$ for the length of $b$, and the quadrilateral reduces from Type $a^2bc$ to Type $a^3b$ when $b=\arctan(\sqrt 3 -1)$, to Type $a^2b^2$ when $b=\frac{\pi}{4}$. 
	
	\vspace{9pt}
	We conclude that all $a^2bc$-tilings with the $1st$ special tile in Fig. \ref{3444-5} are these quadrilateral subdivisions for $b\in(0,\frac{\pi}{4})\backslash\{\arctan(\sqrt 3 -1)\}$.

\begin{proposition}
	There is no $a^2bc$-tiling with the 2nd special tile in Fig. \ref{3444-5}.
\end{proposition}
\begin{proof}
	\begin{figure}[htp]
		\centering
		\begin{tikzpicture}[>=latex,scale=0.45]
			\draw
			(0,0)--(0,4)
			(-2,0)--(2,0)
			(-2,-2)--(-2,0)
			(-2,4)--(4,4)
			(2,0)--(4,-2)
			(4,-2)--(4,4);
			\draw[line width=1.5]
			(0,2)--(-2,2)
			(2,0)--(2,4)
			(4,-2)--(-2,-2);
			\draw[dashed]
			(-2,2)--(-2,4)
			(0,2)--(4,2)
			(0,0)--(0,-2)
			(-2,0)--(-2,2);
			\node at (0.35,0.3) {\small $\aaa$};
			\node at (0.35,3.65) {\small $\aaa$};
			\node at (0.35,2.35) {\small $\ccc$};
			\node at (0.35,-0.4) {\small $\ccc$};
			\node at (0.3,-1.55) {\small $\ddd$};
			\node at (-0.35,0.3) {\small $\aaa$};
			\node at (-0.3,2.4) {\small $\bbb$};
			\node at (-0.35,3.65) {\small $\aaa$};
			\node at (-0.35,-0.4) {\small $\ccc$};
			\node at (-0.3,-1.55) {\small $\ddd$};
			\node at (0.35,1.65) {\small $\ccc$};
			\node at (1.65,0.45) {\small $\bbb$};
			\node at (1.7,1.6) {\small $\ddd$};
			\node at (-0.3,1.5) {\small $\bbb$};
			\node at (-1.65,0.35) {\small $\ccc$};
			\node at (-1.65,-0.3) {\small $\aaa$};
			\node at (-1.7,-1.55) {\small $\bbb$};
			\node at (-1.7,1.55) {\small $\ddd$};
			\node at (-1.7,2.4) {\small $\ddd$};
			\node at (-1.7,3.6) {\small $\ccc$};
			\node at (1.8,-0.3) {\small $\aaa$};
			\node at (2.35,0.4) {\small $\bbb$};
			\node at (2.9,-1.55) {\small $\bbb$};
			\node at (3.65,-1.2) {\small $\aaa$};
			\node at (3.65,1.65) {\small $\ccc$};
			\node at (3.65,2.35) {\small $\ccc$};
			\node at (3.65,3.65) {\small $\aaa$};
			\node at (2.35,1.6) {\small $\ddd$};
			\node at (2.35,2.4) {\small $\ddd$};
			\node at (1.7,2.4) {\small $\ddd$};
			\node at (2.35,3.5) {\small $\bbb$};
			\node at (1.7,3.5) {\small $\bbb$};
			\node[draw,shape=circle, inner sep=0.5] at (1,1) {\small $1$};
			\node[draw,shape=circle, inner sep=0.5] at (3,3) {\small $3$};
			\node[draw,shape=circle, inner sep=0.5] at (3,1) {\small $4$};
			\node[draw,shape=circle, inner sep=0.5] at (1,3) {\small $2$};
			\node[draw,shape=circle, inner sep=0.5] at (-1,1) {\small $7$};
			\node[draw,shape=circle, inner sep=0.5] at (-1,3) {\small $8$};
			\node[draw,shape=circle, inner sep=0.5] at (-1,-1) {\small $6$};
			\node[draw,shape=circle, inner sep=0.5] at (1,-1) {\small $5$};
			\node at (1,-3.5) {\small $E_{78}=b$};
			\node at (0.3,-2.5) {\small $\ddd$};
			\node at (-0.3,-2.5) {\small $\ddd$};
			\node at (4.25,-2.25) {\small $\bbb$};
		\end{tikzpicture}\hspace{30pt}
		\begin{tikzpicture}[>=latex,scale=0.45]
			\draw
			(0,0)--(0,4)
			(-2,0)--(2,0)
			(-2,-2)--(-2,0)
			(-2,4)--(4,4)
			(2,0)--(4,-2)
			(4,-2)--(4,4);
			\draw[line width=1.5]
			(-2,2)--(-2,4)
			(2,0)--(2,4)
			(4,-2)--(-2,-2)
			(-2,0)--(-2,2);
			\draw[dashed]
			
			(0,2)--(4,2)
			(0,0)--(0,-2)
			(0,2)--(-2,2);
			\node at (0.35,0.3) {\small $\aaa$};
			\node at (0.35,3.65) {\small $\aaa$};
			\node at (0.35,2.35) {\small $\ccc$};
			\node at (0.35,-0.4) {\small $\ccc$};
			\node at (0.3,-1.55) {\small $\ddd$};
			\node at (-0.35,0.3) {\small $\aaa$};
			\node at (-0.3,2.35) {\small $\ccc$};
			\node at (-0.35,3.65) {\small $\aaa$};
			\node at (-0.35,-0.4) {\small $\ccc$};
			\node at (-0.3,-1.55) {\small $\ddd$};
			\node at (0.35,1.65) {\small $\ccc$};
			\node at (1.65,0.45) {\small $\bbb$};
			\node at (1.7,1.6) {\small $\ddd$};
			\node at (-0.3,1.65) {\small $\ccc$};
			\node at (-1.65,0.45) {\small $\bbb$};
			\node at (-1.65,-0.3) {\small $\aaa$};
			\node at (-1.7,-1.55) {\small $\bbb$};
			\node at (-1.7,1.6) {\small $\ddd$};
			\node at (-1.7,2.4) {\small $\ddd$};
			\node at (-1.65,3.5) {\small $\bbb$};
			\node at (1.8,-0.3) {\small $\aaa$};
			\node at (2.35,0.4) {\small $\bbb$};
			\node at (2.9,-1.55) {\small $\bbb$};
			\node at (3.65,-1.2) {\small $\aaa$};
			\node at (3.65,1.65) {\small $\ccc$};
			\node at (3.65,2.35) {\small $\ccc$};
			\node at (3.65,3.65) {\small $\aaa$};
			\node at (2.35,1.6) {\small $\ddd$};
			\node at (2.35,2.4) {\small $\ddd$};
			\node at (1.7,2.4) {\small $\ddd$};
			\node at (2.35,3.5) {\small $\bbb$};
			\node at (1.7,3.5) {\small $\bbb$};
			\node[draw,shape=circle, inner sep=0.5] at (1,1) {\small $1$};
			\node[draw,shape=circle, inner sep=0.5] at (3,3) {\small $3$};
			\node[draw,shape=circle, inner sep=0.5] at (3,1) {\small $4$};
			\node[draw,shape=circle, inner sep=0.5] at (1,3) {\small $2$};
			\node[draw,shape=circle, inner sep=0.5] at (-1,1) {\small $7$};
			\node[draw,shape=circle, inner sep=0.5] at (-1,3) {\small $8$};
			\node[draw,shape=circle, inner sep=0.5] at (-1,-1) {\small $6$};
			\node[draw,shape=circle, inner sep=0.5] at (1,-1) {\small $5$};
			\node at (1,-3.5) {\small $E_{78}=c,E_{28}=a$};
			
			\node at (0.3,-2.5) {\small $\ddd$};
			\node at (-0.3,-2.5) {\small $\ddd$};
			\node at (4.25,-2.25) {\small $\bbb$};
			
		\end{tikzpicture}\hspace{30pt}
		\begin{tikzpicture}[>=latex,scale=0.45]
			\draw
			(0,0)--(0,2)
			(-2,0)--(2,0)
			(-2,-2)--(-2,0)
			(-2,4)--(6,4)
			(2,0)--(4,0)
			(4,0)--(6,2)
			(-2,2)--(-2,4)
			(2,2)--(2,4)
			(6,2)--(6,-2);
			\draw[line width=1.5]
			(6,2)--(6,4)
			(-2,0)--(-2,2)
			(2,0)--(2,2)
			(4,-2)--(-2,-2)
			(0,2)--(0,4)
			(4,-2)--(4,0);
			\draw[dashed]
			
			(0,2)--(4,2)
			(0,0)--(4,-2)
			(0,2)--(-2,2)
			(4,2)--(6,2)
			(6,-2)--(4,-2);
			\node at (0.35,0.3) {\small $\aaa$};
			\node at (0.35,3.5) {\small $\bbb$};
			\node at (0.35,2.35) {\small $\ddd$};
			\node at (1.5,-0.4) {\small $\ccc$};
			\node at (3.7,-1.4) {\small $\ddd$};
			\node at (-0.35,0.3) {\small $\aaa$};
			\node at (-0.3,2.4) {\small $\ddd$};
			\node at (-0.35,3.5) {\small $\bbb$};
			\node at (-0.1,-0.4) {\small $\ccc$};
			\node at (2.1,-1.55) {\small $\ddd$};
			\node at (0.35,1.65) {\small $\ccc$};
			\node at (1.65,0.45) {\small $\bbb$};
			\node at (1.7,1.6) {\small $\ddd$};
			\node at (-0.3,1.65) {\small $\ccc$};
			\node at (-1.65,0.45) {\small $\bbb$};
			\node at (-1.65,-0.3) {\small $\aaa$};
			\node at (-1.7,-1.55) {\small $\bbb$};
			\node at (-1.7,1.6) {\small $\ddd$};
			\node at (-1.7,2.4) {\small $\ccc$};
			\node at (-1.7,3.65) {\small $\aaa$};
			\node at (2.2,-0.3) {\small $\aaa$};
			\node at (2.35,0.45) {\small $\bbb$};
			\node at (4.35,-0.15) {\small $\bbb$};
			\node at (3.7,-0.45) {\small $\bbb$};
			\node at (3.8,0.25) {\small $\aaa$};
			\node at (5.25,1.65) {\small $\ccc$};
			\node at (5.65,2.4) {\small $\ddd$};
			\node at (5.65,3.5) {\small $\bbb$};
			\node at (2.35,1.6) {\small $\ddd$};
			\node at (2.35,2.4) {\small $\ccc$};
			\node at (1.7,2.4) {\small $\ccc$};
			\node at (2.35,3.65) {\small $\aaa$};
			\node at (1.7,3.65) {\small $\aaa$};
			\node at (4.3,-1.6) {\small $\ddd$};
			\node at (5.65,1.2) {\small $\aaa$};	
			\node at (5.65,-1.65) {\small $\ccc$};
			\node[draw,shape=circle, inner sep=0.5] at (5,-1) {\small $9$};
			\node[draw,shape=circle, inner sep=0.5] at (1,1) {\small $1$};
			\node[draw,shape=circle, inner sep=0.5] at (4,3) {\small $3$};
			\node[draw,shape=circle, inner sep=0.5] at (3.5,1) {\small $4$};
			\node[draw,shape=circle, inner sep=0.5] at (1,3) {\small $2$};
			\node[draw,shape=circle, inner sep=0.5] at (-1,1) {\small $7$};
			\node[draw,shape=circle, inner sep=0.5] at (-1,3) {\small $8$};
			\node[draw,shape=circle, inner sep=0.5] at (-1,-1) {\small $6$};
			\node[draw,shape=circle, inner sep=0.5] at (3,-0.7) {\small $5$};
			\node at (2,-3.5) {\small $E_{78}=c,E_{28}=b$};
		\end{tikzpicture}
		
		\caption{Partial neighborhoods of the 2nd special tile.}\label{case3.1.2}
	\end{figure}
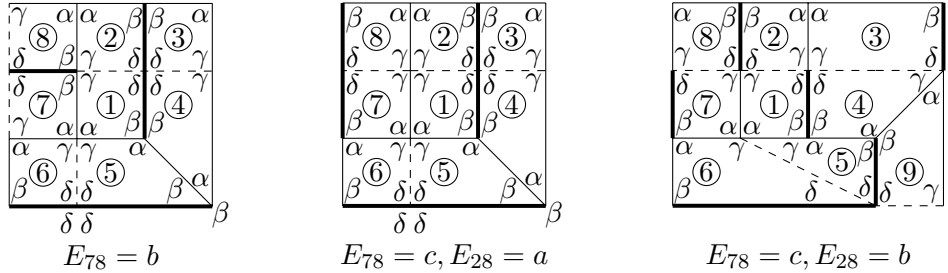
	Let the second of Fig. \ref{3444-5} be the center tile $T_1$ in the partial neighborhoods in Fig. \ref{case3.1.2}. By Lemma \ref{lem3}, $\bbb_1\cdots=\aaa_5\bbb_1\bbb_4$ determines $T_4$. By $\aaa_5\bbb_1\bbb_4$ and Lemma \ref{lem4}, $\aaa_1\cdots=\aaa^2\ccc^2=\aaa_1\aaa_7\ccc_5\ccc_6$, which determines $T_5,T_6$.	
	There are three possibilities shown in Fig. \ref{case3.1.2}, according to $E_{78}=b$ or $c$ and $E_{28}=a$ or $b$.
	
	\subsubsection*{Case $E_{78}=b$}
	This edge determines $T_7$. By Lemma \ref{lem4}, $\bbb_7\ccc_1\cdots=\bbb^2\ccc^2$ determines $T_2,T_8$, which further determine  $\ddd_1\ddd_2\ddd_3\ddd_4$ and $T_3$.	
	The angle sums at $\aaa\bbb^2$, $\aaa^2\ccc^2$, $\bbb^2\ccc^2$, $\ddd^4$ imply $\aaa=\bbb=\frac{2\pi}{3},\ccc=\frac{\pi}{3},\ddd=\frac{\pi}{2}$. By the edge length consideration, we get $\aaa_4\bbb_5\cdots=\aaa\bbb^2,\ddd_5\ddd_6\cdots=\ddd^4$.  So $T_5$ is a $3344$-Tile, which has been handled in Section \ref{3345d}.
	
	\subsubsection*{Case $E_{78}=c,E_{28}=a$}
	
	These two edges determine $T_2,T_7,T_8$, which further determine  $\ddd_1\ddd_2\ddd_3\ddd_4$ and $T_3$. The angle sums at $\aaa\bbb^2$, $\aaa^2\ccc^2$, $\ccc^4$,  $\ddd^4$ imply $\aaa=\frac{\pi}{2},\bbb=\frac{3\pi}{4},\ccc=\frac{\pi}{2},\ddd=\frac{\pi}{2}$.  By  the edge length consideration, we get $\aaa_4\bbb_5\cdots=\aaa\beta^2,\thick\ddd_5\dash\ddd_6\thick\cdots=\ddd^4$. So $T_5$ is a $3344$-Tile, which has been handled in Section \ref{3345d}.

	\subsubsection*{Case $E_{78}=c,E_{28}=b$}
	
	These two edges determine $T_2,T_7,T_8$, which further determine $\ccc_2\ccc_3\ddd_1\ddd_4$ and $T_3$. 
	By $\aaa\bbb^2$ and Parity Lemma, $\aaa_4\bbb_5\cdots=\aaa_4\bbb_5\bbb_9$ determines $T_9$. Then $\aaa_9\ccc_4\ddd_3\cdots$ contradicts Lemma $\ref{lemac2}'$.
\end{proof}

\begin{proposition}
	There is no $a^2bc$-tiling with the 3rd special tile in Fig. \ref{3444-5}.
\end{proposition}

\begin{proof}
	Let the third of Fig. \ref{3444-5} be the center tile $T_1$ in the partial neighborhoods in Fig. \ref{ghj}. 	
	If $E_{67}=c$ in the first picture, then $\aaa_1\cdots=\aaa_1\ccc_6\ccc_7$. This determines $T_7$. By Parity Lemma, $\aaa_7\ccc_1\cdots$ muse be even. This implies $\aaa+2\ccc<2\pi$, contradicting $\aaa_1\ccc_6\ccc_7$. If $E_{67}=b$, we get similar contradiction. So we have $E_{67}=a,\aaa_1\cdots=\aaa^3$ in all other pictures.	
	
	If $E_{56}=b$, then $\bbb_1\bbb_6\cdots=\bbb^4$ or $\bbb^2\ddd^2$ in the 2nd and 3rd picture. 
	
	If $E_{56}=c$, then Lemma \ref{lem4} implies $\bbb_1\ccc_6\cdots=\bbb^2\ccc^2$, which determines $T_4,T_5$.	
	We have $E_{23}=a$ or $b$. If $E_{23}=b$, we get $\ddd_1\ddd_2\ddd_3\ddd_4$, which determines $T_2,T_3$. Then we have $E_{78}=b$ in the fourth picture or $E_{78}=c$ in the fifth picture.  If $E_{23}=a$, we have the sixth picture.	
	
	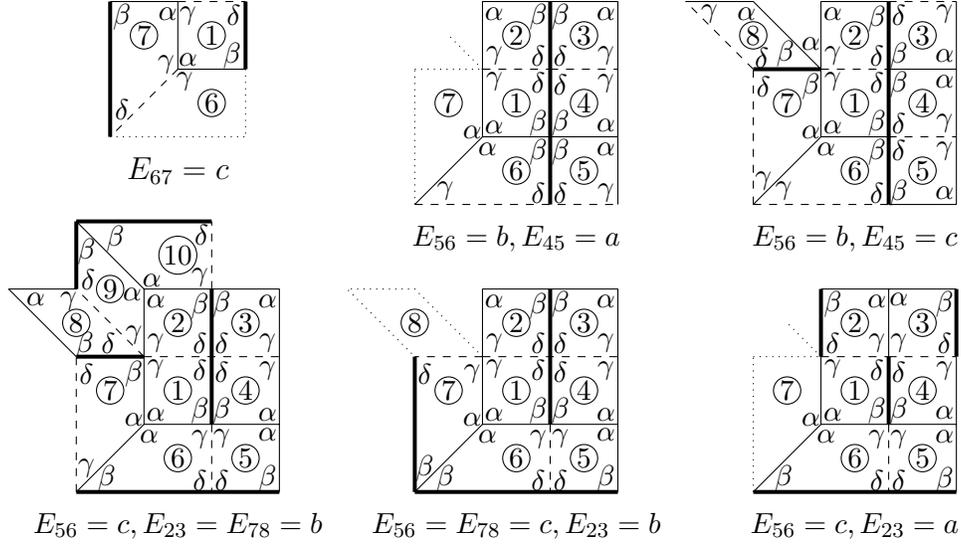
\begin{figure}[htp]
		\centering
		\begin{tikzpicture}[>=latex,scale=0.45]
			
			\begin{scope}[xshift=-9cm,yshift=2cm]						
			\draw (0,0)--(2,0)
			(0,0)--(0,2)			
			(0,2)--(-2,2);
			\draw[line width=1.5]
			(2,0)--(2,2)		
			(-2,-2)--(-2,2);
			\draw[dashed]
			(0,2)--(2,2)
			(0,0)--(-2,-2);
			\draw[dotted]
			(-2,-2)--(2,-2)
			(2,0)--(2,-2);
			\node at (0.3,0.3) {\small $\aaa$};
			\node at (0.3,1.6) {\small $\ccc$};
			\node at (-0.3,1.65) {\small $\aaa$};
			\node at (1.6,0.4) {\small $\bbb$};
			\node at (1.65,1.6) {\small $\ddd$};
			\node at (-0.35,0.15) {\small $\ccc$};
			\node at (-1.65,-1.2) {\small $\ddd$};
			\node at (0.2,-0.35) {\small $\ccc$};
			\node at (-1.65,1.55) {\small $\bbb$};
			\node[draw,shape=circle, inner sep=0.5] at (1,1) {\small $1$};
			\node[draw,shape=circle, inner sep=0.5] at (-1,1) {\small $7$};
			\node[draw,shape=circle,inner sep=0.5] at (1,-1){\small $6$};
			\node at (0,-3){$E_{67}= c$};
		\end{scope}
	
		\begin{scope}[xshift=0cm,yshift=0cm]   
			\draw (0,0) -- (0,2) 
			(0,0)--(-2,-2)
			(0,0)--(2,0)						
			(0,4)--(4,4)		
			(4,-2)--(4,2)
			(2,0)--(4,0)
			(0,2)--(0,4)
			(4,2)--(4,4);
			\draw[line width=1.5]
			(2,0)--(2,2)
			(2,0)--(2,-2)				
			(2,2)--(2,4);
			\draw[dashed]
			(0,2)--(2,2)	
			(2,2)--(4,2)			
			(-2,-2)--(4,-2);
			\draw[dotted]
			(-2,2)--(0,2)
			(-2,-2)--(-2,2)
			(0,2)--(-1,3);
			\node at (2.3,0.4) {\small $\bbb$};
			\node at (2.3,1.55) {\small $\ddd$};
			\node at (2.3,2.4) {\small $\ddd$};
			\node at (2.3,3.55) {\small $\bbb$};
			\node at (2.3,-0.45) {\small $\bbb$};
			\node at (2.3,-1.65) {\small $\ddd$};
			\node at (3.65,0.35) {\small $\aaa$};
			\node at (3.65,1.6) {\small $\ccc$};
			\node at (3.65,2.4) {\small $\ccc$};
			\node at (3.65,3.65) {\small $\aaa$};
			\node at (3.65,-0.35) {\small $\aaa$};
			\node at (3.65,-1.65) {\small $\ccc$};
			\node at (0.35,0.3) {\small $\aaa$};
			\node at (0.35,1.65) {\small $\ccc$};
			\node at (0.35,2.4) {\small $\ccc$};
			\node at (0.35,3.65) {\small $\aaa$};
			\node at (0.15,-0.35) {\small $\aaa$};
			\node at (-0.3,0.15) {\small $\aaa$};
			\node at (1.65,1.6) {\small $\ddd$};
			\node at (1.65,2.4) {\small $\ddd$};
			\node at (1.65,3.55) {\small $\bbb$};
			\node at (1.65,-0.45) {\small $\bbb$};
			\node at (1.65,0.4) {\small $\bbb$};
			\node at (1.65,-1.65) {\small $\ddd$};
			\node at (-1.1,-1.6) {\small $\ccc$};
			\node[draw,shape=circle, inner sep=0.5] at (1,1) {\small $1$};
			\node[draw,shape=circle, inner sep=0.5] at (-1,1) {\small $7$};
			\node[draw,shape=circle, inner sep=0.5] at (1,-1) {\small $6$};
			\node[draw,shape=circle, inner sep=0.5] at (3,-1) {\small $5$};
			\node[draw,shape=circle, inner sep=0.5] at (3,1) {\small $4$};
			\node[draw,shape=circle, inner sep=0.5] at (3,3) {\small $3$};
			\node[draw,shape=circle, inner sep=0.5] at (1,3) {\small $2$};
			\node at (1,-3) {\small $E_{56}=b,E_{45}=a$};
		\end{scope}
	
		\begin{scope}[xshift=10cm,yshift=0cm]       
			\draw (0,0) -- (0,2) 
			(0,0)--(-2,-2)
			(0,0)--(2,0)
			(0,2)--(0,4)
			(0,2)--(-2,4)
			(-2,4)--(-4,4)
			(0,4)--(2,4)
			(2,2)--(4,2)
			(4,-2)--(4,4)
			(2,-2)--(4,-2);
			\draw[line width=1.5]
			(2,0)--(2,2)
			(2,0)--(2,-2)
			(2,2)--(2,4)
			(-2,2)--(0,2);
			\draw[dashed]
			(0,2)--(2,2)
			(-2,-2)--(2,-2)
			(2,0)--(4,0)
			(-2,-2)--(-2,2)
			(-2,2)--(-4,4)
			(2,4)--(4,4);
			\node[draw,shape=circle, inner sep=0.5] at (1,1) {\small $1$};
			\node[draw,shape=circle, inner sep=0.5] at (-1,1) {\small $7$};
			\node[draw,shape=circle, inner sep=0.5] at (1,-1) {\small $6$};
			\node[draw,shape=circle, inner sep=0.5] at (3,-1) {\small $5$};
			\node[draw,shape=circle, inner sep=0.5] at (3,1) {\small $4$};
			\node[draw,shape=circle, inner sep=0.5] at (3,3) {\small $3$};
			\node[draw,shape=circle, inner sep=0.5] at (1,3) {\small $2$};
			\node[draw,shape=circle, inner sep=0.5] at (-2,3) {\small $8$};
			\node at (2.3,0.4) {\small $\ddd$};
			\node at (2.3,1.55) {\small $\bbb$};
			\node at (2.3,2.4) {\small $\bbb$};
			\node at (2.3,3.6) {\small $\ddd$};
			\node at (2.3,-0.4) {\small $\ddd$};
			\node at (2.3,-1.65) {\small $\bbb$};
			\node at (3.65,0.4) {\small $\ccc$};
			\node at (3.65,1.65) {\small $\aaa$};
			\node at (3.65,2.4) {\small $\aaa$};
			\node at (3.65,3.65) {\small $\ccc$};
			\node at (3.65,-0.4) {\small $\ccc$};
			\node at (3.65,-1.65) {\small $\aaa$};
			\node at (0.35,0.3) {\small $\aaa$};
			\node at (0.35,1.65) {\small $\ccc$};
			\node at (0.35,2.4) {\small $\ccc$};
			\node at (0.35,3.65) {\small $\aaa$};
			\node at (0.15,-0.35) {\small $\aaa$};
			\node at (-0.3,0.15) {\small $\aaa$};
			\node at (-0.3,1.5) {\small $\bbb$};
			\node at (1.65,1.6) {\small $\ddd$};
			\node at (1.65,2.4) {\small $\ddd$};
			\node at (1.65,3.5) {\small $\bbb$};
			\node at (1.65,-0.45) {\small $\bbb$};
			\node at (1.65,0.4) {\small $\bbb$};
			\node at (1.65,-1.65) {\small $\ddd$};
			\node at (-1.1,-1.6) {\small $\ccc$};
			\node at (-1.7,-1.3) {\small $\ccc$};
			\node at (-1.7,1.5) {\small $\ddd$};
			\node at (-1.75,2.35) {\small $\ddd$};
			\node at (-1.05,2.45) {\small $\bbb$};
			\node at (-2.2,3.65) {\small $\aaa$};
			\node at (-3.3,3.65) {\small $\ccc$};
			\node at (-0.3,2.75) {\small $\aaa$};
			\node at (1,-3) {\small $E_{56}=b,E_{45}=c$};
		\end{scope}	
	
	 \begin{scope}[xshift=-10cm,yshift=-8.5cm]   
		\draw (0,0) -- (0,2) 
		(0,0)--(-2,-2)
		(0,0)--(2,0)
		(0,2)--(0,4)
		(0,4)--(4,4)
		(2,0)--(4,0)
		(4,0)--(4,4)
		(4,0)--(4,-2)
		(0,4)--(-2,6)
		(-2,2)--(-4,4)
		(-4,4)--(-2,4);
		\draw[line width=1.5]
		(2,0)--(2,4)
		(-2,-2)--(4,-2)
		(0,2)--(-2,2)
		(-2,4)--(-2,6)
		(2,6)--(-2,6);
		\draw[dashed]
		(0,2)--(4,2)
		(2,0)--(2,-2)
		(-2,-2)--(-2,2)
		(0,2)--(-2,4)
		(2,4)--(2,6);
		\node at (2.3,0.4) {\small $\bbb$};
		\node at (2.3,1.55) {\small $\ddd$};
		\node at (2.3,2.4) {\small $\ddd$};
		\node at (2.3,3.55) {\small $\bbb$};
		\node at (2.3,-0.4) {\small $\ccc$};
		\node at (2.3,-1.65) {\small $\ddd$};
		\node at (3.65,0.3) {\small $\aaa$};
		\node at (3.65,1.6) {\small $\ccc$};
		\node at (3.65,2.4) {\small $\ccc$};
		\node at (3.65,3.65) {\small $\aaa$};
		\node at (3.65,-0.3) {\small $\aaa$};
		\node at (3.65,-1.6) {\small $\bbb$};
		\node at (0.35,0.3) {\small $\aaa$};
		\node at (0.35,1.65) {\small $\ccc$};
		\node at (0.35,2.4) {\small $\ccc$};
		\node at (0.35,3.65) {\small $\aaa$};
		\node at (0.15,-0.35) {\small $\aaa$};
		\node at (-0.3,0.15) {\small $\aaa$};
		\node at (-0.3,1.5) {\small $\bbb$};
		\node at (1.65,1.6) {\small $\ddd$};
		\node at (1.65,2.4) {\small $\ddd$};
		\node at (1.65,3.5) {\small $\bbb$};
		\node at (1.65,-0.45) {\small $\ccc$};
		\node at (1.65,0.4) {\small $\bbb$};
		\node at (1.65,-1.65) {\small $\ddd$};
		\node at (-1.1,-1.6) {\small $\bbb$};
		\node at (-1.7,-1.3) {\small $\ccc$};
		\node at (-1.7,1.5) {\small $\ddd$};
		\node at (-1.75,2.35) {\small $\bbb$};
		\node at (-1.05,2.4) {\small $\ddd$};
		\node at (-2.23,3.65) {\small $\ccc$};
		\node at (-3.2,3.7) {\small $\aaa$};
		\node at (-1.65,4.3) {\small $\ddd$};
		\node at (-1.7,5.1) {\small $\bbb$};
		\node at (-0.3,2.65) {\small $\ccc$};
		\node at (-0.35,3.85) {\small $\aaa$}; 
		\node at (0.23,4.25) {\small $\aaa$};
		\node at (-0.85,5.5) {\small $\bbb$};
		\node at (1.7,5.6) {\small $\ddd$};
		\node at (1.65,4.35) {\small $\ccc$};
		\node[draw,shape=circle, inner sep=0.5] at (1,1) {\small $1$};
		\node[draw,shape=circle, inner sep=0.5] at (-1,1) {\small $7$};
		\node[draw,shape=circle, inner sep=0.5] at (1,-1) {\small $6$};
		\node[draw,shape=circle, inner sep=0.5] at (3,-1) {\small $5$};
		\node[draw,shape=circle, inner sep=0.5] at (3,1) {\small $4$};
		\node[draw,shape=circle, inner sep=0.5] at (3,3) {\small $3$};
		\node[draw,shape=circle, inner sep=0.5] at (1,3) {\small $2$};
		\node[draw,shape=circle, inner sep=0.5] at (-1,4) {\small $9$};
		\node[draw,shape=circle, inner sep=0.5] at (-2,3) {\small $8$};
		\node[draw,shape=circle, inner sep=0.5] at (1,5) {\small $10$};
		\node at (1,-3) {\small $E_{56}=c,E_{23}=E_{78}=b$};
	  \end{scope}  
  
		\begin{scope}[xshift=0cm,yshift=-8.5cm]  
			\draw (0,0) -- (0,2) 
			(0,0)--(-2,-2)
			(0,0)--(2,0)
			(0,2)--(0,4)
			(0,4)--(4,4)
			(2,0)--(4,0)
			(4,0)--(4,4)
			(4,0)--(4,-2);
			\draw[line width=1.5]
			(2,0)--(2,4)
			(-2,-2)--(4,-2)
			(-2,-2)--(-2,2);
			\draw[dashed]
			(0,2)--(4,2)
			(2,0)--(2,-2)
			(0,2)--(-2,2);
			\draw[dotted]
			(0,2)--(-2,4)--(-4,4)--(-2,2);
			\node at (2.3,0.4) {\small $\bbb$};
			\node at (2.3,1.55) {\small $\ddd$};
			\node at (2.3,2.4) {\small $\ddd$};
			\node at (2.3,3.55) {\small $\bbb$};
			\node at (2.3,-0.4) {\small $\ccc$};
			\node at (2.3,-1.65) {\small $\ddd$};
			\node at (3.65,0.3) {\small $\aaa$};
			\node at (3.65,1.6) {\small $\ccc$};
			\node at (3.65,2.4) {\small $\ccc$};
			\node at (3.65,3.65) {\small $\aaa$};
			\node at (3.65,-0.3) {\small $\aaa$};
			\node at (3.65,-1.6) {\small $\bbb$};
			\node at (0.35,0.3) {\small $\aaa$};
			\node at (0.35,1.65) {\small $\ccc$};
			\node at (0.35,2.4) {\small $\ccc$};
			\node at (0.35,3.65) {\small $\aaa$};
			\node at (0.15,-0.35) {\small $\aaa$};
			\node at (-0.3,0.15) {\small $\aaa$};
			\node at (-0.3,1.5) {\small $\ccc$};
			\node at (1.65,1.6) {\small $\ddd$};
			\node at (1.65,2.4) {\small $\ddd$};
			\node at (1.65,3.5) {\small $\bbb$};
			\node at (1.65,-0.45) {\small $\ccc$};
			\node at (1.65,0.4) {\small $\bbb$};
			\node at (1.65,-1.65) {\small $\ddd$};
			\node at (-1.1,-1.6) {\small $\bbb$};
			\node at (-1.7,-1.3) {\small $\bbb$};
			\node at (-1.7,1.5) {\small $\ddd$};
			\node[draw,shape=circle, inner sep=0.5] at (1,1) {\small $1$};
			\node[draw,shape=circle, inner sep=0.5] at (-1,1) {\small $7$};
			\node[draw,shape=circle, inner sep=0.5] at (1,-1) {\small $6$};
			\node[draw,shape=circle, inner sep=0.5] at (3,-1) {\small $5$};
			\node[draw,shape=circle, inner sep=0.5] at (3,1) {\small $4$};
			\node[draw,shape=circle, inner sep=0.5] at (3,3) {\small $3$};
			\node[draw,shape=circle, inner sep=0.5] at (1,3) {\small $2$};
			\node[draw,shape=circle, inner sep=0.5] at (-2,3) {\small $8$};
			\node at (1,-3) {\small $E_{56}=E_{78}=c,E_{23}=b$};
 	    \end{scope}  	 
       
	    \begin{scope}[xshift=10cm,yshift=-8.5cm]
	    	\draw (0,0) -- (0,2) 
	    	(0,0)--(-2,-2)
	    	(0,0)--(2,0)					
	    	(0,4)--(4,4)
	    	(2,2)--(2,4)
	    	(4,-2)--(4,2)
	    	(2,0)--(4,0);
	    	\draw[line width=1.5]
	    	(2,0)--(2,2)
	    	(-2,-2)--(4,-2)		
	    	(0,2)--(0,4)
	    	(4,2)--(4,4);
	    	\draw[dashed]
	    	(0,2)--(2,2)
	    	(2,0)--(2,-2)
	    	(2,2)--(4,2)			
	    	;
	    	\draw[dotted]
	    	(-2,2)--(0,2)
	    	(-2,-2)--(-2,2)
	    	(0,2)--(-1,3);
	    	\node at (2.3,0.4) {\small $\bbb$};
	    	\node at (2.3,1.55) {\small $\ddd$};
	    	\node at (2.3,2.4) {\small $\ccc$};
	    	\node at (2.3,3.7) {\small $\aaa$};
	    	\node at (2.3,-0.45) {\small $\ccc$};
	    	\node at (2.3,-1.65) {\small $\ddd$};
	    	\node at (3.65,0.35) {\small $\aaa$};
	    	\node at (3.65,1.6) {\small $\ccc$};
	    	\node at (3.65,2.4) {\small $\ddd$};
	    	\node at (3.65,3.55) {\small $\bbb$};
	    	\node at (3.65,-0.35) {\small $\aaa$};
	    	\node at (3.65,-1.65) {\small $\bbb$};
	    	\node at (0.35,0.3) {\small $\aaa$};
	    	\node at (0.35,1.65) {\small $\ccc$};
	    	\node at (0.35,2.4) {\small $\ddd$};
	    	\node at (0.35,3.55) {\small $\bbb$};
	    	\node at (0.15,-0.35) {\small $\aaa$};
	    	\node at (-0.3,0.15) {\small $\aaa$};
	    	\node at (1.65,1.6) {\small $\ddd$};
	    	\node at (1.65,2.4) {\small $\ccc$};
	    	\node at (1.65,3.65) {\small $\aaa$};
	    	\node at (1.65,-0.45) {\small $\ccc$};
	    	\node at (1.65,0.4) {\small $\bbb$};
	    	\node at (1.65,-1.65) {\small $\ddd$};
	    	\node at (-1.1,-1.6) {\small $\bbb$};
	    	\node[draw,shape=circle, inner sep=0.5] at (1,1) {\small $1$};
	    	\node[draw,shape=circle, inner sep=0.5] at (-1,1) {\small $7$};
	    	\node[draw,shape=circle, inner sep=0.5] at (1,-1) {\small $6$};
	    	\node[draw,shape=circle, inner sep=0.5] at (3,-1) {\small $5$};
	    	\node[draw,shape=circle, inner sep=0.5] at (3,1) {\small $4$};
	    	\node[draw,shape=circle, inner sep=0.5] at (3,3) {\small $3$};
	    	\node[draw,shape=circle, inner sep=0.5] at (1,3) {\small $2$};
	    	\node at (1,-3) {\small $E_{56}=c,E_{23}=a$};
	    \end{scope}  
    	   
	    \end{tikzpicture} 
		\caption{Partial neighborhoods of the 3rd special tile.}
		\label{ghj}
	\end{figure}
	
	\subsubsection*{Case  $E_{56}=b,E_{45}=a$}
	These edges determine $T_4,T_5,T_6$. 
	If $E_{23}=a$, we get $\ccc^2\ddd^2$ which contradicts Proposition \ref{proposition1} by $\aaa^3$ and $\bbb^4$. So $E_{23}=b$, which determines $\ddd_1\ddd_2\ddd_3\ddd_4$ and $T_2,T_3$. By $\aaa^3,\bbb^4,\ddd^4$ and Lemma \ref{anglesum}, we get $\aaa=\frac{2\pi}{3}$, $\bbb=\ddd=\frac{\pi}{2}$, $\ccc=(\frac13+\frac4f)\pi$.	
	By $\aaa_7$, we have the degree $5$ vertex $\thin\ccc_1\dash\ccc_2\thin\cdots=\aaa\bbb^2\ccc^2,\aaa\ccc^4$ or $\bbb\ccc^3\ddd$, which are all $>2\pi$, a contradiction. 
	
	\subsubsection*{Case  $E_{56}=b,E_{45}=c$}
	These edges determine $T_4,T_5,T_6$. By Lemma \ref{lem4}, we get  $\bbb_4\ddd_1\cdots=\bbb^2\ddd^2$, which determines $T_2,T_3$. By $\aaa_7$, we get $\thin\ccc_1\dash\ccc_2\thin\cdots=\aaa\bbb^2\ccc^2,\aaa\ccc^4$ or $\bbb\ccc^3\ddd$.
	By $\aaa^3,\bbb^2\ddd^2$ and Lemma \ref{anglesum}, we get $\aaa=\frac{2\pi}{3}$,  $\ccc=(\frac13+\frac4f)\pi$, $\bbb+\ddd=\pi$. So $\aaa+4\ccc>2\pi$ and $\bbb+\ddd+3\ccc>2\pi$. 
	Therefore $\ccc_1\ccc_2\cdots=\aaa\bbb^2\ccc^2=\thin\ccc_2\thick\ccc_1\thin\bbb_7\thick\bbb_8\thin\aaa\thin$. This determines $T_7,T_8$. By $\aaa\bbb^2\ccc^2$, we get $\bbb=(\frac{1}{3}-\frac{4}{f})\pi$, $\ddd=(\frac{2}{3}+\frac{4}{f})\pi$.
	Then Lemma $\ref{lembd}'$ implies that $\ddd\thick\ddd\cdots$ is not a vertex, contradicting $\ddd_7\thick\ddd_8\cdots$. 
	
	\subsubsection*{Case $E_{56}=c,E_{23}=E_{78}=b$}
	These edges determine $T_k, 2\le k\le7$. By the edge length consideration and $\bbb_1\bbb_4\ccc_5\ccc_6$, we get $\bbb_7\ccc_1\ccc_2\cdots=\bbb\ccc^3\ddd$, which determines $T_8,T_9$.  By $\aaa^3,\bbb^2\ccc^2,\ddd^4$ and $\bbb\ccc^3\ddd$, we get $\aaa=\frac{2\pi}{3},\bbb=\frac{3\pi}{4},\ccc=\frac{\pi}{4},\ddd=\frac{\pi}{2}$. So we have $\bbb_2\bbb_3\cdots=\bbb^2\ccc^2$ and $\aaa_2\aaa_9\cdots=\aaa^3$, which determine $T_{10}$. 
	Then $\bbb_9\thin\bbb_{10}\cdots$ is a vertex, contradicting Lemma \ref{lembd}.
	
	\subsubsection*{Case  $E_{56}=E_{78}=c,E_{23}=b$}
	These edges determine $T_k, 2\le k\le7$. By the edge length consideration, we get $\ccc_1\ccc_2\ccc_7\cdots=\aaa\ccc^4$ or $\bbb\ccc^3\ddd$,	which implies  $\bbb=\frac{2\pi}{3}$ or $\frac{3\pi}{4}$, $\ddd=\frac{\pi}{2}$ by $\aaa^3,\bbb^2\ccc^2,\ddd^4$.
    Then $\bbb_6\thin\bbb_7\cdots$ is a vertex, contradicting Lemma \ref{lembd}.
	
	\subsubsection*{Case  $E_{56}=c,E_{23}=a$}	
	These edges determine $T_k, 2\le k\le6$. By $\aaa^3,\bbb^2\ccc^2,\ccc^2\ddd^2$ and Lemma \ref{anglesum}, we get $\aaa=\frac{2\pi}{3}$, $\bbb=\ddd=(\frac13+\frac4f)\pi$, $\ccc=(\frac23-\frac4f)\pi$. By $\aaa_7$, we have the degree $5$ vertex $\ccc_1\ddd_2\cdots=\bbb^3\ccc\ddd,\bbb\ccc^3\ddd$ or $\bbb\ccc\ddd^3$, which are all $>2\pi$, a contradiction.
\end{proof}
\begin{proposition}
	There is no $a^2bc$-tiling with the 4th special tile in Fig. \ref{3444-5}.
\end{proposition}

\begin{proof}
	Let the fourth of Fig. \ref{3444-5} be the center tile $T_1$ in the partial neighborhoods in Fig. \ref{12345}. For the same reason as Proposition \ref{case1111}, we have $E_{78}=a$. There are three possibilities shown in Fig. \ref{12345}, according to $E_{89}=b$ or $c$ and $E_{29}=a$ or $b$.

	\begin{figure}[htp]
		\centering
		\begin{tikzpicture}[>=latex,scale=0.45]
			\draw
			(0,0)--(2,0)
			(0,0)--(0,2)
			(0,0)--(-2,-2)
			(0,2)--(0,4)		
			(0,4)--(2,4)
			(2,2)--(4,4)			
			(4,4)--(4,6)
			(0,4)--(-2,4)	
			(4,0)--(4,2)
			(2,2)--(4,2);
			\draw[dashed]
			(0,2)--(2,2)
			(-2,-2)--(-2,2)
			(-2,4)--(-2,2)
			(4,6)--(2,4)
			(2,0)--(4,0);
			\draw[dotted]
			(4,4)--(6,4)
			(4,2)--(6,4)
			(-2,-2)--(2,-2)
			(2,0)--(2,-2);
			\draw[line width=1.5]
			(2,0)--(2,2)
			(0,2)--(-2,2)
			(2,2)--(2,4);
			\node at (0.3,0.25) {\small $\aaa$};
			\node at (0.3,1.65) {\small $\ccc$};
			\node at (0.3,2.4) {\small $\ccc$};
			\node at (0.3,3.65) {\small $\aaa$};
			\node at (0.1,-0.3) {\small $\aaa$};
			\node at (-0.3,0.2) {\small $\aaa$};
			\node at (-0.3,1.5) {\small $\bbb$};
			\node at (-0.3,2.4) {\small $\bbb$};
			\node at (-0.3,3.65) {\small $\aaa$};
			\node at (-1.65,1.5) {\small $\ddd$};
			\node at (-1.65,-1.25) {\small $\ccc$};
			\node at (-1.65,2.4) {\small $\ddd$};
			\node at (-1.65,3.65) {\small $\ccc$};
			\node at (1.65,0.4) {\small $\bbb$};
			\node at (1.65,1.6) {\small $\ddd$};
			\node at (1.65,2.4) {\small $\ddd$};
			\node at (1.65,3.5) {\small $\bbb$};
			\node at (2.3,0.4) {\small $\ddd$};
			\node at (2.3,1.55) {\small $\bbb$};
			\node at (2.8,2.35) {\small $\aaa$};
			\node at (2.3,2.65) {\small $\bbb$};
			\node at (2.3,3.9) {\small $\ddd$};
			\node at (3.7,5.2) {\small $\ccc$};
			\node at (3.7,4.2) {\small $\aaa$};
			\node at (3.7,1.65) {\small $\aaa$};
			\node at (3.7,0.4) {\small $\ccc$};
			
			\node[draw,shape=circle, inner sep=0.5] at (1,1) {\small $1$};
			\node[draw,shape=circle, inner sep=0.5] at (3,4) {\small $3$};
			\node[draw,shape=circle, inner sep=0.5] at (4,3) {\small $4$};
			\node[draw,shape=circle, inner sep=0.5] at (3,1) {\small $5$};
			\node[draw,shape=circle, inner sep=0.5] at (1,3) {\small $2$};
			\node[draw,shape=circle, inner sep=0.5] at (-1,1) {\small $8$};
			\node[draw,shape=circle, inner sep=0.5] at (-1,3) {\small $9$};
			\node[draw,shape=circle, inner sep=0.5] at (1,-1) {\small $7$};
			\node at (1,-3) {\small $E_{89}=b$};
		\end{tikzpicture}\hspace{25pt}
		\begin{tikzpicture}[>=latex,scale=0.45]
			\draw
			(0,0)--(2,0)
			(0,0)--(0,2)
			(0,0)--(-2,-2)
			(0,2)--(0,4)
			(0,4)--(-2,4)												
			(0,4)--(2,4);
			\draw[dashed]
			(0,2)--(2,2)
			(0,2)--(-2,2);
			\draw[line width=1.5]
			(2,0)--(2,2)
			(-2,4)--(-2,2)
			(2,2)--(2,4)
			(-2,-2)--(-2,2);
			\draw[dotted]
			(2,0)--(2,-2)					
			(-2,-2)--(2,-2)
			(2,0)--(3,0)
			(2,2)--(3,3)
			(2,2)--(3,2);
			
			\node at (0.3,0.25) {\small $\aaa$};
			\node at (0.3,1.65) {\small $\ccc$};
			\node at (0.3,2.4) {\small $\ccc$};
			\node at (0.3,3.65) {\small $\aaa$};
			\node at (0.1,-0.3) {\small $\aaa$};
			\node at (-0.3,0.2) {\small $\aaa$};
			\node at (-0.3,1.65) {\small $\ccc$};
			\node at (-0.3,2.4) {\small $\ccc$};
			\node at (-0.3,3.65) {\small $\aaa$};
			\node at (-1.65,1.5) {\small $\ddd$};
			\node at (-1.65,-1.3) {\small $\bbb$};
			\node at (-1.65,2.4) {\small $\ddd$};
			\node at (-1.65,3.5) {\small $\bbb$};
			\node at (1.65,0.4) {\small $\bbb$};
			\node at (1.65,1.6) {\small $\ddd$};
			\node at (1.65,2.4) {\small $\ddd$};
			\node at (1.65,3.5) {\small $\bbb$};
			\node[draw,shape=circle, inner sep=0.5] at (1,1) {\small $1$};
			\node[draw,shape=circle, inner sep=0.5] at (1,3) {\small $2$};
			\node[draw,shape=circle, inner sep=0.5] at (-1,1) {\small $8$};
			\node[draw,shape=circle, inner sep=0.5] at (-1,3) {\small $9$};
			\node[draw,shape=circle, inner sep=0.5] at (1,-1) {\small $7$};
			\node at (0,-3) {\small $E_{89}=c,E_{29}=a$};
		\end{tikzpicture}\hspace{35pt}
		\begin{tikzpicture}[>=latex,scale=0.45]
			\draw
			(0,0)--(2,0)
			(0,0)--(0,2)
			(0,0)--(-2,-2)
			(2,2)--(2,4)
			(0,4)--(-2,4)
			(0,4)--(2,4)			
			(-2,4)--(-2,2)		
			(4,0)--(4,-2)
			(4,0)--(4,2)			
			(2,0)--(4,0);
			\draw[dashed]
			(0,2)--(2,2)				
			(2,2)--(4,2)			
			(0,2)--(-2,2);
			\draw[line width=1.5]
			(2,0)--(2,2)
			(-2,-2)--(-2,2)		
			(0,2)--(0,4)
			;
			\draw[dotted]
			(-2,-2)--(2,-2)
			(2,0)--(2,-2)
			(4,-2)--(2,-2)
			(2,2)--(3,3);
			
			\node at (0.3,0.25) {\small $\aaa$};
			\node at (0.3,1.65) {\small $\ccc$};
			\node at (0.3,2.4) {\small $\ddd$};
			\node at (0.3,3.5) {\small $\bbb$};
			\node at (0.1,-0.3) {\small $\aaa$};
			\node at (-0.3,0.2) {\small $\aaa$};
			\node at (-0.3,1.65) {\small $\ccc$};
			\node at (-0.3,2.4) {\small $\ddd$};
			\node at (-0.3,3.5) {\small $\bbb$};
			\node at (-1.65,1.5) {\small $\ddd$};
			\node at (-1.65,-1.3) {\small $\bbb$};
			\node at (-1.65,2.4) {\small $\ccc$};
			\node at (-1.65,3.65) {\small $\aaa$};
			\node at (1.65,0.4) {\small $\bbb$};
			\node at (1.65,1.6) {\small $\ddd$};
			\node at (1.65,2.4) {\small $\ccc$};
			\node at (1.65,3.65) {\small $\aaa$};
			\node at (2.3,0.4) {\small $\bbb$};
			\node at (2.3,-0.4) {\small $\ccc$};
			\node at (1.7,-0.4) {\small $\ccc$};
			\node at (2.3,1.6) {\small $\ddd$};
			\node at (3.7,1.6) {\small $\ccc$};
			\node at (3.7,0.3) {\small $\aaa$};
			\node[draw,shape=circle, inner sep=0.5] at (1,1) {\small $1$};
			\node[draw,shape=circle, inner sep=0.5] at (3,1) {\small $5$};
			\node[draw,shape=circle, inner sep=0.5] at (1,3) {\small $2$};
			\node[draw,shape=circle, inner sep=0.5] at (3,-1) {\small $6$};
			\node[draw,shape=circle, inner sep=0.5] at (-1,1) {\small $8$};
			\node[draw,shape=circle, inner sep=0.5] at (-1,3) {\small $9$};
			\node[draw,shape=circle, inner sep=0.5] at (1,-1) {\small $7$};
			\node at (1,-3) {\small $E_{89}=c,E_{29}=b$};
		\end{tikzpicture}
		\caption{Partial neighborhoods of the 4th special tile.}
		\label{12345}
	\end{figure}
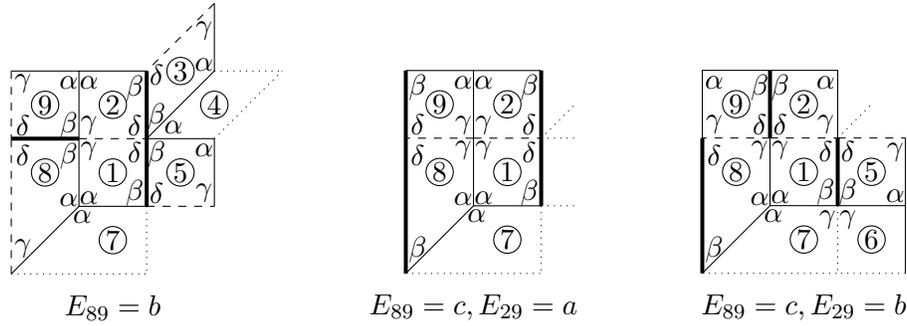
	
	\subsubsection*{Case $E_{89}=b$}
	This edge determines $T_8$. By Lemma \ref{lem4}, $\bbb_8\ccc_1\cdots=\bbb^2\ccc^2$ determines $T_2,T_9$. By $E_{23}=b$, we get $E_{34}=a$ or $c$, which implies $\thick\ddd_1\dash\ddd_2\thick\cdots=\bbb\ccc\ddd^3$ or $\aaa\bbb^2\ddd^2$. But $\bbb\ccc\ddd^3$, $\aaa^3$ and $\bbb^2\ccc^2$ imply $3(\aaa+\bbb+\ccc+\ddd)=6\pi$, contradicting Lemma \ref{anglesum}. So $\ddd_1\ddd_2\cdots=\aaa\bbb^2\ddd^2$, which determines $T_3,T_5$. By Lemma \ref{lem4}, $\bbb_1\ddd_5\cdots=\bbb^2\ddd^2$, contradicting $\aaa\bbb^2\ddd^2$. 
	
	\subsubsection*{Case  $E_{89}=c,E_{29}=a$}
	These two edges determine $T_2,T_8,T_9$. By Lemma \ref{lem4}, $\bbb_1\cdots=\bbb^4,\bbb^2\ccc^2$ or $\bbb^2\ddd^2$. If $\bbb^4$ or $\bbb^2\ccc^2$ is a vertex, we get $\bbb=\ccc$ by $\ccc^4$, contradicting Lemma \ref{lemd4}. 
	So $\aaa^3,\ccc^4,\bbb^2\ddd^2$ are vertices, contradicting Proposition \ref{proposition1}$'$.
	\subsubsection*{Case  $E_{89}=c,E_{29}=b$}
	These two edges determine $T_2,T_8,T_9$. By Lemma \ref{lem4}, $\bbb_1\cdots=\bbb^4,\bbb^2\ccc^2$ or $\bbb^2\ddd^2$. But $\aaa^3,\ccc^2\ddd^2,\bbb^4$ contradicts Proposition \ref{proposition1}. If $\bbb^2\ddd^2$ is a vertex, we get $\bbb=\ccc$ by $\ccc^2\ddd^2$, contradicting Lemma \ref{lemd4}. So $\bbb_1\cdots=\bbb^2\ccc^2$, which determines $T_5$. Then $\thin\ccc_2\dash\ddd_1\thick\ddd_5\dash\cdots=\aaa\ccc^2\ddd^2$ or $\bbb\ccc\ddd^3$.  By $\aaa^3,\bbb^2\ccc^2$ and $\ccc^2\ddd^2$, we get $\aaa=\frac{2\pi}{3}$, $\bbb=\ddd=(\frac13+\frac4f)\pi$, $\ccc=(\frac23-\frac4f)\pi$. So $\aaa+2\ccc+2\ddd>2\pi$ and $\bbb+\ccc+3\ddd>2\pi$, a contradiction.
\end{proof}
\begin{proposition}
	There is no $a^2bc$-tiling with the 5th special tile in Fig. \ref{3444-5}.
\end{proposition}

\begin{proof}
	Let the fifth of Fig. \ref{3444-5} be the center tile $T_1$ in the partial neighborhoods in Fig. \ref{3345 case3.1}. By Lemma \ref{lem3}, we get  $\bbb_1\cdots=\aaa_5\bbb_1\bbb_4$, which determines $T_4$. By $\aaa_5\bbb_1\bbb_4$ and Lemma \ref{lem4}, $\aaa_1\cdots=\aaa^2\ccc^2=\aaa_1\aaa_7\ccc_5\ccc_6$, which determines $T_5,T_6$. There are two cases shown in Fig. \ref{3345 case3.1}, according to $E_{23}=a$ or $b$.
	
	\begin{figure}[htp]
		\centering
		\begin{tikzpicture}[>=latex,scale=0.45]
			\draw (0,0) -- (0,2) 
			(0,0) -- (2,0)
			(2,0)--(4,-2)
			(4,-2)--(4,2)
			(2,2)--(2,4)
			(2,4)--(4,4)
			(0,4)--(2,4)
			(0,0)--(-2,0)
			(0,2)--(-2,4)
			(-2,0)--(-2,-2)
			(-4,4)--(-2,4)
			(-2,4)--(-2,6);
			\draw[dashed]  (0,2)--(2,2)
			(2,2)--(4,2)
			(0,0)--(0,-2)
			(0,2)--(-2,2)
			(0,4)--(-2,6);
			\draw[line width=1.5] (2,0)--(2,2)
			(4,2)--(4,4)
			(0,2)--(0,4)
			(0,-2)--(4,-2)
			(0,-2)--(-2,-2)
			(-2,0)--(-2,2)
			(-2,2)--(-4,4);
			\node[draw,shape=circle, inner sep=0.5] at (1,1) {\small $1$};
			\node[draw,shape=circle, inner sep=0.5] at (1,3) {\small $2$};
			\node[draw,shape=circle, inner sep=0.5] at (3,3) {\small $3$};
			\node[draw,shape=circle, inner sep=0.5] at (3,1) {\small $4$};
			\node[draw,shape=circle, inner sep=0.5] at (-1,1) {\small $7$};
			\node[draw,shape=circle, inner sep=0.5] at (-1,4) {\small $9$};
			\node[draw,shape=circle, inner sep=0.5] at (-2.2,3.1) {\small $8$};
			\node[draw,shape=circle, inner sep=0.5] at (1,-1) {\small $5$};
			\node[draw,shape=circle, inner sep=0.5] at (-1,-1) {\small $6$};
			\node at (-1.65,0.35){\small $\bbb$};
			\node at (-1.65,5.15){\small $\ccc$};
			\node at (-2.8,3.5){\small $\bbb$};		
			\node at (0.35,0.35){\small $\aaa$};
			\node at (1.65,0.35){\small $\bbb$};
			\node at (1.65,1.65){\small $\ddd$}; \node at (2.35,1.65){\small $\ddd$};
			\node at (0.35,1.65){\small $\ccc$};
			\node at (2.35,0.3){\small $\bbb$}; \node at (1.8,-0.35){\small $\aaa$};
			\node at (3.65,1.65){\small $\ccc$}; \node at (3.7,-1.2){\small $\aaa$};
			\node at (0.35,-0.35){\small $\ccc$}; \node at (0.35,-1.6){\small $\ddd$};
			\node at (-0.35,-1.6){\small $\ddd$}; \node at (-0.35,-0.35){\small $\ccc$};
			\node at (-0.35,0.35){\small $\aaa$}; \node at (2.9,-1.55){\small $\bbb$};
			\node at (-1.65,-1.55){\small $\bbb$}; \node at (-1.65,-0.35){\small $\aaa$};
			\node at (0.45,2.35){\small $\ddd$}; 
			\node at (-0.35,2.95){\small $\bbb$}; 
			\node at (-0.95,2.25){\small $\ccc$}; 
			\node at (-0.35,1.65){\small $\ccc$};  
			\node at (1.65,2.25){\small $\ccc$};
			\node at (2.25,2.25){\small $\ccc$};
			\node at (3.65,2.35){\small $\ddd$};
			\node at (3.65,3.45){\small $\bbb$};
			\node at (2.35,3.65){\small $\aaa$};
			\node at (1.65,3.65){\small $\aaa$};
			\node at (0.45,3.55){\small $\bbb$};
			\node at (-0.35,3.8){\small $\ddd$};
			\node at (-1.65,4.25){\small $\aaa$};
			\node at (-2.15,3.75){\small $\aaa$};
			\node at (-1.65,1.55){\small $\ddd$};
			\node at (-1.85,2.45){\small $\ddd$};
			\node at (1,-4) {\small $E_{23}=a$};
		\end{tikzpicture}\hspace{70pt}
		\begin{tikzpicture}[>=latex,scale=0.45]  
			\draw
			(0,0)--(2,0)
			(0,0)--(0,2)
			(0,0)--(-2,0)
			(4,2)--(4,-2)
			(4,-2)--(2,0)
			(-2,0)--(-2,-2)
			(0,2)--(0,4)
			(0,4)--(4,4)
			(4,2)--(4,4);
			\draw[dashed]
			(0,2)--(4,2)
			(0,0)--(0,-2);
			\draw[line width=1.5]
			(2,0)--(2,2)
			(-2,-2)--(4,-2)
			(2,2)--(2,4);
			\draw[dotted]
			(0,2)--(-2,2)
			(-2,0)--(-2,2)
			(0,2)--(-2,4)
			(-2,2)--(-4,4)
			(-4,4)--(-2,2)
			(-2,4)--(-2,6)
			(-2,6)--(0,4)
			(-4,4)--(-2,4);
			\node[draw,shape=circle, inner sep=0.5] at (1,1) {\small $1$};
			\node[draw,shape=circle, inner sep=0.5] at (1,3) {\small $2$};
			\node[draw,shape=circle, inner sep=0.5] at (3,3) {\small $3$};
			\node[draw,shape=circle, inner sep=0.5] at (3,1) {\small $4$};\node[draw,shape=circle, inner sep=0.5] at (1,-1) {\small $5$};
			\node[draw,shape=circle, inner sep=0.5] at (-1,-1) {\small $6$};
			\node[draw,shape=circle, inner sep=0.5] at (-1,1) {\small $7$};
			\node[draw,shape=circle, inner sep=0.5] at (-1,4) {\small $9$};
			\node[draw,shape=circle, inner sep=0.5] at (-2.2,3.1) {\small $8$};

			\node at (0.35,0.35){\small $\aaa$};
			\node at (1.65,0.35){\small $\bbb$};
			\node at (1.65,1.65){\small $\ddd$}; \node at (2.35,1.65){\small $\ddd$};
			\node at (0.35,1.65){\small $\ccc$};
			
			\node at (2.35,0.3){\small $\bbb$}; \node at (1.8,-0.35){\small $\aaa$};
			
			\node at (3.65,1.65){\small $\ccc$}; \node at (3.7,-1.2){\small $\aaa$};
			
			\node at (0.35,-0.35){\small $\ccc$}; \node at (0.35,-1.6){\small $\ddd$};
			\node at (-0.35,-1.6){\small $\ddd$}; \node at (-0.35,-0.35){\small $\ccc$};
			\node at (-0.35,0.35){\small $\aaa$}; \node at (2.9,-1.55){\small $\bbb$};
			\node at (-1.65,-1.55){\small $\bbb$}; \node at (-1.65,-0.35){\small $\aaa$};
			\node at (1.65,2.35){\small $\ddd$};
			\node at (2.35,2.35){\small $\ddd$};
			\node at (3.65,2.25){\small $\ccc$};
			\node at (3.65,3.45){\small $\aaa$};
			\node at (2.35,3.55){\small $\bbb$};
			\node at (1.65,3.55){\small $\bbb$};
			\node at (0.35,2.25){\small $\ccc$};
			\node at (0.35,3.65){\small $\aaa$};

			\node at (4.25,-2.25){\small $\bbb$};
			\node at (1.6,-4) {\small $E_{23}=b$};
			\node at (0.3,-2.5) {\small $\ddd$};
			\node at (-0.3,-2.5) {\small $\ddd$};
		\end{tikzpicture}
		\caption{Partial neighborhoods of the 5th special tile.} \label{3345 case3.1}
	\end{figure}
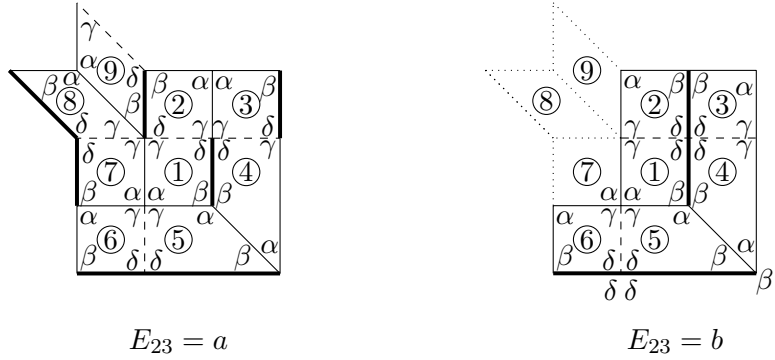
	\subsubsection*{Case $E_{23}=a$}	
	This edge determines $T_2,T_3$. By $\aaa\bbb^2$, $\aaa^2\ccc^2$, $\ccc^2\ddd^2$ and Lemma \ref{anglesum}, we get $\aaa=\ddd=\frac{8\pi}{f}$, $\bbb=(1-\frac4f)\pi$, $\ccc=(1-\frac8f)\pi$, $f\ge10$. 
	By $\aaa_7$, $\ccc_1\ddd_2\cdots=\bbb^3\ccc\ddd,\bbb\ccc^3\ddd$ or $\bbb\ccc\ddd^3$. 
	But $3\bbb+\ccc+\ddd>2\pi$, $\bbb+\ccc+3\ddd>2\pi$, so $\ccc_1\ddd_2\cdots=\bbb\ccc^3\ddd$, which determines $T_7,T_8,T_9$ and $f=10$, 
	$\aaa=\ddd=\frac{4}{5}\pi$, $\bbb=\frac{3}{5}\pi$, $\ccc=\frac{1}{5}\pi$.
	Then $\thick\ddd_5\dash\ddd_6\thick\cdots$ is a vertex, contradicting Lemma \ref{lembd}.

	\subsubsection*{Case  $E_{23}=b$}
	
	This edge determines $T_2,T_3$. By $\aaa_7$ and $\aaa\bbb^2$, we get $\ccc_1\ccc_2\cdots=\aaa\ccc^4$ or $\bbb\ccc^3\ddd$. 
	If $\ccc_1\ccc_2\cdots=\aaa\ccc^4$, the angle sums at $\aaa\bbb^2,\aaa^2\ccc^2,\ddd^4,\aaa\ccc^4$ imply $\aaa=\bbb=\frac{2\pi}{3}$, $\ccc=\frac{\pi}{3}$, $\ddd=\frac{\pi}{2}$, $f=24$.
	If $\ccc_1\ccc_2\cdots=\bbb\ccc^3\ddd$, we get $\aaa=\frac{5}{7}\pi$, $\bbb=\frac{9}{14}\pi$, $\ccc=\frac{2}{7}\pi$, $\ddd=\frac{\pi}{2}$, $f=28$. Both implies $\ddd_5\ddd_6\cdots=\ddd^4$ and $\aaa_4\bbb_5\cdots=\aaa\bbb^2$. So $T_5$ is a $3344$-Tile, which has been handled in Section \ref{3345d}.	
\end{proof}

\begin{proposition}
	There is no $a^2bc$-tiling with the 6th special tile in Fig. \ref{3444-5}.
\end{proposition}

\begin{proof}
	Let the sixth of Fig. \ref{3444-5} be the center tile $T_1$ in the partial neighborhoods in Fig. \ref{3445case5.1}. By Lemma \ref{lem3}, we get  $\bbb_1\cdots=\aaa_5\bbb_1\bbb_4$, which determines $T_4$.  There are two possibilities shown in Fig. \ref{3445case5.1}, according to $E_{23}=a$ or $b$. 
	\begin{figure}[htp]
		\centering
		\begin{tikzpicture}[>=latex,scale=0.45]
			\draw (0,0) -- (0,2) 
			(0,0) -- (2,0)
			(2,0)--(4,-2)
			(4,-2)--(4,2)
			(2,2)--(2,4)
			(2,4)--(4,4)
			(0,4)--(2,4)
			(-2,2)--(-2,4)
			(-2,4)--(0,4)
			(0,0)--(-2,0)
			(0,0)--(-2,-2)
			(-2,-2)--(-2,-4);
			\draw[dashed]  (0,2)--(2,2)
			(2,2)--(4,2)
			(0,2)--(-2,2)
			(0,0)--(0,-2);
			\draw[line width=1.5] (2,0)--(2,2)
			(4,2)--(4,4)
			(0,2)--(0,4)
			(-2,0)--(-2,2)
			(0,-2)--(4,-2)
			(0,-2)--(-2,-4);
			\draw[dotted] (-2,0)--(-4,-2)
			(-4,-2)--(-2,-2);

			\node[draw,shape=circle, inner sep=0.5] at (1,1) {\small $1$};
			\node[draw,shape=circle, inner sep=0.5] at (1,3) {\small $2$};
			\node[draw,shape=circle, inner sep=0.5] at (3,3) {\small $3$};
			\node[draw,shape=circle, inner sep=0.5] at (3,1) {\small $4$};
			\node[draw,shape=circle, inner sep=0.5] at (-1,3) {\small $9$};
			\node[draw,shape=circle, inner sep=0.5] at (-1,1) {\small $8$};
			\node[draw,shape=circle, inner sep=0.5] at (-2,-1) {\small $7$};
			\node[draw,shape=circle, inner sep=0.5] at (1,-1) {\small $5$};
			\node[draw,shape=circle, inner sep=0.5] at (-1,-2) {\small $6$};
			%
			\node at (0.35,0.35){\small $\aaa$};
			\node at (-0.35,0.35){\small $\aaa$};
			\node at (-0.85,-0.35){\small $\aaa$};
			\node at (-0.35,-0.85){\small $\ccc$};
			\node at (0.35,-0.45){\small $\ccc$};
			
			\node at (1.65,0.35){\small $\bbb$};
			\node at (2.35,0.35){\small $\bbb$};
			\node at (2,-0.35){\small $\aaa$};
			
			\node at (1.65,1.65){\small $\ddd$};
			\node at (2.35,1.65){\small $\ddd$};
			\node at (1.65,2.35){\small $\ccc$};
			\node at (2.35,2.35){\small $\ccc$};
			
			\node at (0.35,1.65){\small $\ccc$};   
			\node at (-0.35,1.65){\small $\ccc$}; 
			\node at (0.35,2.35){\small $\ddd$}; 
			\node at (-0.35,2.35){\small $\ddd$};

			
			\node at (3.65,2.35){\small $\ddd$};
			\node at (3.65,3.45){\small $\bbb$};
			\node at (2.35,3.55){\small $\aaa$};
			\node at (1.65,3.55){\small $\aaa$};
			\node at (0.35,3.55){\small $\bbb$};
			\node at (0.35,-1.6){\small $\ddd$};
			\node at (-0.35,-1.9){\small $\ddd$};
			\node at (-0.35,-1.9){\small $\ddd$};
			
			\node at (-1.65,-2.25){\small $\aaa$};	
			\node at (-1.65,-3.25){\small $\bbb$};	
			\node at (3.65,1.65){\small $\ccc$};	
			\node at (-1.65,0.35){\small $\bbb$};	
			\node at (-1.65,1.65){\small $\ddd$};
			\node at (-1.65,2.25){\small $\ccc$};
			\node at (-1.65,3.65){\small $\aaa$};
			\node at (-0.35,3.55){\small $\bbb$};
			\node at (3.7,-1.2){\small $\aaa$};
			\node at (2.9,-1.55){\small $\bbb$};
			\node at (1,-4) {\small $E_{23}=a$};
		\end{tikzpicture}\hspace{70pt}
		\begin{tikzpicture}[>=latex,scale=0.45]
			\draw (0,0) -- (0,2) 
			(0,0) -- (2,0)
			(2,0)--(4,-2)
			(4,-2)--(4,2)
			(2,4)--(4,4)
			(4,2)--(4,4)
			(0,4)--(2,4)
			(0,2)--(0,4)
			(0,0)--(-2,-2)
			(-2,-2)--(-2,-4);
			\draw[dashed]  (0,2)--(2,2)
			(2,2)--(4,2)
			(0,0)--(0,-2);
			\draw[line width=1.5] (2,0)--(2,2)
			(2,2)--(2,4)
			(0,-2)--(4,-2)
			(-2,-4)--(0,-2);
			\draw[dotted] (0,0)--(-2,0)
			(-2,0)--(-2,2)
			(-2,2)--(0,2)
			(-2,0)--(-4,-2)
			(-4,-2)--(-2,-2)
			(-2,2)--(-2,4)
			(-2,4)--(0,4);
			
			\node[draw,shape=circle, inner sep=0.5] at (1,1) {\small $1$};
			\node[draw,shape=circle, inner sep=0.5] at (1,3) {\small $2$};
			\node[draw,shape=circle, inner sep=0.5] at (3,3) {\small $3$};
			\node[draw,shape=circle, inner sep=0.5] at (3,1) {\small $4$};
			\node[draw,shape=circle, inner sep=0.5] at (-1,3) {\small $9$};
			\node[draw,shape=circle, inner sep=0.5] at (-1,1) {\small $8$};
			\node[draw,shape=circle, inner sep=0.5] at (-2,-1) {\small $7$};
			\node[draw,shape=circle, inner sep=0.5] at (1,-1) {\small $5$};
			\node[draw,shape=circle, inner sep=0.5] at (-1,-2) {\small $6$};
			
			%
			\node at (0.35,0.35){\small $\aaa$};
			\node at (0.35,-0.35){\small $\ccc$};
			\node at (-0.25,-0.75){\small $\ccc$};
			
			\node at (1.65,0.35){\small $\bbb$};
			\node at (2.35,0.35){\small $\bbb$};
			\node at (2,-0.35){\small $\aaa$};
			
			\node at (1.65,1.65){\small $\ddd$};
			\node at (2.35,1.65){\small $\ddd$};
			\node at (1.65,2.35){\small $\ddd$};
			\node at (2.35,2.35){\small $\ddd$};
			
			\node at (0.35,1.65){\small $\ccc$};   
			\node at (0.35,2.35){\small $\ccc$}; 
			
			\node at (3.65,2.35){\small $\ccc$};
			\node at (3.65,3.55){\small $\aaa$};
			\node at (2.35,3.55){\small $\bbb$};
			\node at (1.65,3.55){\small $\bbb$};
			\node at (0.35,3.55){\small $\aaa$};
			\node at (0.35,-1.6){\small $\ddd$};
			\node at (-0.35,-1.9){\small $\ddd$};
			\node at (-0.35,-1.9){\small $\ddd$};
			\node at (3.65,1.65){\small $\ccc$};
			\node at (-1.65,-3.25){\small $\bbb$};
			\node at (-1.65,-2.25){\small $\aaa$};
			\node at (3.7,-1.2){\small $\aaa$};
			\node at (2.9,-1.55){\small $\bbb$};
			\node at (4.2,-2.25){\small $\bbb$};
			\node at (-0.2,-2.8) {\small $\ddd$};
			\node at (0.4,-2.5) {\small $\ddd$};
			\node at (1.8,-4) {\small $E_{23}=b$};
		\end{tikzpicture}     
		\caption{Partial neighborhoods of the 6th special tile}\label{3445case5.1}
	\end{figure}
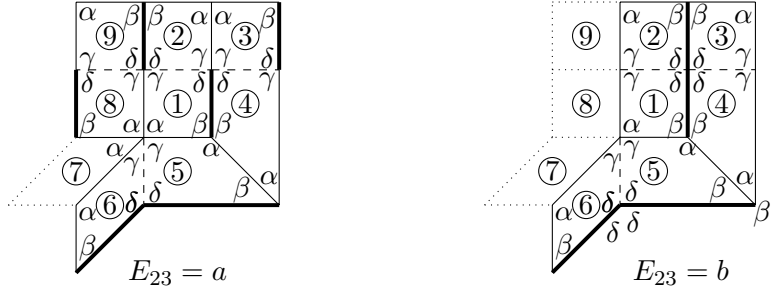
	\subsubsection*{Case $E_{23}=a$}		
	This edge determines $T_2,T_3$. By Lemma \ref{lem4}, $\ccc_1\ddd_2\cdots=\ccc^2\ddd^2$ determines $T_8,T_9$. 
	By $\aaa_5\bbb_1\bbb_4$ and Lemma \ref{anglesum}, we get $\aaa_1\aaa_8\cdots=\aaa^3\ccc^2$, which determines $T_5,T_6$.
	Then $\aaa\bbb^2,\ccc^2\ddd^2,\aaa^3\ccc^2$ imply $\aaa=\frac{8\pi}{f}$, $\bbb=(1-\frac{4}{f})\pi$, $\ccc=(1-\frac{12}{f})\pi$, $\ddd=\frac{12\pi}{f}$. Then Lemma \ref{lembd}, $\thick\ddd_5\dash\ddd_6\thick\cdots$ and $\bbb+\ddd>\pi$ imply $\ddd\le\frac{\pi}{2}$. So we have $f\ge24$ and $\aaa\le\frac{\pi}{3},\bbb\ge\frac{5\pi}{6},\ccc\ge\frac{\pi}{2}$. 
	By Lemma \ref{lembd}, there is no $\bbb\thin\bbb\cdots$. So we get $\bbb_8\thin\ccc_7\cdots=\theta\thick\bbb_8\thin\ccc_7\dash\rho\cdots$, where $\theta=\bbb$ or $\ddd$ and $\rho=\ccc$ or $\ddd$. But we always have $\theta+\bbb>\pi$ and $\rho+\ccc\ge\pi$, a contradiction.

	\subsubsection*{Case $E_{23}=b$}		
	This edge determines $T_2,T_3$. By $\aaa\bbb^2$ and Lemma $\ref{lemac2}'$, we get $E_{56}=c$, which determines $T_5$. By Lemma $\ref{lemac2}'$, we get $\aaa_1\ccc_5\cdots=\aaa\ccc^4$,  or $\aaa^3\ccc^2$. 
	
	If $\aaa_1\ccc_5\cdots=\aaa\ccc^4$, then $E_{89}=a$ and $\ccc_1\ccc_2\cdots=\aaa^2\ccc^2$. The angle sums at  $\aaa\ccc^4,\aaa^2\ccc^2,\ddd^4,\aaa\bbb^2$ imply $\aaa=\bbb=\frac{2\pi}{3}$, $\ccc=\frac{\pi}{3}$, $\ddd=\frac{\pi}{2}$, $f=24$. 
	
	If $\aaa_1\ccc_5\cdots=\aaa^3\ccc^2$, then $\aaa=(\frac12-\frac4f)\pi$, $\bbb=(\frac34+\frac2f)\pi$, $\ccc=(\frac14+\frac6f)\pi$, $\ddd=\frac12$ by $\aaa\bbb^2$ and $\ddd^4$. Since $\bbb+\ccc>\pi$, we get $E_{89}=c$ and $\ccc_1\ccc_2\cdots=\ccc^4$.  Then $\aaa=\frac{\pi}{3}$, $\bbb=\frac{5\pi}{6}$, $\ccc=\ddd=\frac{\pi}{2}$ and $f=24$. 
	
	Both imply $\thick\ddd_5\dash\ddd_6\thick\cdots=\ddd^4$  and $\aaa_4\bbb_5\cdots=\aaa\bbb^2$. So $T_5$ is a $3344$-Tile, which has been handled in Section \ref{3345d}. 
\end{proof}
\begin{proposition}
	There is no $a^2bc$-tiling with the 7th special tile in Fig. \ref{3444-5}.
\end{proposition}

\begin{proof}

	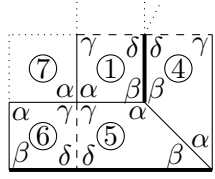
\begin{figure}[htp]
		\centering
		\begin{tikzpicture}[>=latex,scale=0.45]
			\draw (0,0) -- (0,2) 
			(0,0) -- (2,0)
			(2,0)--(4,-2)
			(4,-2)--(4,2)
			(0,0)--(-2,0)
			(-2,0)--(-2,-2);
			\draw[dashed]  (0,2)--(2,2)
			(2,2)--(4,2)
			(0,0)--(0,-2);
			\draw[line width=1.5] (2,0)--(2,2)
			(0,-2)--(4,-2)
			(-2,-2)--(0,-2);
			\draw[dotted] (0,2)--(-2,2)
			(-2,2)--(-2,0)
			(2,2)--(2,3)
		    (2,2)--(2.5,3)
			(0,2)--(0,3);
			
			\node at (0.35,0.35){\small $\aaa$};
			\node at (-0.35,0.35){\small $\aaa$};
			\node at (0.35,-0.35){\small $\ccc$};
			\node at (-0.35,-0.35){\small $\ccc$};
			
			\node at (1.65,0.35){\small $\bbb$};
			\node at (2.35,0.35){\small $\bbb$};
			\node at (1.8,-0.3){\small $\aaa$};
			
			\node at (1.65,1.65){\small $\ddd$};
			\node at (2.35,1.65){\small $\ddd$};
			
			\node at (0.35,1.65){\small $\ccc$};
			
			\node at (3.65,1.65){\small $\ccc$};
			\node at (2.9,-1.6){\small $\bbb$};
			\node at (3.65,-1.15){\small $\aaa$};
			\node at (0.35,-1.55){\small $\ddd$};
			\node at (-0.35,-1.55){\small $\ddd$};
			\node at (-1.65,-1.55){\small $\bbb$};
			\node at (-1.65,-0.35){\small $\aaa$};
			
			\node[draw,shape=circle, inner sep=0.5] at (1,1) {\small $1$};
			\node[draw,shape=circle, inner sep=0.5] at (3,1) {\small $4$};
			\node[draw,shape=circle, inner sep=0.5] at (-1,1) {\small $7$};
			\node[draw,shape=circle, inner sep=0.5] at (1,-1) {\small $5$};
			\node[draw,shape=circle, inner sep=0.5] at (-1,-1) {\small $6$};

		\end{tikzpicture}  
		\caption{Partial neighborhood of the 7th special tile.} \label{3445subcase7}
	\end{figure}
	
	Let the seventh of Fig. \ref{3444-5} be the center tile $T_1$ in the partial neighborhoods in Fig. \ref{3445subcase7}. By Lemma \ref{lem3}, $\bbb_1\cdots=\aaa_5\bbb_1\bbb_4$ determines $T_4$.
	By $\aaa_5\bbb_1\bbb_4$ and Lemma \ref{lem4}, $\aaa_1\cdots=\aaa^2\ccc^2$ which determines $T_5,T_6$. By the edge length consideration and Lemma $\ref{lemac2}'$, we have  $\dash\ddd_1\thick\ddd_4\dash\cdots=\bbb\ccc\ddd^3$. Then by $\aaa\bbb^2$ and $\aaa^2\ccc^2$, we get  $3(\aaa+\bbb+\ccc+\ddd)=(\aaa+2\bbb)+(2\aaa+2\ccc)+(\bbb+\ccc+3\ddd)=6\pi$, contradicting Lemma \ref{anglesum}. 
\end{proof}

\vspace{2ex}
	
In summary, except quadrilateral subdivisions of the octahedron (every tile is $3444$-Tile), all other $a^2bc$-tilings with a $3444$-Tile would also have a $334d$-Tile and can be picked out from Section \ref{3345d} as follows: 
\begin{itemize}
	\item the flip of a special quadrilateral subdivision of the octahedron with $f=24$:
	$T(2\aaa^3,6\aaa\ccc^2,6\ddd^4,6\bbb^2\ccc^2,6\aaa^2\bbb^2)$ ; 
	\item the $3$-layer earth map tiling with $f=16$: $T(8\aaa\ccc^2,2\bbb^4,4\ddd^4,4\aaa^2\bbb^2)$ ;
	\item the second flip of the $3$-layer earth map tiling  with $f=24$:
	
	$T(10\aaa\ccc^2,2\aaa\bbb^4,2\bbb^2\ccc^2,6\ddd^4,6\aaa^2\bbb^2)$.
\end{itemize}

All $a^2bc$-tilings with a $3445$-Tile would also have a $334d$-Tile and can be picked out from Section \ref{3345d} as follows: 
\begin{itemize}
	\item the first flip of the $3$-layer earth map tiling with $f=24$:
	
	 $T(12\aaa\ccc^2,4\aaa\bbb^4,6\ddd^4,4\aaa^2\bbb^2)$;
	\item the second flip of the $3$-layer earth map tiling with $f=24$:
	
	 $T(10\aaa\ccc^2,2\aaa\bbb^4,2\bbb^2\ccc^2,6\ddd^4,6\aaa^2\bbb^2)$.
\end{itemize}

\vspace{9pt}
The classification for Type $a^2bc$ has been completed and three classes of tilings are summarized in the main theorem in the introduction.

\end{document}